# Efficient Alternative Finite Difference WENO Schemes for Hyperbolic Systems with Non-Conservative Products

By

Dinshaw S. Balsara[1,2], Deepak Bhoriya[1], Chi-Wang Shu[3] and Harish Kumar[4]

[1]Physics Department, [2]ACMS Department, University of Notre Dame

[3]Division of Applied Mathematics, Brown University

[4]Department of Mathematics, Indian Institute of Technology, Delhi

**Abstract**

Higher order finite difference Weighted Essentially Non-Oscillatory (WENO) schemes for conservation laws represent a technology that has been reasonably consolidated. They are extremely popular because, when applied to multidimensional problems, they offer high order accuracy at a fraction of the cost of finite volume WENO or DG schemes. They come in two flavors. There is the classical finite difference WENO (FD-WENO) method (Shu and Osher, J. Comput. Phys., 83 (1989) 32-78). However, in recent years there is also an alternative finite difference WENO (AFD-WENO) method which has recently been formalized into a very useful general-purpose algorithm for conservation laws (Balsara *et al*., *Efficient Alternative Finite Difference WENO Schemes for Hyperbolic Conservation Laws*, submitted to CAMC (2023)). However, the FD-WENO algorithm has only very recently been formulated for hyperbolic systems with non-conservative products (Balsara *et al*., *Efficient Finite Difference WENO Scheme for Hyperbolic Systems with Non-Conservative Products*, to appear CAMC (2023)). In this paper we show that there are substantial advantages in obtaining an AFD-WENO algorithm for hyperbolic systems with non-conservative products. Such an algorithm is documented in this paper. We present an AFD-WENO formulation in fluctuation form that is carefully engineered to retrieve the flux form when that is warranted and nevertheless extends to non-conservative products. The method is flexible because it allows any Riemann solver to be used. The formulation we arrive at is such that when non-conservative products are absent it reverts exactly to the formulation in the second citation above which is in exact flux conservation form. The ability to transition to a precise conservation form when non-conservative products are absent ensures, via the Lax-Wendroff theorem, that shock locations will be exactly captured by the method.

We present two formulations of AFD-WENO that can be used with hyperbolic systems with non-conservative products and stiff source terms with slightly differing computational complexities. The speeds of our new AFD-WENO schemes are compared to the speed of the classical FD-WENO algorithm from the first of the above-cited papers. At all orders, AFD-WENO outperforms FD-WENO. We also show a very desirable result that higher order variants of AFD-WENO schemes don't cost that much more than their lower order variants. This is because the larger number of floating point operations associated with larger stencils are almost very efficiently amortized by the CPU when the AFD-WENO code is designed to be cache friendly. This should have great, and very beneficial, implications for the role of our AFD-WENO schemes in Peta- and Exascale computing.

We apply the method to several stringent test problems drawn from the Baer-Nunziato system, two layer shallow water equations and multicomponent debris flow. The method meets its design accuracy for smooth flow and can handle stringent problems in one and multiple dimensions. Because of the pointwise nature of its update, AFD-WENO for hyperbolic systems with non-conservative products is also shown to be a very efficient performer on problems with stiff source terms.

**I) Introduction**

In a landmark paper, Harten et al. [27] developed finite volume Essentially Non-Oscillatory (ENO) schemes which could evolve conservation laws with better than second order of accuracy. Soon thereafter, Shu and Osher [40], [41] developed finite difference versions of ENO schemes which were substantially faster than their finite volume counterparts in multi-dimensions. Early versions of ENO schemes suffered from the deficiency that rapid switching of the stencils could diminish the order of accuracy. Weighted Essentially Non-Oscillatory (WENO) schemes overcame this deficiency by making a non-linearly hybridized weighting of the reconstruction polynomials from multiple stencils (Liu, Osher & Chan [33], Jiang & Shu [29]). The methods were extended to seventh, ninth and eleventh orders by Balsara & Shu [4] and much later to seventeenth order by Gerolymos, Sénéchal & Vallet [26]. Some of the early deficiencies of WENO schemes stemmed from a loss of accuracy at critical points, and a way out of this problem was presented in Henrick, Aslam & Powers [28], Borges *et al*. [12] and Castro *et al*. [13]. For a comprehensive review of WENO schemes, see Shu [42] [43].

For a long time now, progress on finite difference WENO schemes has focused on improving their performance for hyperbolic systems of conservation laws. The original approach by Shu and Osher [40], [41] was based on making a flux vector splitting of the LLF flux into left- and right-going contributions and then carrying out an upwinded reconstruction of those fluxes. However, in recent years, newer hyperbolic PDE systems have come to be of interest which have non-conservative products. As a result, it is very desirable to have efficient finite difference WENO methods that can evolve hyperbolic systems with non-conservative products. In the literature, the update of PDEs with non-conservative products is mainly done in a fluctuation form, which is quite different from a flux form. Since the original Shu and Osher [40], [41] formulation of finite difference WENO was based on a flux form, it is understandable that the extension to fluctuation form must have seemed difficult. As far as we know, the first effort to produce a finite difference WENO scheme for treating hyperbolic systems with non-conservative products was made by Balsara *et al*. [10]. The method we presented was cast in fluctuation form and was, therefore, capable of treating PDEs that are in conservation form as well as PDEs with non-conservative products. The problem with a fluctuation form is that unless one engineers the situation perfectly it is hard to retrieve a flux conservative form from a fluctuation form. Thanks to the high accuracy of the WENO methods, we were able to show that the method in Balsara *et al*. [10] was effectively conservative for all intents and purposes. In other words, even if a problem was dominated by strong shocks, conservation was preserved to a high level of accuracy in all situations where it should have been preserved. However, the method was not exactly conservative down to machine accuracy. In this paper we precisely engineer the transcription from flux form to non-conservative form and vice versa. Consequently, this paper presents the first finite difference WENO methods that are exactly conservative, i.e. they retain a flux form, when conservation is warranted; but the methods presented here can nevertheless handle PDEs with non-conservative products.

The paper by Shu and Osher [41] contained a thorough description of the classical finite difference WENO (FD-WENO) that has by now become very successful and extremely well-known. This well-known FD-WENO formulation is based on the aforementioned flux vector splitting. However, the first paper by Shu and Osher [40] also contained the description of another alternative finite difference WENO (AFD-WENO) algorithm that was not followed up to any significant extent later, until quite recently. A paper by Merriman [34] made some headway in clarifying the AFD-WENO algorithm. Jiang, Shu and Zhang [30] labeled the AFD-WENO algorithm in Shu and Osher [40] as the "alternative formulation of the finite difference WENO scheme". (To explain further, it is an alternative algorithm to the original, well-known and well-used finite difference WENO algorithm from Shu and Osher [41].) Interest in AFD-WENO has been sporadic (Jiang, Shu and Zhang [30], [31], Zheng, Shu, Qiu [48], Gao *et al*. [25]). In Balsara *et al*. [11] we presented an AFD-WENO formulation for hyperbolic conservation laws that overcame its three major weaknesses. First, it was shown that the entire algorithm in conservation form could be derived by using a computer algebra system, which demystifies the derivation of the equations. Second, AFD-WENO relies on WENO interpolation rather than on the much better known WENO reconstruction. Therefore, in Balsara *et al*. [11] we provided all the explicit formulae needed for implementing AFD-WENO with the use of WENO interpolation. Third, the AFD-WENO algorithm requires higher order derivatives of the fluxes to be available at zone boundaries. Since those derivatives are usually obtained by finite differencing the zone-centered fluxes, they can become a source of spurious oscillations when the solution is non-smooth. However, please note that the inclusion of those fluxes is also crucially important for preserving the order property when the solution is smooth. Balsara *et al*. [11] invented a novel WENO interpolation which takes the first derivatives of the fluxes at zone centers as its inputs and returns the requisite non-linearly hybridized higher order derivatives of flux-like terms at the zone boundaries as its output. This non-linear hybridization stabilizes the evaluation of the higher derivatives of the fluxes at zone boundaries. Removal of these three barriers makes the implementation of AFD-WENO for conservation laws much more accessible to the greater community.

A good formulation of AFD-WENO has many desirable traits. In recent years we have seen many different Riemann solvers emerge which have special attributes that make them very useful in various application areas. Those Riemann solvers do not fit well into the strictures of classical FD-WENO because of its focus on reconstructing the fluxes. The AFD-WENO algorithm is free of such strictures – any type of Riemann solver can be invoked in a pointwise fashion at the zone boundaries. This makes a well-designed AFD-WENO very broadly applicable to many application areas. Classical FD-WENO also does not take well to preserving the free stream condition on curvilinear meshes; whereas AFD-WENO can indeed take well to curvilinear meshes (Jiang, Shu and Zhang [30], [31]). As a result, we see that it would be very desirable to arrive at an AFD-WENO formulation in fluctuation form that is carefully engineered to retrieve the flux form when that is warranted and nevertheless extends to non-conservative products. The method we arrive at is such that when non-conservative products are absent it reverts exactly to the method in Balsara *et al*. [11] which is in exact conservation form. Arriving at such a formulation for AFD-WENO is indeed the primary *goal* of this paper. Because the method reverts to the algorithm in Balsara *et al*. [11] when non-conservative products are absent, we do not display any results from conservation laws in this paper. Also note that the ability to transition to a precise flux conservation form when non-conservative products are absent ensures, via the Lax-Wendroff theorem, that shock locations will be exactly captured by the method. For many hyperbolic systems some of the components of the solution vector are in flux conservation form while other components might be affected by non-conservative products. In such situations, the AFD-WENO algorithm in this paper is so designed that the components of the solution vector that are in conservation form will indeed be evolved by the scheme in a fully conservative fashion. For this reason, all the examples that we will present in the later sections will consist of hyperbolic systems with at least some non-conservative products. When non-conservative products are completely absent, the AFD-WENO scheme we present here will perform identically to the AFD-WENO scheme in Balsara *et al*. [11]. In that paper we show numerous examples that are drawn from several different conservation laws.

Section II derives the AFD-WENO for hyperbolic systems with non-conservative products. Section III provides further detail for evaluating the non-conservative products in a stable but nevertheless higher order fashion. Section IV provides a pointwise plan for implementing the algorithm. Section V presents accuracy analysis, showing that the method meets its design accuracy. It also provides a speed comparison of different finite difference WENO formulations. Section VI shows one dimensional tests. Section VII presents multidimensional tests. Section VIII shows that the method takes well to hyperbolic systems that have non-conservative products as well as stiff source terms. Section IX presents some conclusions.

**II) Deriving AFD-WENO for Hyperbolic Systems with Non-Conservative Products**

AFD-WENO for conservation laws, along with the reconstruction strategies that support it, has been described extensively in Balsara *et al*. [11]. We consider that paper to be a background paper for this paper. A majority of hyperbolic systems with non-conservative products will have some components of the solution vector that are indeed in conservation form while the remaining components of the solution vector may have non-conservative products. In this Section we wish to arrive at a scheme that retains conservation form for the conserved components while fully accommodating the non-conservative products. In Section II.a we give some background on AFD-WENO schemes. In Section II.b we present a strategy for transitioning from a hyperbolic PDE in conservation form to a PDE with non-conservative products. In Section II.c, we show how we use this strategy to derive an update strategy for AFD-WENO that is naturally conservative for those components of the hyperbolic PDE that are in conservation form and can simultaneously accommodate a PDE that might have some non-conservative products.

**II.a) AFD-WENO for PDEs in Conservation Form**

Let us start by recalling that classical FD-WENO of the type described in Jiang and Shu [29] or Balsara and Shu [4] relies on a fundamental trick that was invented in Shu and Osher [40]. The trick consists of realizing that if one applies reconstruction to a smooth (Lipschitz continuous) flux then the finite difference scheme just reduces to a dimension-by-dimension reconstruction applied to the right- and left-going components of a numerical flux. Traditionally, the Locally Lax-Friedrichs (LLF) flux splitting is used for this. Despite the convenience and simplicity of traditional FD-WENO, it turns out that many popular Riemann solvers are not amenable to such a flux splitting. For this reason, subsequent attention has focused on another less-used idea from Shu and Osher [40] which was subsequently elaborated on by Merriman [34]. The idea consists of realizing that one can use any pointwise Riemann solver at zone boundaries as long as one can add higher order flux derivatives at zone boundaries to restore high order pointwise accuracy at the zone centers of the overall finite difference WENO scheme. This scheme is referred to as Alternative Finite Difference WENO, or AFD-WENO (Jiang, Shu and Zhang [30]). In AFD-WENO one can apply any type of Riemann solver at zone boundaries, making the method more flexible for different uses. For example, low cost Riemann solvers, like the HLLI Riemann solver of Dumbser and Balsara [24] which accurately preserves all stationary linearly degenerate discontinuities, can also be used. AFD-WENO schemes can also retain accuracy on body fitted curvilinear meshes where the preservation of the free stream condition seems to be beyond the capabilities of classical FD-WENO (Jiang, Shu and Zhang [31]). A corollary of AFD-WENO is that it is not based on dimension-by-dimension finite volume *reconstruction* of the upwinded fluxes but rather on a dimension-by-dimension pointwise *interpolation* of the conserved variables. This distinction between reconstruction and interpolation is very important. In Balsara, Garain and Shu [7] several WENO reconstruction formulae at different orders have been provided that are useful for classical FD-WENO. In Balsara *et al*. [11] analogous pointwise WENO interpolation formulae have been provided at different orders for use in AFD-WENO.

We start by writing an AFD-WENO scheme for a conservation law. We focus on the solution of a one-dimensional PDE system given by

$$\partial_t \mathbf{U} + \partial_x \mathbf{F} = 0 \tag{1}$$

Let us establish some notation. Please see Fig. 1. It shows a small sub-section of the mesh function in a few adjacent zones. The zones are labeled by "$i-1, i, i+1$" etc. and their zone centers are denoted by "$x_{i-1}, x_i, x_{i+1}$" etc. The zone boundaries of each zone "$i$" are denoted by "$x_{i-1/2}$" and "$x_{i+1/2}$" with $\Delta x = x_{i+1/2} - x_{i-1/2}$ being a constant because we have assumed a uniform mesh. The associated mesh functions are specified in pointwise fashion at the zone centers and are labeled by "$\mathbf{U}_{i-1}, \mathbf{U}_i, \mathbf{U}_{i+1}$" etc. Here "$\mathbf{U}$" is a vector of primal variables for the hyperbolic PDE that we are considering. For the moment we consider a one-dimensional mesh, but because this is a finite difference scheme, the method can be extended dimension-by-dimension to multiple dimensions. Using our familiar WENO-AO interpolation strategy, as applied to the point values of the mesh function, we can obtain a suitably high order interpolation within each zone. (The "AO" in WENO-AO stands for adaptive order, Balsara Garain and Shu [7]. See also Balsara et al. [5], [8], [9] for additional insights.) The interpolation within zone "$i$" gives us interpolated values of the mesh function, $\hat{\mathbf{U}}^-_{i+1/2}$ and $\hat{\mathbf{U}}^+_{i-1/2}$, at the right and left boundaries of the zone being considered, see Fig. 1. We will use a caret to denote such interpolated variables. Note that $\hat{\mathbf{U}}^-_{i+1/2}$ is available at the left side of the zone boundary $x_{i+1/2}$ and $\hat{\mathbf{U}}^+_{i-1/2}$ is available at the right side of the zone boundary $x_{i-1/2}$. Doing the same in zone "$i+1$", we obtain $\hat{\mathbf{U}}^+_{i+1/2}$ at the right side of the zone boundary $x_{i+1/2}$. Similarly from zone "$i-1$", we obtain $\hat{\mathbf{U}}^-_{i-1/2}$ at the left side of the zone boundary $x_{i-1/2}$. we assume that a Riemann solver with left and right states given by $\hat{\mathbf{U}}^-_{i+1/2}$ and $\hat{\mathbf{U}}^+_{i+1/2}$ is applied at zone boundary $x_{i+1/2}$ and it yields a resolved state $\mathbf{U}^*_{i+1/2}$ that overlies the zone boundary, as seen in Fig. 1. Because we wish to eventually work with systems that have non-conservative products, we can use the resolved state, and the structure of the Riemann fan, to obtain left- and right-going fluctuations at that zone boundary. We denote those fluctuations by $\mathbf{D}^{*-}\left(\hat{\mathbf{U}}^-_{i+1/2},\hat{\mathbf{U}}^+_{i+1/2}\right)$ and $\mathbf{D}^{*+}\left(\hat{\mathbf{U}}^-_{i+1/2},\hat{\mathbf{U}}^+_{i+1/2}\right)$ respectively. (When non-conservative products are absent, these fluctuations carry an amount of upwinding-related information that is comparable to a resolved flux which we denote by $\mathbf{F}^*\left(\hat{\mathbf{U}}^-_{i+1/2},\hat{\mathbf{U}}^+_{i+1/2}\right)$.) Likewise, we assume that a Riemann solver with left and right states given by $\hat{\mathbf{U}}^-_{i-1/2}$ and $\hat{\mathbf{U}}^+_{i-1/2}$ is applied at zone boundary $x_{i-1/2}$ and it produces left- and right-going fluctuations at that zone boundary which are denoted by $\mathbf{D}^{*-}\left(\hat{\mathbf{U}}^-_{i-1/2},\hat{\mathbf{U}}^+_{i-1/2}\right)$ and $\mathbf{D}^{*+}\left(\hat{\mathbf{U}}^-_{i-1/2},\hat{\mathbf{U}}^+_{i-1/2}\right)$ respectively. (When non-conservative products are absent, these fluctuations carry an amount of upwinding-related information that is comparable to a resolved flux which we denote by $\mathbf{F}^*\left(\hat{\mathbf{U}}^-_{i-1/2},\hat{\mathbf{U}}^+_{i-1/2}\right)$.) One of the advantages of the AFD-WENO algorithm is that it is agnostic to the type of Riemann solver that is used. This completes our description of Fig. 1.

Let us start from a simpler starting point by considering a hyperbolic conservation law, see eqn. (1). The discrete in space but continuous in time update equation for the AFD-WENO scheme can then be written as

$$\partial_t \mathbf{U}_i = -\frac{1}{\Delta x}\left\{\mathbf{F}^*\left(\hat{\mathbf{U}}_{i+1/2}^-, \hat{\mathbf{U}}_{i+1/2}^+\right) - \mathbf{F}^*\left(\hat{\mathbf{U}}_{i-1/2}^-, \hat{\mathbf{U}}_{i-1/2}^+\right)\right\}$$
$$-\frac{1}{\Delta x}\left\{\begin{bmatrix} -\frac{1}{24}(\Delta x)^2 \left[\partial_x^2 \mathbf{F}\right]_{i+1/2} + \frac{7}{5760}(\Delta x)^4 \left[\partial_x^4 \mathbf{F}\right]_{i+1/2} \\ -\frac{31}{967680}(\Delta x)^6 \left[\partial_x^6 \mathbf{F}\right]_{i+1/2} + \frac{127}{154828800}(\Delta x)^8 \left[\partial_x^8 \mathbf{F}\right]_{i+1/2} \end{bmatrix} - \begin{bmatrix} -\frac{1}{24}(\Delta x)^2 \left[\partial_x^2 \mathbf{F}\right]_{i-1/2} + \frac{7}{5760}(\Delta x)^4 \left[\partial_x^4 \mathbf{F}\right]_{i-1/2} \\ -\frac{31}{967680}(\Delta x)^6 \left[\partial_x^6 \mathbf{F}\right]_{i-1/2} + \frac{127}{154828800}(\Delta x)^8 \left[\partial_x^8 \mathbf{F}\right]_{i-1/2} \end{bmatrix}\right\} \quad (2)$$

The above equation was originally set down in Shu and Osher [40] but those who want a simpler computer algebra system-based derivation of it with more explanations can also see Balsara *et al*. [11]. The first curly bracket of eqn. (2) would yield a scheme whose accuracy is restricted to second order. The higher order flux derivatives are needed for raising the accuracy of the update equation to its design accuracy. The black flux derivative terms ( $\left[\partial_x^2 \mathbf{F}\right]_{i+1/2}$ and $\left[\partial_x^2 \mathbf{F}\right]_{i-1/2}$) in the above equation yield a third order scheme. In that case, these second derivatives have to be third order accurate. If the red flux derivative terms ( $\left[\partial_x^4 \mathbf{F}\right]_{i+1/2}$ and $\left[\partial_x^4 \mathbf{F}\right]_{i-1/2}$) are also included, in addition to the black terms, the scheme becomes fifth order accurate. In that case, all the applicable derivatives have to be fifth order accurate. If the blue flux derivative terms ( $\left[\partial_x^6 \mathbf{F}\right]_{i+1/2}$ and $\left[\partial_x^6 \mathbf{F}\right]_{i-1/2}$) are also included, in addition to the black and red terms, we get a seventh order scheme. In that case, all the applicable derivatives have to be seventh order accurate. If the magenta flux derivative terms ( $\left[\partial_x^8 \mathbf{F}\right]_{i+1/2}$ and $\left[\partial_x^8 \mathbf{F}\right]_{i-1/2}$) are included, in addition to the black, red and blue terms, we get a ninth order scheme. In that case, all the applicable derivatives have to be ninth order accurate. We will retain this coloring scheme in all subsequent equations that describe the update equation for this scheme. Eqn. (2) is to be used to derive an AFD-WENO scheme that can include hyperbolic systems with non-conservative products. At this stage of the discussion in this paper, we do not specify how the higher derivatives of the flux are to be obtained. For that reason, we keep the notation for the higher derivatives of the flux somewhat relaxed. We will tighten up the notation once we get to the part where the algorithmic aspects are to be discussed.

Please observe that eqn. (2) is in flux form and should, therefore, be able to capture shocks accurately. The first line in eqn. (2) comes from applying Riemann solvers at zone boundaries. These Riemann solvers contain the stabilization, i.e. the numerical dissipation, that is needed for handling discontinuous solutions. However, if only the first line in eqn. (2) is used, the scheme would be restricted to second order of accuracy. The subsequent higher order derivatives of the flux in eqn. (2) are indeed *essential* for restoring higher order accuracy to eqn. (2). It should also be noted that the higher order derivatives can cause Gibbs oscillation and should only be used when the solution is smooth. We see, therefore, that the higher order derivatives of the flux in eqn. (2) constitute a two-edged sword. When the solution is smooth, they are needed for higher order accuracy. But, when the solution is non-smooth, their suppression is needed for the numerical stabilization of the AFD-WENO scheme. In Balsara *et al.* [11] a variant of WENO interpolation was presented which can take the zone-centered point values of the fluxes as input and return non-linearly hybridized values of the higher derivatives of the fluxes at the zone boundaries. It was found that this interpolation was very useful in retaining high order accuracy when it is justified but nevertheless avoiding spurious numerical oscillations from developing in the vicinity of non-smooth solutions. This completes our broad-brush description of the AFD-WENO scheme for conservation laws.

## II.b) Insight on Transitioning to Hyperbolic PDEs with Non-Conservative Products

Let us begin with the same conservation law as in eqn. (1) but write it as

$$\partial_t \mathbf{U} + \partial_x \mathbf{F} = 0 \quad \text{with the flux splitting} \quad \mathbf{F}(\mathbf{U}) = \mathbf{F}_C(\mathbf{U}) + \mathbf{F}_{NC}(\mathbf{U}) \tag{3}$$

The subscripts "C" and "NC" indicate that in the upcoming discussion we will treat the flux $\mathbf{F}_C(\mathbf{U})$ in conservation form and we will write $\mathbf{F}_{NC}(\mathbf{U})$ as if it is only available as a non-conservative product. We can now write the characteristic matrix as

$$\mathbf{A}(\mathbf{U}) = \frac{\partial \left( \mathbf{F}_C(\mathbf{U}) + \mathbf{F}_{NC}(\mathbf{U}) \right)}{\partial \mathbf{U}} = \mathbf{B}(\mathbf{U}) + \mathbf{C}(\mathbf{U}) \quad \text{with} \quad \mathbf{B}(\mathbf{U}) \equiv \frac{\partial \mathbf{F}_C(\mathbf{U})}{\partial \mathbf{U}} \quad \text{and} \quad \mathbf{C}(\mathbf{U}) \equiv \frac{\partial \mathbf{F}_{NC}(\mathbf{U})}{\partial \mathbf{U}} \; .$$
$$\tag{4}$$

As a result, the hyperbolic PDE system in eqn. (3) can be written in any of three equivalent forms:-

$$\partial_t \mathbf{U} + \partial_x \mathbf{F}(\mathbf{U}) = 0 \quad \Leftrightarrow \quad \partial_t \mathbf{U} + \mathbf{A}(\mathbf{U}) \partial_x \mathbf{U} = 0 \quad \Leftrightarrow \quad \partial_t \mathbf{U} + \partial_x \mathbf{F}_C(\mathbf{U}) + \mathbf{C}(\mathbf{U}) \partial_x \mathbf{U} = 0 \quad (5)$$

It is the last form of eqn. (5) that gives us the essential insight because it looks exactly like a general hyperbolic PDE system that has some non-conservative products, which we write as

$$\partial_t \mathbf{U} + \partial_x \mathbf{F}(\mathbf{U}) + \mathbf{C}(\mathbf{U}) \partial_x \mathbf{U} = 0. \tag{6}$$

Compared to the last equation in (5), we have obtained eqn. (6) by just erasing the "*C*" subscript for the flux. This is going to be our strategy for starting with eqn. (2) and recasting it in a form that accommodates non-conservative products.

## II.c) Deriving an AFD-WENO Scheme with a Conservative Limit and also Accommodation for Non-Conservative Products

There is a good physics-based reason for wanting to write the final numerical update of eqn. (6) in fluctuation form. It stems from the fact that the matrix "**A**" does have real eigenvectors and a complete set of eigenvalues; but it cannot be guaranteed that the matrices "**B**" and "**C**" have real eigenvalues or a complete set of eigenvectors. As a result, when non-conservative products are present, the fluctuation form is the only way to go. Many Riemann solvers for conservation laws have been formulated so that they can be written in flux form or an entirely equivalent fluctuation form that respects conservation. Many of those Riemann solvers have been extended so that they can provide a fluctuation form even when non-conservative products are present.

When the hyperbolic PDE is in conservation form, it can be written in an entirely equivalent fluctuation form. However, the previous paragraph shows us that once we have obtained a fluctuation form, we can use it also for a hyperbolic system with non-conservative products. To that end, let us write the fluctuation form for the left-going fluctuations at zone boundary "$i+1/2$" as:-

$$\mathbf{D}^{*-}\left(\hat{\mathbf{U}}_{i+1/2}^{-}, \hat{\mathbf{U}}_{i+1/2}^{+}\right) = \mathbf{F}^{*}\left(\hat{\mathbf{U}}_{i+1/2}^{-}, \hat{\mathbf{U}}_{i+1/2}^{+}\right) - \mathbf{F}\left(\hat{\mathbf{U}}_{i+1/2}^{-}\right)$$
$$\Rightarrow \mathbf{F}^{*}\left(\hat{\mathbf{U}}_{i+1/2}^{-}, \hat{\mathbf{U}}_{i+1/2}^{+}\right) = \mathbf{D}^{*-}\left(\hat{\mathbf{U}}_{i+1/2}^{-}, \hat{\mathbf{U}}_{i+1/2}^{+}\right) + \mathbf{F}_C\left(\hat{\mathbf{U}}_{i+1/2}^{-}\right) + \mathbf{F}_{NC}\left(\hat{\mathbf{U}}_{i+1/2}^{-}\right). \tag{7}$$

We can also write the fluctuation form for the right-going fluctuations at zone boundary "$i-1/2$" as:-

$$\mathbf{D}^{*+}\left(\hat{\mathbf{U}}_{i-1/2}^{-},\hat{\mathbf{U}}_{i-1/2}^{+}\right)=\mathbf{F}\left(\hat{\mathbf{U}}_{i-1/2}^{+}\right)-\mathbf{F}^{*}\left(\hat{\mathbf{U}}_{i-1/2}^{-},\hat{\mathbf{U}}_{i-1/2}^{+}\right)$$
$$\Rightarrow \mathbf{F}^{*}\left(\hat{\mathbf{U}}_{i-1/2}^{-},\hat{\mathbf{U}}_{i-1/2}^{+}\right)=-\mathbf{D}^{*+}\left(\hat{\mathbf{U}}_{i-1/2}^{-},\hat{\mathbf{U}}_{i-1/2}^{+}\right)+\mathbf{F}_{C}\left(\hat{\mathbf{U}}_{i-1/2}^{+}\right)+\mathbf{F}_{NC}\left(\hat{\mathbf{U}}_{i-1/2}^{+}\right). \quad (8)$$

We should recall $\mathbf{F}(\mathbf{U})=\mathbf{F}_{C}(\mathbf{U})+\mathbf{F}_{NC}(\mathbf{U})$ and its usage in the above two equations. To arrive at our derivation, let us simply replace the resolved flux terms (the ones with the star superscripts) in the first curly bracket eqn. (2) with the fluctuation terms from eqns. (7) and (8). From eqn. (2), we get

$$\begin{aligned}\partial_{t}\mathbf{U}_{i}=&-\frac{1}{\Delta x}\left\{\mathbf{D}^{*-}\left(\hat{\mathbf{U}}_{i+1/2}^{-},\hat{\mathbf{U}}_{i+1/2}^{+}\right)+\mathbf{D}^{*+}\left(\hat{\mathbf{U}}_{i-1/2}^{-},\hat{\mathbf{U}}_{i-1/2}^{+}\right)\right\}\\ &-\frac{1}{\Delta x}\left\{\mathbf{F}_{C}\left(\hat{\mathbf{U}}_{i+1/2}^{-}\right)-\mathbf{F}_{C}\left(\hat{\mathbf{U}}_{i-1/2}^{+}\right)\right\}-\frac{1}{\Delta x}\left\{\mathbf{F}_{NC}\left(\hat{\mathbf{U}}_{i+1/2}^{-}\right)-\mathbf{F}_{NC}\left(\hat{\mathbf{U}}_{i-1/2}^{+}\right)\right\}\\ &-\frac{1}{\Delta x}\left\{\begin{array}{l}\left[\begin{array}{l}-\frac{1}{24}(\Delta x)^{2}\left[\partial_{x}^{2}\mathbf{F}\right]_{i+1/2}+\frac{7}{5760}(\Delta x)^{4}\left[\partial_{x}^{4}\mathbf{F}\right]_{i+1/2}\\ -\frac{31}{967680}(\Delta x)^{6}\left[\partial_{x}^{6}\mathbf{F}\right]_{i+1/2}+\frac{127}{154828800}(\Delta x)^{8}\left[\partial_{x}^{8}\mathbf{F}\right]_{i+1/2}\end{array}\right]\\ -\left[\begin{array}{l}-\frac{1}{24}(\Delta x)^{2}\left[\partial_{x}^{2}\mathbf{F}\right]_{i-1/2}+\frac{7}{5760}(\Delta x)^{4}\left[\partial_{x}^{4}\mathbf{F}\right]_{i-1/2}\\ -\frac{31}{967680}(\Delta x)^{6}\left[\partial_{x}^{6}\mathbf{F}\right]_{i-1/2}+\frac{127}{154828800}(\Delta x)^{8}\left[\partial_{x}^{8}\mathbf{F}\right]_{i-1/2}\end{array}\right]\end{array}\right\}\end{aligned} \quad (9)$$

Now we realize that the $-\left\{\mathbf{F}_{C}\left(\hat{\mathbf{U}}_{i+1/2}^{-}\right)-\mathbf{F}_{C}\left(\hat{\mathbf{U}}_{i-1/2}^{+}\right)\right\}/\Delta x$ term, along with the fluctuation term $-\left\{\mathbf{D}^{*-}\left(\hat{\mathbf{U}}_{i+1/2}^{-},\hat{\mathbf{U}}_{i+1/2}^{+}\right)+\mathbf{D}^{*+}\left(\hat{\mathbf{U}}_{i-1/2}^{-},\hat{\mathbf{U}}_{i-1/2}^{+}\right)\right\}/\Delta x$, in the above equation will give us a very nice conservation property for those components of the hyperbolic PDE that are genuinely in conservation form. For that reason, we leave those terms as they are. However, in our current philosophy, the term $-\left\{\mathbf{F}_{NC}\left(\hat{\mathbf{U}}_{i+1/2}^{-}\right)-\mathbf{F}_{NC}\left(\hat{\mathbf{U}}_{i-1/2}^{+}\right)\right\}/\Delta x$ is just a placeholder for a non-conservative product and it will eventually have to be represented in terms of the matrix of non-conservative products "$\mathbf{C}$" in eqn. (9). To that end, we perform a Taylor series expansion for the placeholder term and write it as follows:-

$$-\frac{1}{\Delta x}\left\{\mathbf{F}_{NC}\left(\hat{\mathbf{U}}_{i+1/2}^{-}\right)-\mathbf{F}_{NC}\left(\hat{\mathbf{U}}_{i-1/2}^{+}\right)\right\}$$

$$\cong -\left(\partial_x \mathbf{F}_{NC}\right)_i - \left\{\begin{array}{l}+\dfrac{\Delta x^2}{24}\left(\partial_x^3 \mathbf{F}_{NC}\right)_i + \dfrac{\Delta x^4}{1920}\left(\partial_x^5 \mathbf{F}_{NC}\right)_i \\ +\dfrac{\Delta x^6}{322560}\left(\partial_x^7 \mathbf{F}_{NC}\right)_i + \dfrac{\Delta x^8}{92897280}\left(\partial_x^9 \mathbf{F}_{NC}\right)_i\end{array}\right\} \qquad (10)$$

The above equation can also be written in terms of the higher order derivatives of $\mathbf{F}_{NC}$ that are evaluated at the zone boundaries. This is useful because it will enable us to make a simplification in eqn. (9). We, therefore, write the above equation as

$$-\frac{1}{\Delta x}\left\{\mathbf{F}_{NC}\left(\hat{\mathbf{U}}_{i+1/2}^{-}\right)-\mathbf{F}_{NC}\left(\hat{\mathbf{U}}_{i-1/2}^{+}\right)\right\}$$

$$\cong -\mathbf{C}(\mathbf{U}_i)\left(\partial_x \hat{\mathbf{U}}\right)_i - \frac{1}{\Delta x}\left\{\begin{array}{l}\left[\dfrac{1}{24}(\Delta x)^2\left[\partial_x^2 \mathbf{F}_{NC}\right]_{i+1/2} - \dfrac{7}{5760}(\Delta x)^4\left[\partial_x^4 \mathbf{F}_{NC}\right]_{i+1/2}\right. \\ \left.+\dfrac{31}{967680}(\Delta x)^6\left[\partial_x^6 \mathbf{F}_{NC}\right]_{i+1/2} - \dfrac{127}{154828800}(\Delta x)^8\left[\partial_x^8 \mathbf{F}_{NC}\right]_{i+1/2}\right] \\ -\left[\dfrac{1}{24}(\Delta x)^2\left[\partial_x^2 \mathbf{F}_{NC}\right]_{i-1/2} - \dfrac{7}{5760}(\Delta x)^4\left[\partial_x^4 \mathbf{F}_{NC}\right]_{i-1/2}\right. \\ \left.+\dfrac{31}{967680}(\Delta x)^6\left[\partial_x^6 \mathbf{F}_{NC}\right]_{i-1/2} - \dfrac{127}{154828800}(\Delta x)^8\left[\partial_x^8 \mathbf{F}_{NC}\right]_{i-1/2}\right]\end{array}\right\} \qquad (11)$$

In the above two equations we use the "$\cong$" sign instead of the "$=$" sign because the two sides of the equation are not exactly identical; but they are only identical to the level of discretization error. This also shows us that this is the step where exact conservation can be lost. However, by localizing this loss of conservation to those components of the solution vector that have non-zero values for $\left(\mathbf{C}\partial_x \hat{\mathbf{U}}\right)$, we see that the components of the solution vector that are in conservation form will indeed retain exact flux conservation. This is a very important property because it ensures that when exact conservation is guaranteed by the PDE, our AFD-WENO scheme will indeed retain exact flux conservation in a numerical simulation up to machine precision.

Now by amalgamating eqns. (9) and (11) we get

$$\partial_t \mathbf{U}_i = -\frac{1}{\Delta x}\left\{\mathbf{D}^{*-}\left(\hat{\mathbf{U}}_{i+1/2}^-, \hat{\mathbf{U}}_{i+1/2}^+\right) + \mathbf{D}^{*+}\left(\hat{\mathbf{U}}_{i-1/2}^-, \hat{\mathbf{U}}_{i-1/2}^+\right)\right\} - \frac{1}{\Delta x}\left\{\mathbf{F}_C\left(\hat{\mathbf{U}}_{i+1/2}^-\right) - \mathbf{F}_C\left(\hat{\mathbf{U}}_{i-1/2}^+\right)\right\} - \mathbf{C}(\mathbf{U}_i)\left(\partial_x \hat{\mathbf{U}}\right)_i$$

$$-\frac{1}{\Delta x}\left\{\begin{bmatrix} -\frac{1}{24}(\Delta x)^2 \left[\partial_x^2 \mathbf{F}_C\right]_{i+1/2} + \frac{7}{5760}(\Delta x)^4 \left[\partial_x^4 \mathbf{F}_C\right]_{i+1/2} \\ -\frac{31}{967680}(\Delta x)^6 \left[\partial_x^6 \mathbf{F}_C\right]_{i+1/2} + \frac{127}{154828800}(\Delta x)^8 \left[\partial_x^8 \mathbf{F}_C\right]_{i+1/2} \end{bmatrix} \\ -\begin{bmatrix} -\frac{1}{24}(\Delta x)^2 \left[\partial_x^2 \mathbf{F}_C\right]_{i-1/2} + \frac{7}{5760}(\Delta x)^4 \left[\partial_x^4 \mathbf{F}_C\right]_{i-1/2} \\ -\frac{31}{967680}(\Delta x)^6 \left[\partial_x^6 \mathbf{F}_C\right]_{i-1/2} + \frac{127}{154828800}(\Delta x)^8 \left[\partial_x^8 \mathbf{F}_C\right]_{i-1/2} \end{bmatrix}\right\}$$

(12)

The above equation still pertains to eqn. (5). We now realize that in order to transition from eqn. (5) to eqn. (6) we only need to make the transcription $\mathbf{F}_C \to \mathbf{F}$. Therefore, when working with eqn. (6) we have

$$\partial_t \mathbf{U}_i = -\frac{1}{\Delta x}\left\{\mathbf{D}^{*-}\left(\hat{\mathbf{U}}_{i+1/2}^-, \hat{\mathbf{U}}_{i+1/2}^+\right) + \mathbf{D}^{*+}\left(\hat{\mathbf{U}}_{i-1/2}^-, \hat{\mathbf{U}}_{i-1/2}^+\right)\right\} - \frac{1}{\Delta x}\left\{\mathbf{F}\left(\hat{\mathbf{U}}_{i+1/2}^-\right) - \mathbf{F}\left(\hat{\mathbf{U}}_{i-1/2}^+\right)\right\} - \mathbf{C}(\mathbf{U}_i)\left[\partial_x \hat{\mathbf{U}}\right]_i$$

$$-\frac{1}{\Delta x}\left\{\begin{bmatrix} -\frac{1}{24}(\Delta x)^2 \left[\partial_x^2 \hat{\mathbf{F}}\right]_{i+1/2} + \frac{7}{5760}(\Delta x)^4 \left[\partial_x^4 \hat{\mathbf{F}}\right]_{i+1/2} \\ -\frac{31}{967680}(\Delta x)^6 \left[\partial_x^6 \hat{\mathbf{F}}\right]_{i+1/2} + \frac{127}{154828800}(\Delta x)^8 \left[\partial_x^8 \hat{\mathbf{F}}\right]_{i+1/2} \end{bmatrix} \\ -\begin{bmatrix} -\frac{1}{24}(\Delta x)^2 \left[\partial_x^2 \hat{\mathbf{F}}\right]_{i-1/2} + \frac{7}{5760}(\Delta x)^4 \left[\partial_x^4 \hat{\mathbf{F}}\right]_{i-1/2} \\ -\frac{31}{967680}(\Delta x)^6 \left[\partial_x^6 \hat{\mathbf{F}}\right]_{i-1/2} + \frac{127}{154828800}(\Delta x)^8 \left[\partial_x^8 \hat{\mathbf{F}}\right]_{i-1/2} \end{bmatrix}\right\}$$

(13)

Till eqn. (12) we had been willing to accept a relaxed notation for the higher derivatives of the fluxes at the zone boundaries because we were just engaged in a derivation. In eqn. (13) and onwards, we have tightened the notation. Therefore, the caret in $\left[\partial_x^2 \hat{\mathbf{F}}\right]_{i+1/2}$, and terms like it, indicates that is term is to be obtained via some form of WENO interpolation and the interpolant is to be evaluated at zone boundary "$i+1/2$". The exact form of the interpolation will be discussed in the next Section. The above equation is potentially very useful in several situations where there is a weak coupling between the conservative terms and the non-conservative products. For example, it could be very useful for large combinations of PDE systems like general relativistic hydrodynamics where there is a clear split between the Einstein field equations and the equations of relativistic hydrodynamics. However, for tightly coupled systems, i.e. systems where the flux terms and non-conservative terms interact strongly with one another, it still has a small deficiency which we explain below.

Please look at eqn. (13) and notice the derivative $\left[\partial_x \hat{\mathbf{U}}\right]_i$ in the term $-\mathbf{C}(\mathbf{U}_i)\left[\partial_x \hat{\mathbf{U}}\right]_i$. That derivative will be evaluated using WENO interpolation. However, the non-linear stabilization in WENO is applied to the interpolant, not to its derivative. Therefore, it is slightly advantageous to write eqn. (13) in an equivalent form given by

$$\partial_t \mathbf{U}_i = -\frac{1}{\Delta x}\left\{\mathbf{D}_{i+1/2}^{*-}\left(\hat{\mathbf{U}}_{i+1/2}^-, \hat{\mathbf{U}}_{i+1/2}^+\right) + \mathbf{D}_{i-1/2}^{*+}\left(\hat{\mathbf{U}}_{i-1/2}^-, \hat{\mathbf{U}}_{i-1/2}^+\right)\right\} - \frac{1}{\Delta x}\left\{\mathbf{F}\left(\hat{\mathbf{U}}_{i+1/2}^-\right) - \mathbf{F}\left(\hat{\mathbf{U}}_{i-1/2}^+\right)\right\}$$

$$-\frac{1}{\Delta x}\mathbf{C}(\mathbf{U}_i)\left(\hat{\mathbf{U}}_{i+1/2}^- - \hat{\mathbf{U}}_{i-1/2}^+\right) + \frac{1}{\Delta x}\mathbf{C}(\mathbf{U}_i)\left\{\begin{array}{l}\frac{1}{24}(\Delta x)^3\left[\partial_x^3 \hat{\mathbf{U}}\right]_i + \frac{1}{1920}(\Delta x)^5\left[\partial_x^5 \hat{\mathbf{U}}\right]_i \\ +\frac{1}{322560}(\Delta x)^7\left[\partial_x^7 \hat{\mathbf{U}}\right]_i + \frac{1}{92897280}(\Delta x)^9\left[\partial_x^9 \hat{\mathbf{U}}\right]_i\end{array}\right\}$$

$$-\frac{1}{\Delta x}\left\{\begin{array}{l}\left[\begin{array}{l}-\frac{1}{24}(\Delta x)^2\left[\partial_x^2 \hat{\mathbf{F}}\right]_{i+1/2} + \frac{7}{5760}(\Delta x)^4\left[\partial_x^4 \hat{\mathbf{F}}\right]_{i+1/2} \\ -\frac{31}{967680}(\Delta x)^6\left[\partial_x^6 \hat{\mathbf{F}}\right]_{i+1/2} + \frac{127}{154828800}(\Delta x)^8\left[\partial_x^8 \hat{\mathbf{F}}\right]_{i+1/2}\end{array}\right] \\ -\left[\begin{array}{l}-\frac{1}{24}(\Delta x)^2\left[\partial_x^2 \hat{\mathbf{F}}\right]_{i-1/2} + \frac{7}{5760}(\Delta x)^4\left[\partial_x^4 \hat{\mathbf{F}}\right]_{i-1/2} \\ -\frac{31}{967680}(\Delta x)^6\left[\partial_x^6 \hat{\mathbf{F}}\right]_{i-1/2} + \frac{127}{154828800}(\Delta x)^8\left[\partial_x^8 \hat{\mathbf{F}}\right]_{i-1/2}\end{array}\right]\end{array}\right\}$$

(14)

Notice from the above equation that the terms $\hat{\mathbf{U}}_{i+1/2}^-$ and $\hat{\mathbf{U}}_{i-1/2}^+$ in $-\mathbf{C}(\mathbf{U}_i)\left(\hat{\mathbf{U}}_{i+1/2}^- - \hat{\mathbf{U}}_{i-1/2}^+\right)/\Delta x$ can now be obtained by a good non-linearly hybridized WENO interpolation process. Having seen that the flux terms can be written as finite difference approximation of suitable combinations of higher order derivatives at zone boundaries in the above equation, we want to use the same idea for the higher order derivatives given by $\left[\partial_x^3 \hat{\mathbf{U}}\right]_i$, $\left[\partial_x^5 \hat{\mathbf{U}}\right]_i$, $\left[\partial_x^7 \hat{\mathbf{U}}\right]_i$, and $\left[\partial_x^9 \hat{\mathbf{U}}\right]_i$ in eqn. (14). This is a reasonable thing to do because we have the intuition that a finite difference approximation is numerically more stable than the evaluation of higher derivatives at zone centers. For that reason, we rewrite eqn. (14) as

$$\partial_t \mathbf{U}_i = -\frac{1}{\Delta x}\left\{\mathbf{D}^{*-}_{i+1/2}\left(\hat{\mathbf{U}}^-_{i+1/2}, \hat{\mathbf{U}}^+_{i+1/2}\right) + \mathbf{D}^{*+}_{i-1/2}\left(\hat{\mathbf{U}}^-_{i-1/2}, \hat{\mathbf{U}}^+_{i-1/2}\right)\right\} - \frac{1}{\Delta x}\left\{\mathbf{F}\left(\hat{\mathbf{U}}^-_{i+1/2}\right) - \mathbf{F}\left(\hat{\mathbf{U}}^+_{i-1/2}\right)\right\}$$

$$-\frac{1}{\Delta x}\mathbf{C}(\mathbf{U}_i)\left(\hat{\mathbf{U}}^-_{i+1/2} - \hat{\mathbf{U}}^+_{i-1/2}\right)$$

$$-\frac{1}{\Delta x}\mathbf{C}(\mathbf{U}_i)\left\{\begin{array}{l}\left[\begin{array}{l}-\frac{1}{24}(\Delta x)^2\left[\partial_x^2 \hat{\mathbf{U}}\right]_{i+1/2} + \frac{7}{5760}(\Delta x)^4\left[\partial_x^4 \hat{\mathbf{U}}\right]_{i+1/2} \\ -\frac{31}{967680}(\Delta x)^6\left[\partial_x^6 \hat{\mathbf{U}}\right]_{i+1/2} + \frac{127}{154828800}(\Delta x)^8\left[\partial_x^8 \hat{\mathbf{U}}\right]_{i+1/2}\end{array}\right] \\ -\left[\begin{array}{l}-\frac{1}{24}(\Delta x)^2\left[\partial_x^2 \hat{\mathbf{U}}\right]_{i-1/2} + \frac{7}{5760}(\Delta x)^4\left[\partial_x^4 \hat{\mathbf{U}}\right]_{i-1/2} \\ -\frac{31}{967680}(\Delta x)^6\left[\partial_x^6 \hat{\mathbf{U}}\right]_{i-1/2} + \frac{127}{154828800}(\Delta x)^8\left[\partial_x^8 \hat{\mathbf{U}}\right]_{i-1/2}\end{array}\right]\end{array}\right\}$$

$$-\frac{1}{\Delta x}\left\{\begin{array}{l}\left[\begin{array}{l}-\frac{1}{24}(\Delta x)^2\left[\partial_x^2 \hat{\mathbf{F}}\right]_{i+1/2} + \frac{7}{5760}(\Delta x)^4\left[\partial_x^4 \hat{\mathbf{F}}\right]_{i+1/2} \\ -\frac{31}{967680}(\Delta x)^6\left[\partial_x^6 \hat{\mathbf{F}}\right]_{i+1/2} + \frac{127}{154828800}(\Delta x)^8\left[\partial_x^8 \hat{\mathbf{F}}\right]_{i+1/2}\end{array}\right] \\ -\left[\begin{array}{l}-\frac{1}{24}(\Delta x)^2\left[\partial_x^2 \hat{\mathbf{F}}\right]_{i-1/2} + \frac{7}{5760}(\Delta x)^4\left[\partial_x^4 \hat{\mathbf{F}}\right]_{i-1/2} \\ -\frac{31}{967680}(\Delta x)^6\left[\partial_x^6 \hat{\mathbf{F}}\right]_{i-1/2} + \frac{127}{154828800}(\Delta x)^8\left[\partial_x^8 \hat{\mathbf{F}}\right]_{i-1/2}\end{array}\right]\end{array}\right\}$$

(15)

The color coding in eqns. (9) to (15), as it pertains to order of accuracy, is identical to the color coding in eqn. (2). Eqn. (15) is the AFD-WENO update equation that we are seeking. It contains all the update terms that will be needed in AFD-WENO schemes that are up to ninth order accurate. To obtain a third order AFD-WENO scheme from eqn. (15), please retain only the black terms and eliminate the red, blue and magenta terms. To obtain a fifth order AFD-WENO scheme from eqn. (15), please retain only the black and red terms and eliminate the blue and magenta terms. To obtain a seventh order AFD-WENO scheme from eqn. (15), please retain only the black, red and blue terms and eliminate the magenta terms. To obtain a ninth order AFD-WENO scheme from eqn. (15), please retain all the terms in that equation. The computer algebra system script in Appendix A of Balsara *et al*. [11] can be extended to give even higher orders. In Appendix A of this paper we provide an alternative to eqn. (15) which could very slightly reduce the cost of evaluating the last curly bracket in eqn. (15).

Eqn. (15) has such a nice structure that when the matrix of non-conservative products is zero, i.e. when we have "$\mathbf{C} = 0$", then eqn. (15) reduces identically to eqn. (2) which is in manifestly flux conservation form. When only a few rows of the matrix "$\mathbf{C}$" are non-zero, eqn. (15) retains a nice structure which ensures that the remaining rows of the hyperbolic PDE will indeed remain in flux conservation form. Therefore, eqn. (15) retains flux conservation form when such a conservation form is present in the governing PDE; and eqn. (15) nevertheless permits the incorporation of non-conservative products when the PDE has such products. Furthermore, eqn. (15) is in finite difference form, which makes it very efficient for high order accurate computations in multiple dimensions. Also notice that eqn. (15) can be used with any Riemann solver that can be cast in fluctuation form. Since many of the well-known Riemann solvers have a flux form as well as a fluctuation form, they can all be used in eqn. (15). Different Riemann solvers may have some special properties that make them useful in certain fields of study. This makes eqn. (15) broadly applicable to many PDEs in many different fields of study.

To arrive at a non-linearly stabilized version of eqn. (15), we will need three WENO interpolation steps. Notice that the first WENO interpolation is needed to obtain the non-linearly stabilized and interpolated solution at various locations on the mesh. This will give us the first line (i.e. the first three terms) of eqn. (15), which indeed retains up to second order of accuracy. It is very desirable that this interpolation be carried out in the eigenspace of the governing equation. The first three terms of eqn. (15), if they are used by themselves, also guarantee that we will get a very high quality second order scheme if a high quality WENO interpolation is used. The last two curly brackets of eqn. (15) only contain higher order derivatives that are needed for raising the accuracy of the AFD-WENO scheme to the desired accuracy when the solution is smooth. But those higher order derivatives are also a two-edged sword because they can introduce spurious oscillations, thereby damaging the accuracy, when the solution is non-smooth. We seek to stabilize them using a different style of WENO interpolation from Balsara *et al*. [11]. The second WENO interpolation is needed for obtaining the non-linearly stabilized higher order derivatives, $\left[\partial_x^2 \hat{\mathbf{F}}\right]_{i+1/2}$, $\left[\partial_x^4 \hat{\mathbf{F}}\right]_{i+1/2}$, $\left[\partial_x^6 \hat{\mathbf{F}}\right]_{i+1/2}$ and $\left[\partial_x^8 \hat{\mathbf{F}}\right]_{i+1/2}$ (as needed) at zone boundary "$i+1/2$"; with analogous evaluations needed at the zone boundary "$i-1/2$". The third WENO interpolation is needed for obtaining the non-linearly stabilized higher order derivatives, $\left[\partial_x^2 \hat{\mathbf{U}}\right]_{i+1/2}$, $\left[\partial_x^4 \hat{\mathbf{U}}\right]_{i+1/2}$, $\left[\partial_x^6 \hat{\mathbf{U}}\right]_{i+1/2}$ and $\left[\partial_x^8 \hat{\mathbf{U}}\right]_{i+1/2}$ (as needed) at the zone boundaries. The second and third WENO interpolations can be carried out in physical space, so that their cost can be minimized. It is also useful to emphasize that, eqn. (15) is not the only formulation that is possible. In some circumstances, and for some PDEs where the non-conservative products are not strongly non-linear, eqn. (13) can yield a simpler formulation that requires only two WENO interpolation steps.

Our experience has been that the WENO interpolation is sufficient for controlling the higher order derivative terms in eqn. (15). WENO interpolation represents a *numerics-based modulation* of the solution and applies generally to all PDEs. But for a very small number of very stringent hyperbolic systems with non-conservative products a *physics-based modulation* of the interpolation requires the use of a flattener algorithm (Colella and Woodward [17], Balsara [3], [6]). These physics-based flattener algorithms are PDE-specific and they are designed to do nothing when the solution is smooth or only mildly non-smooth; and activate themselves only when the solution is very non-smooth. So, the higher derivative terms in eqn. (15) have to be used with care.

**III) Detailed Strategy for Obtaining the Higher Order Derivatives in Eqn. (15)**

Taken by itself, the first three terms of eqn. (15) would only achieve second order of accuracy in space. Eqn. (15) shows that the final AFD-WENO scheme requires several higher order derivatives in order to achieve its higher order spatial accuracy. These derivatives are essential if higher order accuracy has to be achieved by the overall scheme when it is used to treat a solution which is smooth on the mesh. However, if the solution is non-smooth on the mesh, these higher order derivatives can even be a source of unphysical oscillations. We seek an automatic process that evaluates these derivatives with sufficient accuracy to meet the design accuracy of the scheme when the solution is smooth. However, this automatic process should also suppress these derivatives when they are likely to trigger spurious oscillations.

In Balsara *et al*. [11] a novel type of WENO interpolation was designed which does just that. Fig. 2, which is modified from Balsara *et al*. [11], shows how the set of zone-centered fluxes $\{\mathbf{F}(\mathbf{U}_i)\}$ are use to obtain the higher order derivatives of the flux at the zone boundaries. But realize from examining eqn. (15) that higher derivatives of $\{\mathbf{F}(\mathbf{U}_i)\}$ are needed at the zone boundaries where they provide higher order contributions to the flux terms when the solution is smooth. Section IV of Balsara *et al.* [11] provides such a WENO interpolation strategy and provides all the explicit formulae for obtaining the result. This interpolation strategy is illustrated in Fig. 2. For a third order AFD-WENO scheme, $\left[\partial_x^2 \hat{\mathbf{F}}\right]_{i+1/2}$ is needed at zone boundary "i+1/2" with third order of accuracy. Therefore, the magenta colored left-biased stencil and the blue colored right-biased stencil in Fig. 2 can be non-linearly hybridized in order to get $\left[\partial_x^2 \hat{\mathbf{F}}\right]_{i+1/2}$. For a fifth order AFD-WENO scheme, $\left[\partial_x^2 \hat{\mathbf{F}}\right]_{i+1/2}$ and $\left[\partial_x^4 \hat{\mathbf{F}}\right]_{i+1/2}$ both need to be evaluated at zone boundary "i+1/2" with fifth order of accuracy. To achieve that, the magenta colored left-biased stencil, the blue colored right-biased stencil and the red colored large sixth order accurate central stencil in Fig. 2 need to be non-linearly hybridized in order to get $\left[\partial_x^2 \hat{\mathbf{F}}\right]_{i+1/2}$ and $\left[\partial_x^4 \hat{\mathbf{F}}\right]_{i+1/2}$ at the zone boundaries. The exact details of this novel WENO interpolation are given in Section IV of Balsara *et al.* [11]. This paragraph, along with Fig. 2, has shown us how the last curly bracket term in eqn. (15) is obtained. It gives us the higher order contributions to the flux terms when the solution is smooth.

Extensive numerical experimentation in Balsara *et al.* [11] shows that the AFD-WENO schemes in conservation form are indeed very robust when applied to conservation laws. This is partly because the Lax-Wendroff theorem comes to the rescue. This is also why we wrote the higher order non-conservative products in eqn. (15) in a form that is as close as possible to a finite difference-like form. We see from Fig. 3 that the same strategy that was used in Fig. 2 for obtaining higher order derivatives of the fluxes at the zone boundaries can also be used for obtaining the higher derivatives of the state "$\mathbf{U}$" at the zone boundaries.

This completes our description of the non-linearly hybridized WENO interpolation processes that will be used to extract all the higher derivatives that we see in eqn. (15).

**IV) Pointwise Implementation of Our AFD-WENO Scheme for Hyperbolic PDEs with Non-Conservative Products**

The AFD-WENO scheme that we present here is based on the WENO interpolation methods described in Section III of Balsara *et al*. [11]. In that section, we present Adaptive Order as well as Multiresolution methods for starting from point values at zone centers and using them to obtain interpolated values at the zone boundaries for feeding to a Riemann solver; as shown in Fig. 1. We will also draw on WENO interpolation methods described in Section IV of Balsara *et al*. [11] where we show that the set of flux terms $\{\mathbf{F}(\mathbf{U}_i)\}$ that are evaluated at zone centers can have their higher derivatives interpolated to the zone boundaries, as shown in Fig. 2. We will also use similar methods to take the zone-centered state vector "$\mathbf{U}$" and evaluate its higher derivatives at zone boundaries; see Fig. 3. While the WENO interpolation in Sections III and IV of Balsara *et al*. [11] is couched in the language of WENO-AO, all the results there can be easily transcribed to the Multiresolution WENO from (Zhu and Shu [50]); as a result we will show results from WENO-AO and Multiresolution WENO in this paper, displaying that in all cases their effective accuracies are almost identical.

The Riemann solver that we will use in this entire work will be the HLLI Riemann solver from Dumbser and Balsara [24] where the philosophy of such Riemann solvers is explained. In Section V of Balsara *et al*. [10] we also catalogue results from Dumbser and Balsara [24] in a notation that is more suited for the WENO interpolation shown in Fig. 1.

We provide a pointwise implementation of our AFD-WENO scheme for treating non-linear hyperbolic PDEs that have non-conservative products. Some of the steps are similar to Balsara *et al*. [11] and some are different. AFD-WENO schemes are always implemented in dimension-by-dimension fashion, so we only describe one of the dimensional updates here. We realize that the update equation, i.e. eqn. (15), has a lot of terms. The optimal sequence of steps given below is designed so that at the end of each step we catalogue the parts of eqn. (15) that are in hand. Consequently, by the end of these steps, we can finally assemble the entire update equation. The pointwise implementation of our AFD-WENO scheme into a numerical code goes according to the following steps:-

**1)** We start with the mesh function as shown in Fig. 1. This means that at each zone center $x_i$ we have a pointwise value for the conserved variable $\mathbf{U}_i$.

**2)** Starting from the conserved variables in each zone, obtain the primitive variables. Use the conserved and primitive variables, as needed, to obtain the normalized right and left eigenvectors in the conserved variables. Also evaluate $\mathbf{F}(\mathbf{U}_i)$ and $\mathbf{C}(\mathbf{U}_i)$ pointwise at the center of each zone "$i$".

**3)** As shown in Fig. 1, we use the WENO-AO algorithm from Section III of Balsara *et al.* [11]. (We have also used the Multiresolution WENO interpolation from Zhu and Shu [50] and found no difference.) That Section includes all closed form expressions that are needed for WENO interpolation in one dimension. This consists of making a non-linear hybridization between a large high order accurate stencil and smaller lower order accurate stencils. The neighboring zones around zone "$i$" are projected into the characteristic space of zone "$i$". The third and fifth order cases are explicitly shown in Fig. 1. Once the variables in the neighboring zones around zone "$i$" are projected into the characteristic space of zone "$i$", WENO-AO interpolation is carried out in the characteristic space. Projecting the interpolated characteristic variables back into the space of right eigenvectors gives us high order accurate $\hat{\mathbf{U}}^-_{i+1/2}$ and $\hat{\mathbf{U}}^+_{i-1/2}$ within each zone "$i$", as shown in Fig. 1. Since this step involves projecting all the zones in all the stencils of interest into the characteristic space of each zone "$i$" using eigenvectors, it is one of the three computationally expensive steps of the algorithm. By the end of this step we should have the WENO interpolation-based $\hat{\mathbf{U}}^-_{i+1/2}$ and $\hat{\mathbf{U}}^+_{i+1/2}$ at each zone boundary.

**4)** At each zone boundary $x_{i+1/2}$, use the left and right states $\hat{\mathbf{U}}^-_{i+1/2}$ and $\hat{\mathbf{U}}^+_{i+1/2}$ to obtain the left-most and right-most going speeds of the Riemann fan; these are denoted by $S_{L;i+1/2}$ and $S_{R;i+1/2}$. Please note that we are not yet seeking the resolved state within the Riemann fan.

**5)** Now, at each zone boundary $x_{i+1/2}$, we hand in the speeds $S_{L;i+1/2}$ and $S_{R;i+1/2}$ as well as the states $\hat{\mathbf{U}}^-_{i+1/2}$ and $\hat{\mathbf{U}}^+_{i+1/2}$ to the Riemann solver. Unlike the situation with conservation laws (where the Riemann solver returns a resolved flux), we now ask the Riemann solver to return fluctuations $\mathbf{D}^{*-}(\hat{\mathbf{U}}^-_{i+1/2}, \hat{\mathbf{U}}^+_{i+1/2})$ and $\mathbf{D}^{*+}(\hat{\mathbf{U}}^-_{i+1/2}, \hat{\mathbf{U}}^+_{i+1/2})$ as well as the resolved state $\mathbf{U}^*_{i+1/2}$ at each zone boundary "$i+1/2$". This is shown in Fig. 1. By the end of this step we should have the resolved state $\mathbf{U}^*_{i+1/2}$ as well as the corresponding fluctuations $\mathbf{D}^{*-}\left(\hat{\mathbf{U}}_{i+1/2}^{-}, \hat{\mathbf{U}}_{i+1/2}^{+}\right)$ and $\mathbf{D}^{*+}\left(\hat{\mathbf{U}}_{i+1/2}^{-}, \hat{\mathbf{U}}_{i+1/2}^{+}\right)$ at each zone boundary.

**6)** If one wants to make a characteristic projection of the higher derivatives of the flux variables, we can do that using $\mathbf{U}_{i+1/2}^{*}$. This could be useful in the next step. Therefore, we find the matrices of right and left eigenvectors corresponding to the resolved state $\mathbf{U}_{i+1/2}^{*}$ at each zone boundary "$i+1/2$". Please notice that if the HLLI Riemann solver is used, then we will naturally be constructing the left and right eigenvectors from the resolved state of that Riemann solver. Therefore, it is worthwhile to derive the maximum use from those eigenvectors.

**7)** Use the boundary-centered WENO-AO interpolation scheme from Section IV of Balsara *et al.* [11] and Fig. 2 of this paper to start with the zone-centered flux variables and interpolate their higher derivatives to the zone boundaries. This gives us suitably high order derivatives of the flux variables at each zone boundary. With these high order derivatives of the flux variables, we can evaluate all the higher order flux derivatives that contribute to each zone boundary; see the last term eqn. (15). (Since this step may involve projecting all the zones in all the stencils of interest into the characteristic space of each zone boundary "$i+1/2$" using eigenvectors, it is the second of the three computationally expensive steps of the algorithm. However, we have found this characteristic projection to be unnecessary, with the result that WENO interpolation can be applied directly to the flux components.) By the end of this step we should have higher order derivatives of the flux like $\left[\partial_x^2 \hat{\mathbf{F}}\right]_{i+1/2}$, $\left[\partial_x^4 \hat{\mathbf{F}}\right]_{i+1/2}$, $\left[\partial_x^6 \hat{\mathbf{F}}\right]_{i+1/2}$ and $\left[\partial_x^8 \hat{\mathbf{F}}\right]_{i+1/2}$ (as needed) at each of the zone boundaries. (If the PDE is a conservation law, then by this step we have all the pieces in hand that are needed for the assembly of the update equation in eqn. (15).)

**8)** This step is closely analogous to the previous step; except that we now consider the zone centered state "$\mathbf{U}$" instead of the flux. Use the boundary-centered WENO-AO interpolation scheme from Section IV of Balsara *et al.* [11] and Fig. 3 of this paper to start with the zone-centered state variables and interpolate their higher derivatives to the zone boundaries. By the end of this step we should have higher order derivatives of the state like $\left[\partial_x^2 \hat{\mathbf{U}}\right]_{i+1/2}$, $\left[\partial_x^4 \hat{\mathbf{U}}\right]_{i+1/2}$, $\left[\partial_x^6 \hat{\mathbf{U}}\right]_{i+1/2}$ and $\left[\partial_x^8 \hat{\mathbf{U}}\right]_{i+1/2}$ (as needed) at each of the zone boundaries.

**9)** Now realize from the previous steps that we have acquired all the terms that will contribute to eqn. (15). We assemble eqn. (15) which is our update equation. That gives us one spatially higher order update stage of a multistage RK update strategy.

**10)** The above points have only shown one stage of the scheme. It can be coupled with an SSP-RK update strategy, say from Shu and Osher [40] or Spiteri and Ruuth [44], [45], to achieve higher order in time.

**11)** Some of the PDEs also have stiff source terms; these are usually relaxation terms that enable the system to relax to several useful physical limits. The AFD-WENO method makes it very simple to treat stiff source terms because the source terms are treated pointwise and are collocated at the exact same location as the primal variables. For this reason, when stiff source terms are present, we recommend using the Runge-Kutta IMEX methods from Pareschi and Russo [35]; see also Kupka *et al*. [32].

We remind the reader that eqn. (15) has been explicitly written out in Appendix A as it would pertain to the implementation of third, fifth, seventh and ninth order AFD-WENO schemes. When making an implementation in code we advise the computational scientist to first try and implement the first line of eqn. (15) and make it work for several test problems. The subsequent parts of eqn. (15) can be implemented after that.

**V) Accuracy Analysis and Speed Comparisons**

In this Section we do not consider hyperbolic conservation laws. The reason is that the current algorithm reduces exactly to the algorithm in Balsara *et al*. [11] where we indeed presented many accuracy analyses as applied to several conservation laws. Sub-section V.a presents accuracy analysis for the Baer-Nunziato system. Sub-section V.b presents accuracy analysis for the two layer shallow water system. Sub-section V.c presents accuracy analysis for the multiphase debris flow model of Pitman and Le. In each instance, the test problems presented here are identical to the ones presented in Balsara *et al*. [10]. For that reason, we do not describe the test problems in great detail in this paper. In Sub-section V.d we cross-compare the speeds of classical FD-WENO, conservative AFD-WENO and the AFD-WENO scheme which includes non-conservative products.

Many of the results that follow are tabulated for both the WENO-AO-based interpolation method (see Balsara, Garain and Shu [7], Balsara *et al.* [11]) as well as the Multiresolution WENO-based method (Zhu and Shu [50], Zhu and Qiu [49]) showing that in all instances the two interpolation options yield almost identical accuracies.

**V.a) Accuracy Analysis for the two-dimensional Baer-Nunziato Model for Compressible Multi-Phase Flows**

In this Sub-section we focus on the Baer-Nunziato compressible multiphase flow proposed by Baer and Nunziato [2] and extensively studied at the numerical level by Saurel and Abgrall [39], Adrianov and Warnecke [1], Schwendeman *et al.* [53], Dumbser *et al.* [23], Tokareva and Toro [46], Coquel et al. [18] and Chiochetti and Müller [16]. A very useful set of eigenvectors have been presented in Tokareva and Toro [46].

The PDE system assumes two phases, a solid phase denoted by density $\rho_1$, volume fraction $\phi_1$, velocity $\mathbf{v}_1 = (u_1, v_1, w_1)$ and a pressure $p_1$ and a gas phase denoted by density $\rho_2$, volume fraction $\phi_2$, velocity $\mathbf{v}_2 = (u_2, v_2, w_2)$ and a pressure $p_2$. The two phases have an interfacial pressure $P_I$ and an interfacial velocity $\mathbf{V}_I$, but it is suggested in Baer and Nunziato [2] to set $P_I = p_2$ and $\mathbf{V}_I = \mathbf{v}_1$. The total energy density for phase "*j*" is related to the specific internal energy $e_j$ by $\rho_j E_j = \rho_j e_j + \rho_j \mathbf{v}_j^2 / 2$.

$$\partial_t (\phi_1 \rho_1) + \nabla \cdot (\phi_1 \rho_1 \mathbf{v}_1) = 0$$
$$\partial_t (\phi_1 \rho_1 \mathbf{v}_1) + \nabla \cdot (\phi_1 (\rho_1 \mathbf{v}_1 \otimes \mathbf{v}_1 + \mathbf{I} p_1)) - P_I \nabla \phi_1 = 0$$
$$\partial_t (\phi_1 \rho_1 E_1) + \nabla \cdot (\phi_1 \mathbf{v}_1 (\rho_1 E_1 + p_1)) + P_I \partial_t \phi_1 = 0$$
$$\partial_t (\phi_2 \rho_2) + \nabla \cdot (\phi_2 \rho_2 \mathbf{v}_2) = 0$$
$$\partial_t (\phi_2 \rho_2 \mathbf{v}_2) + \nabla \cdot (\phi_2 (\rho_2 \mathbf{v}_2 \otimes \mathbf{v}_2 + \mathbf{I} p_2)) - P_I \nabla \phi_2 = 0$$
$$\partial_t (\phi_2 \rho_2 E_2) + \nabla \cdot (\phi_2 \mathbf{v}_2 (\rho_2 E_2 + p_2)) + P_I \partial_t \phi_2 = 0$$
$$\partial_t \phi_1 + \mathbf{V}_I \cdot \nabla \phi_1 = 0$$

The system requires that the phases volume fractions add up to unity, $\phi_1 + \phi_2 = 1$. The closure relations for each phase are also given by

$$\rho_j e_j = \frac{p_j + \gamma_j \pi_j}{\gamma_j - 1}$$

Here $\gamma_j$ is the ratio of specific heats and $\pi_j$ is a constant. For the above EOS, the sound speed $c_j$ in each phase is given by

$$c_j = \sqrt{\gamma_j \frac{p_j + \pi_j}{\rho_j}}$$

In Section 4 of Dumbser *et al*. [20] a 2D smooth vortex test problem was designed for Baer-Nunziato flow. The problem is an analogue of the hydrodynamic vortex. In Balsara *et al*. [10] the problem was described in detail, so we do not describe it again here. For fifth, seventh and ninth orders, we also had to double the size of the computational domain in order to obtain results that were independent of the exponential fall off in the velocity of the vortex. Table I is based on the algorithm described in eqn. (15). It uses three WENO interpolation steps per RK stage and per dimension. In Table I we show the accuracy results on using WENO-AO interpolation and also Multiresolution WENO interpolation. We see that both algorithms reach their design accuracies and obtain completely comparable results. Table II is based on the algorithm described in eqn. (13). It uses two WENO interpolation steps per RK stage and per dimension. We see that it too meets the design accuracy, with the exception that at third order its accuracy degrades by some amount.

**Table I shows the accuracy of the 2D Baer-Nunziato Vortex problem using the AFD-WENO algorithm that is described in eqn. (15); the solid volume fraction is shown. The first half of the table shows the results from WENO-AO interpolation with non-linear limiting. The second half shows the same results when Multiresolution WENO interpolation is used.**

| WENO-AO-3 | $L_1$ Error | $L_1$ Accuracy | $L_{inf}$ Error | $L_{inf}$ Accuracy |
|---|---|---|---|---|
| $64^2$ | 4.88083E-04 | | 2.66278E-02 | |
| $128^2$ | 9.53635E-05 | 2.36 | 9.96971E-03 | 1.42 |
| $256^2$ | 1.66069E-05 | 2.52 | 3.51389E-03 | 1.50 |
| $512^2$ | 2.63532E-06 | 2.66 | 1.20590E-03 | 1.54 |
| WENO-AO-(5,3) | | | | |
| $32^2$ | 3.36706E-04 | | 2.87333E-02 | |
| $64^2$ | 3.99396E-05 | 3.08 | 3.38373E-03 | 3.09 |
| $128^2$ | 1.33499E-06 | 4.90 | 1.57672E-04 | 4.42 |
| $256^2$ | 4.26364E-08 | 4.97 | 5.36602E-06 | 4.88 |
| WENO-AO-(7,3) | | | | |
| $32^2$ | 2.52404E-04 | | 1.26094E-02 | |
| $64^2$ | 4.94039E-06 | 5.67 | 4.40379E-04 | 4.84 |
| $128^2$ | 4.35096E-08 | 6.83 | 5.68763E-06 | 6.27 |
| $256^2$ | 3.58042E-10 | 6.93 | 4.76953E-08 | 6.90 |
| WENO-AO-(9,3) | | | | |
| $32^2$ | 3.25316E-04 | | 1.97561E-02 | |
| $64^2$ | 8.49805E-07 | 8.58 | 7.63963E-05 | 8.01 |
| $96^2$ | 2.64728E-08 | 8.56 | 3.09700E-06 | 7.91 |
| $128^2$ | 2.22138E-09 | 8.61 | 2.70035E-07 | 8.48 |
| | | | | |
| Order 3 Multires WENO | | | | |
| $64^2$ | 2.72629E-04 | | 1.64969E-02 | |
| $128^2$ | 5.03731E-05 | 2.44 | 5.60625E-03 | 1.56 |
| $256^2$ | 8.83383E-06 | 2.51 | 1.87335E-03 | 1.58 |
| $512^2$ | 1.35812E-06 | 2.70 | 6.18262E-04 | 1.60 |
| Order 5 Multires WENO | | | | |
| $32^2$ | 6.03313E-04 | | 3.29950E-02 | |
| $64^2$ | 6.90399E-05 | 3.13 | 8.98365E-03 | 1.88 |
| $128^2$ | 8.84793E-06 | 2.96 | 1.52644E-03 | 2.56 |
| $256^2$ | 2.00259E-07 | 5.47 | 6.05237E-05 | 4.66 |
| Order 7 Multires WENO | | | | |
| $32^2$ | 3.22961E-04 | | 1.11033E-02 | |
| $64^2$ | 5.44778E-06 | 5.89 | 5.08236E-04 | 4.45 |
| $128^2$ | 6.54879E-08 | 6.38 | 1.08583E-05 | 5.55 |
| $256^2$ | 3.81689E-10 | 7.42 | 5.20716E-08 | 7.70 |
| Order 9 Multires WENO | | | | |
| $32^2$ | 2.34361E-04 | | 8.23576E-03 | |
| $64^2$ | 8.65817E-06 | 4.76 | 7.15367E-04 | 3.53 |
| $96^2$ | 2.74230E-08 | 14.19 | 3.25011E-06 | 13.30 |
| $128^2$ | 2.16289E-09 | 8.83 | 2.54077E-07 | 8.86 |

**Table II shows the accuracy of the 2D Baer-Nunziato Vortex problem using the AFD-WENO algorithm that is described in eqn. (13); the solid volume fraction is shown. The table shows the results from WENO-AO interpolation with non-linear limiting.**

| WENO-AO-3 | $L_1$ Error | $L_1$ Accuracy | $L_{inf}$ Error | $L_{inf}$ Accuracy |
|---|---|---|---|---|
| $64^2$ | 5.14071E-04 | | 2.66499E-02 | |
| $128^2$ | 1.09275E-04 | 2.23 | 1.00424E-02 | 1.41 |
| $256^2$ | 2.18926E-05 | 2.32 | 3.57657E-03 | 1.49 |
| $512^2$ | 4.26144E-06 | 2.36 | 1.24636E-03 | 1.52 |
| WENO-AO-(5,3) | | | | |
| $32^2$ | 3.91723E-04 | | 3.50496E-02 | |
| $64^2$ | 4.42581E-05 | 3.15 | 3.39514E-03 | 3.37 |
| $128^2$ | 1.77426E-06 | 4.64 | 1.61984E-04 | 4.39 |
| $256^2$ | 8.15939E-08 | 4.44 | 5.82518E-06 | 4.80 |
| WENO-AO-(7,3) | | | | |
| $32^2$ | 2.11416E-04 | | 1.14172E-02 | |
| $64^2$ | 5.53998E-06 | 5.25 | 4.18318E-04 | 4.77 |
| $128^2$ | 5.65413E-08 | 6.61 | 5.31257E-06 | 6.30 |
| $256^2$ | 6.74573E-10 | 6.39 | 4.96397E-08 | 6.74 |
| WENO-AO-(9,3) | | | | |
| $32^2$ | 2.37479E-04 | | 1.73189E-02 | |
| $64^2$ | 9.31242E-07 | 7.99 | 7.17748E-05 | 7.91 |
| $96^2$ | 3.02017E-08 | 8.46 | 2.84721E-06 | 7.96 |
| $128^2$ | 2.67208E-09 | 8.43 | 2.55134E-07 | 8.39 |

**V.b) Accuracy analysis for the two-dimensional Two-layer Shallow Water Equations**

In this Sub-section we focus on the two layer shallow water equations from Castro *et al*. [15]. The PDE system in 2D is given by defining $h_1$, $u_1$, $v_1$ as the height of the upper fluid, its x-velocity and its y-velocity respectively and by defining $h_2$, $u_2$, $v_2$ as the height of the lower fluid, its x-velocity and its y-velocity respectively. The bottom topography is denoted by "$b$" and "$g$" is the gravity. The ratio $\rho \equiv \rho_1/\rho_2$ denotes the ratio of the fluid densities. The total surface height, therefore, becomes $\eta = \eta_1 = b + h_1 + h_2$ and the surface elevation of the interior layer is denoted by $\eta_2 = b + h_2$. The conservation law is given by

$$\frac{\partial}{\partial t}\begin{pmatrix} h_1 \\ h_1 u_1 \\ h_1 v_1 \\ h_2 \\ h_2 u_2 \\ h_2 v_2 \\ b \end{pmatrix} + \frac{\partial}{\partial x}\begin{pmatrix} h_1 u_1 \\ h_1 u_1^2 + g h_1^2/2 \\ h_1 u_1 v_1 \\ h_2 u_2 \\ h_2 u_2^2 + g h_2^2/2 \\ h_2 u_2 v_2 \\ 0 \end{pmatrix} + \frac{\partial}{\partial y}\begin{pmatrix} h_1 v_1 \\ h_1 u_1 v_1 \\ h_1 v_1^2 + g h_1^2/2 \\ h_2 v_2 \\ h_2 u_2 v_2 \\ h_2 v_2^2 + g h_2^2/2 \\ 0 \end{pmatrix} + \begin{pmatrix} 0 \\ g h_1 \partial_x h_2 + g h_1 \partial_x b \\ 0 \\ 0 \\ \rho g h_2 \partial_x h_1 + g h_2 \partial_x b \\ 0 \\ 0 \end{pmatrix} + \begin{pmatrix} 0 \\ 0 \\ g h_1 \partial_y h_2 + g h_1 \partial_y b \\ 0 \\ 0 \\ \rho g h_2 \partial_y h_1 + g h_2 \partial_y b \\ 0 \end{pmatrix} = 0$$

All the linearly degenerate eigenvectors in both directions can be evaluated analytically. The non-linear eigenvectors and their eigenvalues have to be found numerically.

In Section 4.1 of Dumbser *et al*. [19] a 2D smooth vortex test problem was designed for the two-layer shallow water flow model. In Balsara *et al*. [10] the problem was described in detail, so we do not describe it again here. For fifth, seventh and ninth orders, we also had to double the size of the computational domain in order to obtain results that were independent of the exponential fall off in the velocity of the vortex. Table III is based on the algorithm described in eqn. (15). It uses three WENO interpolation steps per RK stage and per dimension. In Table III we show the accuracy results on using WENO-AO interpolation and also Multiresolution WENO interpolation. We see that both algorithms reach their design accuracies and obtain comparable results.

**Table III shows the accuracy of the 2D Two-layer Shallow Water Equations Vortex problem using the AFD-WENO algorithm that is described in eqn. (15); the variable $h_1$ is shown. The first half of the table shows the results from WENO-AO interpolation with non-linear limiting. The second half shows the same results when Multiresolution WENO interpolation is used.**

| WENO-AO-3 | $L_1$ Error | $L_1$ Accuracy | $L_{inf}$ Error | $L_{inf}$ Accuracy |
|---|---|---|---|---|
| $64^2$ | 3.44583E-03 | | 7.62937E-02 | |
| $128^2$ | 7.43327E-04 | 2.21 | 2.13210E-02 | 1.84 |
| $256^2$ | 1.17072E-04 | 2.67 | 5.04281E-03 | 2.08 |
| $512^2$ | 1.66279E-05 | 2.82 | 1.28914E-03 | 1.97 |
| WENO-AO-(5,3) | | | | |
| $64^2$ | 3.53368E-04 | | 2.43011E-02 | |
| $128^2$ | 3.51796E-05 | 3.33 | 9.02115E-03 | 1.43 |
| $256^2$ | 1.30811E-06 | 4.75 | 3.89311E-04 | 4.53 |
| $512^2$ | 4.20783E-08 | 4.96 | 1.28275E-05 | 4.92 |
| WENO-AO-(7,3) | | | | |
| $64^2$ | 2.26027E-04 | | 3.91491E-02 | |
| $128^2$ | 3.50401E-06 | 6.01 | 1.38103E-03 | 4.83 |
| $256^2$ | 3.20243E-08 | 6.77 | 1.35288E-05 | 6.67 |
| $512^2$ | 2.63013E-10 | 6.93 | 1.11319E-07 | 6.93 |
| WENO-AO-(9,3) | | | | |
| $32^2$ | 6.27801E-03 | | 4.50278E-01 | |
| $64^2$ | 2.19548E-04 | 4.84 | 2.12822E-02 | 4.40 |
| $128^2$ | 4.42786E-07 | 8.95 | 1.99559E-04 | 6.74 |
| $256^2$ | 1.07759E-09 | 8.68 | 5.09173E-07 | 8.61 |
| | | | | |
| Order 3 Multires WENO | | | | |
| $64^2$ | 2.01362E-03 | | 5.07692E-02 | |
| $128^2$ | 3.77690E-04 | 2.41 | 1.01869E-02 | 2.32 |
| $256^2$ | 5.47456E-05 | 2.79 | 2.66140E-03 | 1.94 |
| $512^2$ | 7.65931E-06 | 2.84 | 7.11564E-04 | 1.90 |
| Order 5 Multires WENO | | | | |
| $64^2$ | 7.38835E-04 | | 8.02979E-02 | |
| $128^2$ | 8.01588E-05 | 3.20 | 1.25310E-02 | 2.68 |
| $256^2$ | 3.34184E-06 | 4.58 | 8.72471E-04 | 3.84 |
| $512^2$ | 6.08324E-08 | 5.78 | 2.78790E-05 | 4.97 |
| Order 7 Multires WENO | | | | |
| $64^2$ | 4.21898E-04 | | 4.34015E-02 | |
| $128^2$ | 1.37073E-05 | 4.94 | 3.32110E-03 | 3.71 |
| $256^2$ | 1.15696E-07 | 6.89 | 7.37685E-05 | 5.49 |
| $512^2$ | 6.28799E-10 | 7.52 | 4.49976E-07 | 7.36 |
| Order 9 Multires WENO | | | | |
| $32^2$ | 6.27152E-03 | | 4.41670E-01 | |
| $64^2$ | 7.00717E-04 | 3.16 | 1.07767E-01 | 2.04 |
| $128^2$ | 6.08838E-06 | 6.85 | 1.52932E-03 | 6.14 |
| $256^2$ | 9.38793E-09 | 9.34 | 5.70332E-06 | 8.07 |

**V.c) Accuracy analysis for the two-dimensional Multiphase Debris Flow Model of Pitman and Le**

In this Sub-section we consider the multiphase debris flow model of Pitman and Le [37]. We use the formulation of Pelanti *et al*. [36] instead of the original formulation by Pitman and Le [37]. The PDE system in 2D is given by defining $h_s$, $u_s$, $v_s$ as the solid height, the solid x-velocity and the solid y-velocity respectively and by defining $h_f$, $u_f$, $v_f$ as the fluid height, the fluid x-velocity and the solid y-velocity respectively. We also have $h_f = (1-\phi)h$ where $\phi$ is the solid volume fraction and $h$ is the total height. The variable "$b$" refers to the bottom topography and is kept constant in all our test problems. The variable "g" refers to the gravitational acceleration and "$\rho \equiv \rho_f / \rho_s$" refers to the density ratio of the fluid and the solid. The time evolutionary equations for this model are:-

$$\frac{\partial}{\partial t}\begin{pmatrix} h_s \\ h_s u_s \\ h_s v_s \\ h_f \\ h_f u_f \\ h_f v_f \\ b \end{pmatrix} + \frac{\partial}{\partial x}\begin{pmatrix} h_s u_s \\ h_s u_s^2 + g h_s^2 \rho / 2 + g(1-\rho) h_s (h_f + h_s)/2 \\ h_s u_s v_s \\ h_f u_f \\ h_f u_f^2 + g h_f^2 / 2 \\ h_f u_f v_f \\ 0 \end{pmatrix}$$

$$+ \frac{\partial}{\partial y}\begin{pmatrix} h_s v_s \\ h_s u_s v_s \\ h_s v_s^2 + g h_s^2 \rho / 2 + g(1-\rho) h_s (h_f + h_s)/2 \\ h_f v_f \\ h_f u_f v_f \\ h_f v_f^2 + g h_f^2 / 2 \\ 0 \end{pmatrix} + \begin{pmatrix} 0 \\ g h_s \rho \, \partial_x h_f + g \, h_s \partial_x b \\ 0 \\ 0 \\ g h_f \partial_x h_s + g \, h_f \partial_x b \\ 0 \\ 0 \end{pmatrix} + \begin{pmatrix} 0 \\ 0 \\ g h_s \rho \, \partial_y h_f + g \, h_s \partial_y b \\ 0 \\ 0 \\ g h_f \partial_y h_s + g \, h_f \partial_y b \\ 0 \end{pmatrix} = 0$$

All the linearly degenerate eigenvectors in both directions can be evaluated analytically. The non-linear eigenvectors and their eigenvalues have to be found numerically.

In Section 4.2 of Dumbser *et al*. [19], a two-dimensional smooth vortex test problem was designed for the debris flow model. In Balsara *et al*. [10] the problem was described in detail, so we do not describe it again here. As previously, to minimize the effect of small jumps in the velocity field at the periodic boundaries, we double the computational domain and stopping time for the fifth, seventh and ninth order schemes. In Table IV we show the accuracy results on using WENO-AO interpolation and also Multiresolution WENO interpolation. The algorithm uses three WENO interpolation steps per RK stage and per dimension. As before, we see that the scheme reaches its design accuracy.

**Table IV shows the accuracy of the 2D Multiphase Debris Flow Model Vortex problem using the AFD-WENO algorithm that is described in eqn. (15); the variable $h_s$ is shown. The first half of the table shows the results from WENO-AO interpolation with non-linear limiting. The second half shows the same results when Multiresolution WENO interpolation is used.**

| WENO-AO-3 | L$_1$ Error | L$_1$ Accuracy | L$_{inf}$ Error | L$_{inf}$ Accuracy |
|---|---|---|---|---|
| 128$^2$ | 6.80623E-04 | | 2.42650E-02 | |
| 256$^2$ | 1.40708E-04 | 2.27 | 6.95720E-03 | 1.80 |
| $512^2$ | 2.68923E-05 | 2.39 | 1.77262E-03 | 1.97 |
| $1024^2$ | 4.30445E-06 | 2.64 | 5.52426E-04 | 1.68 |
| WENO-AO-(5,3) | | | | |
| $64^2$ | 1.27772E-03 | | 1.13181E-01 | |
| $128^2$ | 1.00327E-04 | 3.67 | 1.81570E-02 | 2.64 |
| $256^2$ | 3.63859E-06 | 4.79 | 6.86942E-04 | 4.72 |
| $512^2$ | 1.16263E-07 | 4.97 | 2.24444E-05 | 4.94 |
| WENO-AO-(7,3) | | | | |
| $64^2$ | 8.29424E-04 | | 6.91461E-02 | |
| $128^2$ | 1.02160E-05 | 6.34 | 2.36519E-03 | 4.87 |
| $256^2$ | 9.13526E-08 | 6.81 | 2.33227E-05 | 6.66 |
| $512^2$ | 7.41144E-10 | 6.95 | 1.94674E-07 | 6.90 |
| WENO-AO-(9,3) | | | | |
| $32^2$ | 6.15229E-03 | | 4.40454E-01 | |
| $64^2$ | 7.89605E-04 | 2.96 | 5.23724E-02 | 3.07 |
| $128^2$ | 1.36070E-06 | 9.18 | 3.38921E-04 | 7.27 |
| $256^2$ | 3.28992E-09 | 8.69 | 9.66972E-07 | 8.45 |
| | | | | |
| Order 3 Multires WENO | | | | |
| $64^2$ | 4.16105E-04 | | 1.57617E-02 | |
| $128^2$ | 7.17671E-05 | 2.54 | 3.82240E-03 | 2.04 |
| $256^2$ | 1.22607E-05 | 2.55 | 8.93516E-04 | 2.10 |
| $512^2$ | 1.89426E-06 | 2.69 | 2.81595E-04 | 1.67 |
| Order 5 Multires WENO | | | | |
| $64^2$ | 1.26472E-03 | | 7.18582E-02 | |
| $128^2$ | 1.68468E-04 | 2.91 | 2.12744E-02 | 1.76 |
| $256^2$ | 7.52555E-06 | 4.48 | 1.09312E-03 | 4.28 |
| $512^2$ | 2.05692E-07 | 5.19 | 4.15323E-05 | 4.72 |
| Order 7 Multires WENO | | | | |
| $64^2$ | 1.06838E-03 | | 1.10762E-01 | |
| $128^2$ | 2.89278E-05 | 5.21 | 4.47504E-03 | 4.63 |
| $256^2$ | 2.74816E-07 | 6.72 | 7.93913E-05 | 5.82 |
| $512^2$ | 1.98876E-09 | 7.11 | 8.11574E-07 | 6.61 |
| Order 9 Multires WENO | | | | |
| $32^2$ | 7.92609E-03 | | 5.18075E-01 | |
| $64^2$ | 1.24043E-03 | 2.68 | 1.40228E-01 | 1.89 |
| $128^2$ | 1.31660E-05 | 6.56 | 1.86457E-03 | 6.23 |
| $256^2$ | 2.33957E-08 | 9.14 | 6.60566E-06 | 8.14 |

**V.d) Speed Comparisons**

At this point it is worthwhile to make estimates of the computational cost of various well-known finite difference WENO algorithms. We do that in this paragraph and in the next two paragraphs that follow. The FD-WENO scheme described in Balsara *et al*. [10] entails one WENO reconstruction and one call to the Riemann solver. Since the Riemann solvers tend to be light-weight, the majority of the computational cost is in the WENO reconstruction, which absolutely must be done in characteristic space if the problem has any significant discontinuities. The FD-WENO scheme described in Balsara *et al*. [10] can handle hyperbolic systems with non-conservative products, but it does not guarantee the existence of a flux conservative form when a conservation law is used. Since the classical FD-WENO of Shu and Osher [41] also entailed two WENO reconstruction steps in characteristic variables, the FD-WENO scheme from Balsara *et al*. [10] is expected to cost about half as much as the classical FD-WENO of Shu and Osher [41].

The AFD-WENO scheme developed in Balsara *et al*. [11] for conservation laws requires two WENO interpolation steps and one call to a (light-weight) Riemann solver. The first of these WENO interpolation steps must be done in characteristic variables if the physical problem has discontinuities. The cost of the second WENO interpolation step may be reduced by doing it in a component by component fashion; and this is the choice we make in this Sub-section. Since the classical FD-WENO of Shu and Osher [41] entailed two WENO reconstruction steps in characteristic variables, the method in Balsara *et al*. [11] is expected to cost less. Both methods can only treat conservation laws, however the AFD-WENO method is much more versatile in terms of the Riemann solvers it can handle and in terms of its ability to work on curvilinear meshes.

The AFD-WENO scheme developed in this paper is more versatile because it can handle hyperbolic systems with non-conservative products. It does this while guaranteeing that variables that are still in conservation form are updated with a strictly flux conservative update. The method in this paper requires three WENO interpolation steps and one call to a (light-weight) Riemann solver. The first of these WENO interpolation steps must be done in characteristic variables if the physical problem has discontinuities. The cost of the second and third WENO interpolation steps may be substantially reduced by doing them in a component by component fashion; and this is the choice that we make here. It could, therefore, be slightly more expensive than the classical FD-WENO of Shu and Osher [41], but it is much more versatile. It too can handle different types of Riemann solvers (as long as they can be written in fluctuation form) and it too can work on curvilinear meshes.

Lastly, and perhaps most importantly, AFD-WENO schemes have been used in the past with central differences for the higher order derivatives. The interested reader might wish to know whether using central differences for the higher order derivatives would give a scheme that is substantially faster than the current schemes that uses WENO interpolation for their second and third steps? To answer that question, we have also gathered timing statistics for AFD-WENO schemes when only central differences are used for the higher order derivatives. This answers the question of efficiency. One major purpose of scientific computation is to do new problems for which the solution is initially not known. The other major purpose of scientific computation is to use well-tested methods on the newer classes of PDEs that arise in different fields of science and engineering. Therefore, one wants numerical methods for PDEs that have been verified to work on a large range of well-known problems stemming from different types of well-known PDEs. One can, therefore, ask a further question:- For general problems with discontinuities, is it effective to use just an AFD-WENO scheme that uses central differences for its higher order derivatives? We will show that it is not effective in the subsequent Section.

In light of the expectations that were laid out in the above four paragraphs, it is therefore useful to document the actual speeds of all the above-mentioned finite difference WENO formulations on a problem where they are all applicable. For that reason, we took the very popular one-dimensional Sod shock tube problem for Euler flow and ran it from start to finish on a 1000 zone mesh that is evolved for 200 timesteps. Third order accurate in time SSP-RK timestepping was used in all cases. To have a fair comparison across the competing WENO formulations, the LLF Riemann solver was used in all cases. We did this with a code that was stripped of all input and output functions so that we were able to document just the raw number of zones updated per second for this problem for the various algorithms. Also note that characteristic-based interpolation was only used in the first step for the algorithm described in this paper as well as the algorithm described in Balsara *et al*. [11]. This is justified because the additional WENO interpolation steps do not need characteristic decomposition. A single core of a modern Xeon Gold 6248R processor running at 3 GHz was used for all the runs.

In the first four paragraphs of this Section we had catalogued our expectations for the speeds of the various finite difference WENO schemes being considered. The data from Table V shows that the WENO scheme from Balsara *et al*. [10] is the fastest, which bears out our expectations. Classical FD-WENO is somewhat slower than AFD-WENO from Balsara *et al*. [11] because the latter scheme only uses one characteristic decomposition whereas the former scheme uses two characteristic decompositions. The AFD-WENO scheme from this paper is somewhat slower, but not by much; because it only uses one characteristic decomposition along with two interpolations that are done component by component. However, it provides more capabilities in that it can handle hyperbolic systems with conservative terms as well as non-conservative products. The AFD-WENO scheme from eqn. (13) of this paper, which is expected to be faster than the AFD-WENO described in eqn. (15) of this paper, is also shown. The AFD-WENO scheme that uses only central derivatives is nominally comparable in speed to the AFD-WENO scheme from eqn. (13) of this paper. This is understandable because the major cost in high order WENO interpolation is associated with the largest stencil, and the central differences also have to use the same large stencil that is used in WENO interpolation!

Table V also allows us to cross-compare the speeds of the same algorithm at different orders of spatial accuracy. Our WENO implementations are based on a strip-mining philosophy, where a one-dimensional strip of data is extracted into long one-dimensional arrays and the CPU is asked to work on those strips of data (Colella and Woodward [17], Woodward and Colella [47]). This philosophy is very beneficial for modern CPUs with their large caches because the one-dimensional strips can then be loaded into cache very efficiently and automatically by the CPU. The CPU can then do a lot of work on each strip before releasing it back to the computer's main memory. By scanning the rows of Table V, we see that from third to ninth order, the increase in cost is very modest. This gives us another insight into the nature of finite difference WENO schemes on modern CPUs with large caches. All the finite difference WENO schemes in Table V incur some fixed costs, such as constructing eigenvectors or invoking Riemann solvers or copying to and from main memory. When the WENO reconstruction or interpolation has been perfectly optimized, the larger stencils associated with higher order WENO schemes do not seem to significantly reduce the speed (i.e., the number of zones updated per second). Instead, it is now the fixed costs that begin to dominate the overall cost of the algorithm. The insight that we gain from scanning the rows of Table V is that the larger number of floating point operations associated with larger stencils are almost free of charge if the finite difference WENO code is designed to be cache friendly. Our well-implemented FD-WENO and AFD-WENO schemes come pretty close to the promised land of high performance computing where float point operations are almost free and the dominant costs are the other fixed costs in the implementation!

**Table V shows the speed comparison, expressed as number of zones updated per second, for the Sod shock problem on a one-dimensional 1000 zone mesh that is evolved for 200 timesteps. A single core of a modern Xeon Gold 6248R processor running at 3 GHz was used for all the runs. The classical FD-WENO scheme (Shu and Osher [41], Jiang and Shu [29]) is denoted "Classical FD-WENO"; the WENO scheme from Balsara *et al*. [10] is denoted "WENO B23a"; the AFD-WENO scheme from Balsara *et al*. [11] is denoted "AFD-WENO B23b"; the AFD-WENO scheme from eqn. (15) of this paper is denoted as "AFD-WENO This Paper". We also show the AFD-WENO algorithm from eqn. (13) of this paper, which we denote as "AFD-WENO 13". Lastly, we also show the speed of AFD-WENO algorithm from eqn. (15) when only central derivatives are used for the higher order derivatives, and this is denoted as "AFD-WENO Central Derivatives".**

|  | 3$^{rd}$ Order | 5$^{th}$ Order | 7$^{th}$ Order | 9$^{th}$ Order |
|---|---|---|---|---|
| **Classical FD-WENO** | 80260 zones/s | 72672 zones/s | 65374 zones/s | 60685 zones/s |
| **WENO B23a** | 98644 zones/s | 95720 zones/s | 92090 zones/s | 88684 zones/s |
| **AFD-WENO B23b** | 91786 zones/s | 84623 zones/s | 78093 zones/s | 76741 zones/s |
| **AFD-WENO This Paper** | 85396 zones/s | 77907 zones/s | 69871 zones/s | 66751 zones/s |
| **AFD-WENO 13** | 92483 zones/s | 83159 zones/s | 77105 zones/s | 72850 zones/s |
| **AFD-WENO Central Derivatives** | 89302 zones/s | 85474 zones/s | 78654 zones/s | 70042 zones/s |

**VI) One Dimensional Test Problems**

Here we present the same one-dimensional test problems as from Balsara *et al.* [10] for the Baer-Nunziato system, the two layer shallow water system and the multiphase debris flow model. In Sub-section VI.a we focus on the Baer-Nunziato model of compressible multi-phase flows. In Sub-section VI.b we focus on the two-layer shallow water equations. In Sub-section VI.c we present results from the multiphase debris flow model of Pitman and Le. For all the simulations presented in this Section we used a CFL of 0.8 with a third order SSP-RK scheme.

**VI.a) One-dimensional Test Problems for the Baer-Nunziato Model for Compressible Multi-Phase Flows**

For the Baer-Nunziato model described in Section V.a, we first consider the Abgrall problem from Dumbser *et al*. [20]. The same problem was also used in Balsara *et al.* [10]. The Abgrall problem, see Saurel and Abgrall [39], is based on the idea that a mixture of the two phases moving in 1D with uniform velocity and pressure should be able to continue moving in this fashion without generating any wiggles in the velocity or pressure. We present such a test problem here on a 200 zone mesh that spans the domain $[-0.5, 0.5]$. In the solid phase we set a stiffened EOS with $\gamma_1 = 3$ and $\pi_1 = 100$, whereas in the gas phase we set $\gamma_2 = 1.4$ and $\pi_2 = 0$. The pressures and longitudinal velocities in both phases are set to unity. The problem is initialized with $\rho_{1,L} = 800$, $\rho_{2,L} = 2$, $\phi_{1,L} = 0.99$ on the left of $x = 0$ and $\rho_{1,R} = 1000$, $\rho_{2,R} = 1$, $\phi_{1,L} = 0.01$ to its right. As the Abgrall problem contains strong shocks in the volume fractions, therefore we use the flattening algorithm (as described in Appendix B) to detect the shock locations and there by suppressing the higher-order terms in the simulation. For zone *i,* if we have $\eta_i > 10^{-12}$ then we set higher order terms in eqn. (15), i.e. last two curly bracket terms, to zero. Using the available flattener, we run the problem to a final time of $t = 0.25$ using the fifth order accurate LLF-based AFD-WENO scheme. The result for the solid volume fraction is shown in Fig. 4a, and the pressure and velocity profiles are shown in Fig. 4b. We see that the pressure and velocity profiles are absolutely flat, showing that the scheme is able to capture the solution without generating any oscillations in the velocity and pressure profiles. The seventh and ninth order AFD-WENO schemes also show identical results, and therefore they are not presented here.

Next we consider a set of six one-dimensional Riemann problems from Dumbser *et al*. [20]. The same set of problems were also used in Balsara *et al.* [10]. Table VI catalogues the parameters for these six Riemann problems. The results for these six problems using the AFD-WENO scheme are shown in Fig. 5 to Fig. 10. All the Baer-Nunziato Riemann problems were run on a 200 zone mesh.

**Table VI gives the left and right initial states, computational domain $[x_L, x_R] = [-0.5, 0.5]$ and final times $(t_{end})$ of the six Riemann problems for Baer-Nunziato model.**

| Test Case: | $\rho_1$ | $u_1$ | $p_1$ | $\rho_2$ | $u_2$ | $p_2$ | $\phi_1$ | $t_{end}$ |
|---|---|---|---|---|---|---|---|---|
| **RP1** | $\gamma_1 = 1.4,$ | $\pi_1 = 0,$ | $\gamma_2 = 1.4,$ | $\pi_2 = 0$ | | | | |
| Left: | 1.0 | 0.0 | 1.0 | 0.5 | 0.0 | 1.0 | 0.4 | 0.1 |
| Right: | 2.0 | 0.0 | 2.0 | 1.5 | 0.0 | 2.0 | 0.8 | |
| **RP2** | $\gamma_1 = 3.0,$ | $\pi_1 = 100,$ | $\gamma_2 = 1.4,$ | $\pi_2 = 0$ | | | | |
| Left: | 800.0 | 0.0 | 500.0 | 1.5 | 0.0 | 2.0 | 0.4 | 0.1 |
| Right: | 1000.0 | 0.0 | 600.0 | 1.0 | 0.0 | 1.0 | 0.3 | |
| **RP3** | $\gamma_1 = 1.4,$ | $\pi_1 = 0,$ | $\gamma_2 = 1.4,$ | $\pi_2 = 0$ | | | | |
| Left: | 1.0 | 0.9 | 2.5 | 1.0 | 0.0 | 1.0 | 0.9 | 0.1 |
| Right: | 1.0 | 0.0 | 1.0 | 1.2 | 1.0 | 2.0 | 0.2 | |
| **RP4** | $\gamma_1 = 3.0,$ | $\pi_1 = 3400,$ | $\gamma_2 = 1.35,$ | $\pi_2 = 0$ | | | | |
| Left: | 1900.0 | 0.0 | 10.0 | 2.0 | 0.0 | 3.0 | 0.2 | 0.15 |
| Right: | 1950.0 | 0.0 | 1000.0 | 1.0 | 0.0 | 1.0 | 0.9 | |
| **RP5** | $\gamma_1 = 1.4,$ | $\pi_1 = 0,$ | $\gamma_2 = 1.4,$ | $\pi_2 = 0$ | | | | |
| Left: | 1.0 | 0.0 | 1.0 | 0.2 | 0.0 | 0.3 | 0.8 | 0.20 |
| Right: | 1.0 | 0.0 | 1.0 | 1.0 | 0.0 | 1.0 | 0.3 | |
| **RP6** | $\gamma_1 = 1.4,$ | $\pi_1 = 0,$ | $\gamma_2 = 1.4,$ | $\pi_2 = 0$ | | | | |
| Left: | 0.2068 | 1.4166 | 0.0416 | 0.5806 | 1.5833 | 1.375 | 0.1 | 0.1 |
| Right: | 2.2263 | 0.9366 | 6.0 | 0.4890 | -0.70138 | 0.986 | 0.2 | |

Figs. 5a and 5b show the solid and gas densities for the RP1 test when a fifth order HLL-based AFD-WENO scheme was used. Figs. 6a and 6b show the solid and gas densities for the RP2 test when a seventh order HLL-based AFD-WENO scheme was used. Figs. 7a and 7b show the solid and gas densities for the RP3 test when a ninth order HLL-based AFD-WENO scheme was used. Figs. 8a and 8b show the solid and gas densities for the RP4 test when a fifth order LLF-based AFD-WENO scheme was used. Figs. 9a and 9b show the solid and gas densities for the RP5 test when a seventh order LLF-based AFD-WENO scheme was used. For this problem (RP5), we notice that the seventh and ninth-order schemes tend to exaggerate the small oscillations near shocks if a flattener is not used. For this reason, we have invoked the flattener (given in Appendix B) to simulate this problem. Figs. 10a and 10b show the solid and gas densities for the RP6 test when a ninth order LLF-based AFD-WENO scheme was used. The other higher order LLF-based and HLL-based AFD-WENO schemes also show identical results, and therefore they are not presented here. In each of the six Riemann problems a reference solution has been shown in solid lines. The reference solution was obtained using the third order LLF-based AFD-WENO scheme on a mesh with 2000 zones. We can observe that all the results obtained using AFD-WENO schemes provide a clearer representation of the interfaces.

## VI.b) One-dimensional Test Problems for Two-layer Shallow Water Equations

For the Two-layer Shallow Water equations described in Section V.b, a set of three one-dimensional Riemann problems has been presented in Dumbser and Balsara [24], Dumbser *et al*. [19] and Castro *et al*. [14]. The same set of problems was also used in Balsara *et al*. [10]. Table VII gives the relevant information for these Riemann problems. All the two-layer shallow water Riemann problems were run on a 200 zone mesh with $\rho = 0.8$ and $g = 9.8$.

**Table VII gives the left and right initial states, computational domain $[x_L, x_R]$, location of the discontinuity $(x_c)$ and final times $(t_{end})$ of the three Riemann problems for the Two-Layer Shallow Water Model.**

| Test Case: | | $h_1$ | $u_1$ | $v_1$ | $h_2$ | $u_2$ | $v_2$ | $b$ | $x_L$ | $x_R$ | $x_c$ | $t_{end}$ |
|---|---|---|---|---|---|---|---|---|---|---|---|---|
| **RP1** | Left: | 0.5 | 0.0 | 0.5 | 0.8 | 0.0 | -0.2 | 0.2 | -5.0 | 5.0 | 0.0 | 1.0 |
|  | Right: | 0.5 | 0.0 | -0.5 | 0.2 | 0.0 | 0.2 | 0.8 |  |  |  |  |
| **RP2** | Left: | 0.4 | 0.0 | 0.0 | 0.6 | 0.0 | 0.0 | 0.0 | -5.0 | 5.0 | 0.0 | 1.25 |
|  | Right: | 0.6 | 0.0 | 0.0 | 0.4 | 0.0 | 0.0 | 0.0 |  |  |  |  |
| **RP3** | Left: | 1.0 | 0.0 | 0.0 | 1.0 | 0.0 | 0.0 | 0.0 | -5.0 | 5.0 | 0.0 | 1.0 |
|  | Right: | 0.5 | 0.0 | 0.0 | 0.5 | 0.0 | 0.0 | 0.5 |  |  |  |  |

The first Riemann problem from Table VII illustrates the preservation of stationary jump discontinuities in the linearly degenerate intermediate fields, and the results using the fifth order HLLI-based AFD-WENO scheme are shown in the upper panel of Fig. 11. The results from the seventh and ninth order AFD-WENO schemes are identical, and therefore they are not shown here. Fig. 11 shows the ability of our high order schemes to preserve stationary discontinuities in linearly degenerate intermediate fields. Fig. 12a shows the results from the second Riemann problem in Table VII computed using the seventh order accurate HLLI-based AFD-WENO scheme. The fifth and ninth order schemes show the identical results, and therefore they are not shown here. Fig. 12b shows results from the third Riemann problem in Table VII using the ninth order accurate HLLI-based AFD-WENO scheme. The fifth and seventh order schemes show the identical results, and therefore they are not shown here. A third order AFD-WENO scheme was also run on a mesh with 2000 zones in order to generate a reference solution, and the reference solution is shown with solid lines in Fig. 12. We see that the results in Fig. 12 closely match the reference solution.

Let us now address the topic of effectiveness of an AFD-WENO scheme that uses only central differences for its higher order derivatives. RP1 was run when only central differences were used for the higher order derivatives and the results are shown in the lower panel of Fig. 11. Similarly, RP3 was run when only central differences were used for the higher order derivatives and the results are shown in Fig. 12c. In each instance, we see that there are unphysical oscillations at the discontinuities. The first key purpose of scientific computation is to do new problems which have not been attempted before and for which the solution is not known. The second key purpose is that the methods that have been developed should be applicable to newer classes of PDEs that may arise in the future in science and engineering. This requires verification of the algorithms on large classes of known solutions emerging from different types of PDEs. Before an algorithm can give reliable results for newer classes of problems, for which the solution is unknown, it is essential that it should reproduce all known problems. The results in Figs. 11 and 12 tell us that an AFD-WENO scheme that only uses central differences for its higher order derivatives would not constitute a properly verified algorithm.

**VI.c) One-dimensional Test Problems for Multiphase Debris Flow Model of Pitman and Le**

In this Sub-section we present one-dimensional test problems for the multiphase debris flow model of Pitman and Le [37] in the formulation of Pelanti *et al*. [36]. The equations for the multiphase debris flow model are described in Sub-section V.c. The test problems we present are drawn from Pelanti *et al*. [36], Dumbser *et al*. [19] and Rhebergen *et al*. [38] and they have also been integrated in Dumbser and Balsara [24].

Table VIII describes the three Riemann problems that we present here; they are drawn from Dumbser and Balsara [24]. In all three Riemann problems we use $\rho = 0.5$ and $g = 9.8$. In all cases we use a one-dimensional mesh with 200 zones. The computational domain and the final stopping time is also shown in Table VIII. We use the fifth order accurate HLLI-based AFD-WENO scheme for the first problem from Table VIII. The problem corresponds to the superposition of three stationary jump discontinuities in the linearly degenerate fields. The problem was run with the anti-diffusive fluxes turned on for all the linearly degenerate waves. Figs. 13a-e show that all the jumps in the linearly degenerate stationary wave families are exactly preserved on the mesh when the full algorithm is used. The seventh and ninth order AFD-WENO schemes also show identical results, and therefore they are not shown here. Figs. 13f-j show the same results from an AFD-WENO scheme that only used central differences for the higher order derivatives. We can clearly see that there are spurious oscillations; indicating that such a scheme would not pass all the verification test problems. Fig. 14 shows the second Riemann problem from Table VIII obtained using the seventh order accurate HLLI-based AFD-WENO scheme. The fifth and ninth order AFD-WENO schemes also show identical results, and therefore they are not shown here. Fig. 15 shows the third Riemann problem from Table VIII obtained using the ninth order accurate HLLI-based AFD-WENO scheme. The fifth and seventh order AFD-WENO schemes also show identical results, and therefore they are not shown here. A third order AFD-WENO scheme was also run on a mesh with 2000 zones in order to generate a reference solution, and the reference solution is shown with solid lines in Fig. 14 and Fig. 15. We can observe that the results in Fig. 14 and Fig. 15 closely resemble the reference solution.

**Table VIII gives the left and right initial states, computational domain $[x_L, x_R]$, location of the discontinuity $(x_c)$ and final times $(t_{end})$ of the three Riemann problems for the Multiphase Debris Flow Model.**

| Test Case: | | $h_s$ | $u_s$ | $v_s$ | $h_f$ | $u_f$ | $v_f$ | $b$ | $x_L$ | $x_R$ | $x_c$ | $t_{end}$ |
|---|---|---|---|---|---|---|---|---|---|---|---|---|
| **RP1** | Left: | 1.5 | 0.0 | 0.2 | 0.5 | 0.0 | -0.5 | 0.0 | -5.0 | 5.0 | 0.0 | 1.0 |
| | Right: | 1.125 | 0.0 | -0.2 | 0.375 | 0.0 | 0.5 | 0.5 | | | | |
| **RP2** | Left: | 2.1 | 0.0 | 0.0 | 0.9 | 0.0 | 0.0 | 0.0 | -5.0 | 5.0 | 0.0 | 0.5 |
| | Right: | 0.8 | 0.0 | 0.0 | 1.2 | 0.0 | 0.0 | 0.0 | | | | |
| **RP3** | Left: | 2.1 | -1.4 | 0.0 | 0.9 | 0.3 | 0.0 | 0.0 | -5.0 | 5.0 | 0.0 | 0.5 |
| | Right: | 0.8 | -0.9 | 0.0 | 1.2 | 0.1 | 0.0 | 0.0 | | | | |

**VII) Multidimensional Test Problems**

In this Section, we present the same two-dimensional test problems as from Balsara *et al.* [10] for the Baer-Nunziato system, the two-layer shallow water system, and the multiphase debris flow model. In Sub-section VII.a, we focus on the Baer-Nunziato model of compressible multi-phase flows. In Sub-section VII.b, we focus on the two-layer shallow water equations. In Sub-section VII.c, we present results from the multiphase debris flow model of Pitman and Le. For all the simulations presented in this Section, we used a CFL of 0.4 with a third-order SSP-RK scheme.

**VII.a) Two-dimensional Test Problems for the Baer-Nunziato Model for Compressible Multi-Phase Flows**

For the Baer-Nunziato model described in Section V.a, we present two multi-dimensional problems. The first problem is the shock-bubble interaction problem from Dumbser *et al*. [22]. This problem was also used in Balsara *et al.* [10]. The second problem is an analog of Euler's Shock-Vortex interaction problem described in Balsara and Shu [4]. Such a problem for the Baer-Nunziato model has been presented in Section 9.2 of Balsara et al. [10]. The detailed setup for both problems is given in Balsara *et al.* [10]; therefore, we do not describe the setup here.

The first problem (shock-bubble interaction problem) consists of a planar, right-going shock propagating into an ambient medium that contains a bubble with a radius of 0.25. The problem is set up on a domain that spans $[-0.5, 3.0] \times [-0.75, 0.75]$ and is run to a final time of 0.025 using the third, fifth and seventh order accurate HLL-based AFD WENO-AO schemes on a mesh with 700×300 zones. Fig. 16a shows the resulting solid density for the third order accurate WENO-AO-3 scheme, Fig. 16b shows the resulting solid density for the fifth order accurate WENO-AO-(5,3) scheme and Fig. 16c shows the resulting solid density for the seventh order accurate WENO-AO-(7,5,3) scheme. The ninth order scheme shows similar result to seventh order results, therefore, it is not shown here. We observe that the density profile obtained using the third order scheme lacks sharp resolution. In contrast, the fifth and seventh order schemes effectively resolve all flow structures, thus demonstrating the significance of higher-order schemes.

Next, we consider the Shock-Vortex interaction problem for the Baer-Nunziato model. A vortex with the mean velocity of the pre-shocked region propagates diagonally into a stationary shock, and as time increases, we see various rearranged structures in the vortex profile. The computational domain for this problem spans $[-0.5, 1.5] \times [-0.5, 1.5]$, and the problem was run to a final time of 0.84 using the fifth-order accurate HLL-based AFD WENO-AO-(5,3) scheme on a mesh with 600×600 zones. Figs. 17a, b, and c show solid volume fractions at times t=0.0, 0.23, 0.84. The intermediate time of 0.23 corresponds to a time when the vortex has propagated halfway through the shock. Note that the shock does not reveal itself in this variable. In Figs. 17d, 17e, and 17f we show solid x-velocity at times t=0.0, 0.23, 0.84. The shock is visible in this variable. We see that as the vortex propagates through the shock, the structure of the vortex is rearranged by the shock. However, the vortex is not fully destroyed by the shock. As a result, the vortex sheds some of its angular momentum after it has passed through the shock in an attempt to rearrange its structure. This is revealed by the spiral arms that are shed by the vortex in Figs. 17c and 17f. However, the figure shows us that vortices are very robust flow structures because they carry the angular momentum of the fluid. Because angular momentum is conserved by the Baer-Nunziato flow, the vortices will not be fully destroyed as they pass through shocks of modest strength. This observation about vortices not being fully destroyed by shocks was previously established for Euler flow in Balsara and Shu [4]. However, Fig. 17 provides further evidence that this insight also holds true for the Baer-Nunziato flow. The seventh and ninth order AFD-WENO schemes also show identical results, therefore, they are not shown here.

**VII.b) Two-dimensional Test Problems for the Two-Layer Shallow water equations**

For the Two-Layer Shallow water equations described in Section V.b, we present two multi-dimensional problems. The first problem is the shock-bubble interaction, and the second problem is Shock-Vortex interaction problem. Both problems are presented in great details in Sections 9.3 and 9.4 of Balsara *et al.* [10]. As the detailed setup is already given in Balsara *et al.* [10], we do not describe the setup here.

The first problem (shock-"bubble" interaction problem) consists of a planar, right-going shock propagating into an ambient medium which contains a bubble-like structure. The problem is set up on a domain that spans $[-0.5, 3.0] \times [-0.75, 0.75]$ and is run to a final time of 0.3 using the third, fifth and seventh order HLL-based AFD-WENO-AO schemes on a mesh with 700×300 zones. Fig. 18a shows the resulting height of the lower fluid for the third order accurate WENO-AO-3 scheme, Fig. 18b shows the resulting height of the lower fluid for the fifth order accurate WENO-AO-(5,3) scheme and Fig. 18c shows the resulting height of the lower fluid for the seventh order accurate WENO-AO-(7,5,3) scheme. The shock propagates into the unshocked medium without perturbation till it encounters the depression. When it encounters the "bubble" it sets up a left-going bow shock and a wake around the bubble (shown in green color in Fig. 18). A very interesting flow structure develops behind the "bubble" as can be seen in Fig. 18. The obtained results match with the reported results in Balsara *et al.* [10]. In Fig. 18c, we observe that the seventh order accurate schemes captures the shock-front (shown in red) more accurately than the fifth order accurate results (Fig. 18b). The ninth order scheme shows identical result to seventh order results, therefore, it is not shown here. In conclusion, the higher-order simulations show flow structures that are very crisp, highlighting the value of higher order schemes.

Next, we consider the Shock-Vortex interaction problem for the Two-Layer Shallow water equations. The computational domain for this problem spans $[-0.5, 1.5] \times [-0.5, 1.5]$ and the problem was run to a final time of 0.24 using the seventh order accurate HLL-based AFD WENO-AO-(7,3) scheme on a mesh with 600×600 zones. Figs. 19a, 19b and 19c show height of the upper fluid at times t=0.0, 0.06 and 0.24. The intermediate time of 0.06 corresponds to a time when the vortex has propagated halfway through the shock and a little bend in the shock front can be observed at this time. Figs. 19d, 19e and 19f show x-velocity of the upper fluid at times t=0.0, 0.06, 0.24. Similar to the Figs. 19a, 19b and 19c for the height of upper fluid, a stationary shock is clearly visible in the Figs. 19d, 19e, 19f for the x-velocity of the upper fluid. We observe that as the vortex propagates through the standing shock, the vortex is not fully destroyed by the shock flow. Instead, the vortex sheds some of its angular momentum after it has passed through the shock. In doing so, the vortex has rearranged its form. Analogous to the Baer-Nunziato flow equations, the two-layer shallow water equations preserve angular momentum. As a result, the core of the vortex is not destroyed as it passes through the mild shock. The obtained results are consistent with the results given in Balsara *et al.* [10]. The fifth and ninth order AFD-WENO schemes also show identical results, therefore, they are not shown here.

**VII.c) Two-dimensional Test Problems for the Multiphase Debris Flow**

In this Sub-section, we present two multi-dimensional problems for the Multiphase Debris Flow. Both problems are described in Balsara *et al.* [10]. As the detailed setup is already given in Balsara *et al.* [10], we do not repeat the details here.

The first problem is the shock-"bubble" interaction problem. The problem utilizes the Vortex used in Section V.c. Here we consider a domain that spans $[-0.5, 3.0] \times [-1.25, 1.25]$. The problem was run to a final time of 0.3 using the third, fifth and seventh order HLL-based AFD-WENO-AO schemes on a mesh with 700×500 zones. Fig. 20a shows final image of the height of the fluid phase for the third order accurate WENO-AO-3 scheme, Fig. 20b shows the resulting height of the fluid phase for the fifth order accurate WENO-AO-(5,3) scheme and Fig. 20c shows the resulting height of the fluid phase for the seventh order accurate WENO-AO-(7,5,3) scheme. When the shock encounters the "bubble," it creates a left-going bow shock. A ring-shaped structure originating from the area of the depression can be observed in Fig. 20. The ninth order scheme shows identical result to seventh order results, therefore, it is not shown here. Once again, we observe that the higher-order simulations (fifth and seventh order) exhibit flow structures that are highly sharp. This emphasizes the significance and advantages of employing higher-order schemes.

Next, we consider the Shock-Vortex interaction problem for the Multiphase Debris Flow model. The computational domain for this problem spans $[-0.5,1.5]\times[-0.5,1.5]$. The problem was run to a final time of 0.24 using the ninth order accurate HLL-based AFD WENO-AO-(9,3) scheme on a mesh with 600×600 zones. Figs. 21a, b, and c show the height of the upper fluid at times t=0.0, 0.06, and 0.24. The intermediate simulation time of 0.06 corresponds to a time when the vortex has propagated halfway through the shock. At this time, a slight reduction in the strength of the solid height can be observed in the vortex. Figs. 21d, 21e, 21f show the solid x-velocity at times t=0.0, 0.06, 0.24. The standing shock is visible in all of the panels of Fig. 21. The vortex sheds some of its angular momentum after it has passed through the mild shock. In doing so, the vortex has rearranged its form. The obtained results are consistent with the results given in Balsara *et al.* [10]. The fifth and seventh order AFD-WENO schemes also show identical results, therefore, they are not shown here.

**VIII) Test Problems with Stiff Source Terms**

AFD-WENO schemes demonstrate a significant advantage when dealing with stiff source terms. Unlike finite volume methods, these schemes operate pointwise, which includes the zone centers where the primal variables are located. Consequently, the enforcement of source terms aligns with the collocation of variables, leading to enhanced accuracy and low computation costs. We present two test problems that are based on the Baer-Nunziato compressible multi-phase flow with stiff source terms, as given in Dumbser and Boscheri [21]. The same set of problems were also used in Balsara *et al.* [10]. For both the problems, we use a third order accurate IMEX-SSP3(4,3,3) scheme from Pareschi and Russo [35] for the time integration.

The equations for the Baer-Nunziato compressible multi-phase flow with stiff source terms are given by

$$\partial_t(\phi_1\rho_1) + \nabla \cdot (\phi_1\rho_1\mathbf{v}_1) = 0$$
$$\partial_t(\phi_1\rho_1\mathbf{v}_1) + \nabla \cdot (\phi_1(\rho_1\mathbf{v}_1 \otimes \mathbf{v}_1 + \mathbf{I}\, p_1)) - P_I \nabla \phi_1 = -\lambda(\mathbf{v}_1 - \mathbf{v}_2)$$
$$\partial_t(\phi_1\rho_1 E_1) + \nabla \cdot (\phi_1\mathbf{v}_1(\rho_1 E_1 + p_1)) + P_I \partial_t \phi_1 = -\lambda \mathbf{V}_I \cdot (\mathbf{v}_1 - \mathbf{v}_2)$$
$$\partial_t(\phi_2\rho_2) + \nabla \cdot (\phi_2\rho_2\mathbf{v}_2) = 0$$
$$\partial_t(\phi_2\rho_2\mathbf{v}_2) + \nabla \cdot (\phi_2(\rho_2\mathbf{v}_2 \otimes \mathbf{v}_2 + \mathbf{I}\, p_2)) - P_I \nabla \phi_2 = -\lambda(\mathbf{v}_2 - \mathbf{v}_1)$$
$$\partial_t(\phi_2\rho_2 E_2) + \nabla \cdot (\phi_2\mathbf{v}_2(\rho_2 E_2 + p_2)) + P_I \partial_t \phi_2 = -\lambda \mathbf{V}_I \cdot (\mathbf{v}_2 - \mathbf{v}_1)$$
$$\partial_t \phi_1 + \mathbf{V}_I \cdot \nabla \phi_1 = \mu(p_1 - p_2)$$

The set of equations mentioned above contains stiff terms on the right-hand side, and the stiffness is controlled by two parameters, denoted as $\lambda$ and $\mu$.

The first test problem is a one-dimensional problem that spans the domain $[-0.5, 0.5]$. The left and right initial conditions were initialized around $x = 0$ and are specified in Table IX. The problem uses interphase drag $\lambda = 10^3$ and pressure relaxation parameter $\mu = 10^2$. These values turn this problem into a test problem with a moderately stiff source term. The simulation was run till time t=0.2 on a 400 zone mesh with a CFL of 0.8 using the fifth order accurate HLL-based AFD WENO-AO-(5,3) scheme. Figs. 22a, b, and c show the solid density, solid x-velocity, and solid pressure profiles, respectively. We see that the results are comparable to those presented in Dumbser and Boscheri [21]. The seventh and ninth order AFD-WENO schemes also show identical results, therefore, they are not shown here. The reference solution is also shown as the solid line in Fig. 22 and was obtained using a third order scheme with 4000 zones.

**Table IX gives the left states, right states, and final time of the one-dimensional Riemann problem for the Baer-Nunziato compressible multi-phase flow with stiff source terms.**

| | $\rho_1$ | $u_1$ | $p_1$ | $\rho_2$ | $u_2$ | $p_2$ | $\phi_1$ | $t_{end}$ |
|---|---|---|---|---|---|---|---|---|
| $\gamma_1 = 1.4,\ \pi_1 = 0,\ \gamma_2 = 1.67,\ \pi_2 = 0$ | | | | | | | | |
| **Left:** | 1.0 | 0.0 | 1.0 | 1.0 | 0.0 | 1.0 | 0.99 | 0.2 |

| Right: | 0.125 | 0.0 | 0.1 | 0.125 | 0.0 | 0.1 | 0.01 |

The second test problem is a two-dimensional Riemann problem. The computational domain spans $[-0.5, 0.5] \times [-0.5, 0.5]$. The initial states are described in Table X. The problem was run with a CFL of 0.4 to a stopping time of 0.15 on a mesh of $400 \times 400$ zones. Here we have used $\lambda = 10^5$ and $\mu = 10^2$ which corresponds to rather severely stiff source terms. We run the simulation using the seventh-order accurate HLL-based AFD WENO-AO-(7,3) scheme. Fig. 23a shows the solid density, Fig. 23b shows the gas density, and Fig. 23c shows the solid volume fraction. We observe that the obtained results are consistent with the solution presented in Dumbser and Boscheri [21]. The fifth and ninth order AFD-WENO schemes also show identical results, therefore, they are not shown here.

**Table X gives the initial states of the two-dimensional Riemann problem for the Baer-Nunziato compressible multi-phase flow with stiff source terms.**

| $x, y$ | $\rho_1$ | $u_1$ | $v_1$ | $p_1$ | $\rho_2$ | $u_2$ | $v_2$ | $p_2$ | $\phi_1$ |
|---|---|---|---|---|---|---|---|---|---|
| | $\gamma_1 = 1.4, \ \pi_1 = 0, \ \gamma_2 = 1.67, \ \pi_2 = 0$ | | | | | | | | |
| $x > 0, y > 0$ | 2.0 | 0.0 | 0.0 | 2.0 | 1.5 | 0.0 | 0.0 | 2.0 | 0.8 |
| $x < 0, y > 0$ | 1.0 | 0.0 | 0.0 | 1.0 | 0.5 | 0.0 | 0.0 | 1.0 | 0.4 |
| $x < 0, y < 0$ | 2.0 | 0.0 | 0.0 | 2.0 | 1.5 | 0.0 | 0.0 | 2.0 | 0.8 |
| $x > 0, y < 0$ | 1.0 | 0.0 | 0.0 | 1.0 | 0.5 | 0.0 | 0.0 | 1.0 | 0.4 |

## IX) Conclusions

Classical FD-WENO schemes have been available for conservation laws since the early papers by Shu and Osher [40], [41], Jiang and Shu [29] and Balsara and Shu [4]. Until very recently, all variants of the finite difference WENO scheme have indeed been restricted to treating only hyperbolic systems that are in conservation form. The recent emergence of several classes of hyperbolic systems with non-conservative products exposes a dire need for a new class of finite difference WENO schemes that can handle such systems. This has called for a reassessment of the entire design philosophy of finite difference WENO schemes. The first step in such a redesign took place very recently in Balsara *et al*. [10] where a very efficient set of finite difference WENO schemes were designed that could integrate hyperbolic systems with non-conservative products. The schemes in Balsara *et al*. [10] were demonstrated to work well for up to ninth order of accuracy. They can truly be thought of as being the analogues of the schemes in Shu and Osher [40], [41] when the hyperbolic system has non-conservative products. The schemes in Balsara *et al*. [10] were indeed significantly faster than the classical FD-WENO schemes. However, they were written in a fluctuation form which does not simply reduce to a flux conservative form when the system is indeed conservative. The present AFD-WENO schemes are very innovative, and to the best of our knowledge have never been presented in the literature. They are designed to take the field further so that finite difference WENO methods can handle hyperbolic systems with non-conservative products while at the same time respecting the Lax-Wendroff theorem when it is applicable.

The present schemes derive their start from alternative finite difference WENO (AFD-WENO) schemes for conservation laws. A well-developed AFD-WENO scheme offers many extremely desirable advantages over a classical FD-WENO scheme. First, it can work with any type of Riemann solver and different fields of study may have their own special Riemann solver that they find beneficial. Such a flexibility is not available in FD-WENO formulations, whereas it is indeed available in AFD-WENO formulations. This flexibility of invoking the Riemann solver at pointwise locations can also be exploited to ensure that the final scheme respects the preservation of free stream conditions on curvilinear meshes. On such curvilinear meshes, the flux reconstruction of classical FD-WENO methods becomes a liability; this enables AFD-WENO methods come to the fore because they are not restricted by this liability. The stumbling block in the development of AFD-WENO schemes was that general-purpose AFD-WENO schemes that can work with minimal or no changes for large classes of hyperbolic conservation laws were not in hand. This barrier was removed in Balsara *et al*. [11] where a general-purpose AFD-WENO scheme was designed for up to ninth order and shown to work well for large classes of hyperbolic conservation laws. The previously mentioned paper set the stage for the work reported in this paper.

If an AFD-WENO scheme could be general-purpose, available at several high orders and work well for large classes of conservation laws, then one is inclined to ask the following

question:- Why is a general-purpose, high order AFD-WENO scheme not available for large classes of hyperbolic systems that have non-conservative products? This question is the very important motivating question that animates this entire paper. In Section II we show that there is a way to obtain an answer to this very important question. We start with a general flux which we write as having a conservative part and a part which will only be thought of as a leading to a non-conservative product. This allows us to write the AFD-WENO update equation in a fluctuation form. Once this door is opened, we carefully track which steps retain flux conservation and which steps in the derivation force us to relinquish flux conservation. The end result is shown in eqn. (15) which presents a general-purpose, high order AFD-WENO scheme that is well-suited for large classes of hyperbolic systems that have non-conservative products. The scheme we derive in eqn. (15) is one which retains a flux conservative update for those components of the PDE that are indeed in conservation form and judiciously relinquishes it for those components of the PDE that are dominated by non-conservative products. In retaining conservation form as much as possible, the dictates of the Lax-Wendroff theorem are also preserved as much as possible. Two variants are designed, the first one being documented in eqn. (15) and the second one in eqn. (13). Eqn. (15) is based on making three WENO interpolation steps; where the first WENO interpolation is done in characteristic space and the next two WENO interpolations are lower cost and can be done component by component. Eqn. (13) is based on making two WENO interpolation steps; where the first WENO interpolation is done in characteristic space and the next WENO interpolation is lower cost and can be done component by component.

With the formulation in hand, Section III addresses some issues associated with the non-linear hybridization and stabilization of the higher order derivatives. Section IV provides point by point details on how the method is implemented.

Section V shows accuracy tests. It shows that for three very broad classes of hyperbolic systems with non-conservative products – i.e. the Baer-Nunziato system, the two layer shallow water system and the multiphase debris flow model of Pitman and Le – the method achieves its design accuracies. Equally importantly, in Sub-section V.d we show a speed comparison of the classical FD-WENO scheme, the WENO scheme from Balsara *et al*. [10], the AFD-WENO scheme from Balsara *et al*. [11] which only applies to conservation laws, and the two AFD-WENO schemes from this paper which can be applied to any hyperbolic system with fluxes as well as non-conservative products. The speeds of the schemes conform to our expectations and all schemes are shown to be roughly competitive in performance. The differences in speeds are explained thoroughly in Sub-section V.d. In that Sub-section we also provide a discussion on the implementation of finite difference WENO schemes on modern-day CPUs with large caches. We show a very desirable result that higher order variants of finite difference WENO schemes don't cost that much more than their lower order variants. The scientific explanation for that stems from our realization that the larger number of floating point operations associated with larger stencils are almost free of charge if the finite difference WENO code is designed to be cache friendly. This should have great, and very beneficial, implications for the role of finite difference WENO schemes in Peta- and Exascale computing! Simply put:- Our well-implemented AFD-WENO schemes come pretty close to the promised land of high performance computing where float point operations are almost free and the dominant costs are the other fixed costs in the implementation!

Section VI shows several stringent one-dimensional test problems drawn from our three PDE systems of interest. Section VII shows multidimensional test problems, again drawn from our three PDE systems of interest. Section VIII shows test problems with stiff source terms, demonstrating the very favorable handling and features of AFD-WENO in the stiff limit.

**Acknowledgements**

DSB acknowledges support via NSF grant NSF-AST-2009776, NASA grant NASA-2020-1241 and NASA grant 80NSSC22K0628. DSB and HK acknowledge support from a Vajra award, VJR/2018/00129 and also a travel grant from Notre Dame International. CWS acknowledges support via NSF grant DMS-2309249.

**Ethical Statement**
**i. Compliance with Ethical Standards** : This manuscript complies with all ethical standards for scientific publishing.
**ii. (in case of Funding) Funding** : The funding has been acknowledged. DSB acknowledges support via NSF grants NSF-19-04774, NSF-AST-2009776, NASA-2020-1241 and NASA-80NSSC22K0628. DSB and HK acknowledge support from a Vajra award, VJR/2018/00129. CWS acknowledges support via AFOSR grant FA9550-20-1-0055 and NSF grant DMS-2010107.
**iii. Conflict of Interest** : On behalf of all authors, the corresponding author states that there is no conflict of interest.

**iv. Ethical approval** : N/A
**v. Informed consent** : N/A

**vi. Data Statement** : All data that was used in the generation of the figures has been stored and available for later use.

**Appendix A**

The following variant of eqn. (15) may also be used:-

$$\partial_t \mathbf{U}_i = -\frac{1}{\Delta x}\left\{\mathbf{D}^{*-}_{i+1/2}\left(\hat{\mathbf{U}}^-_{i+1/2}, \hat{\mathbf{U}}^+_{i+1/2}\right) + \mathbf{D}^{*+}_{i-1/2}\left(\hat{\mathbf{U}}^-_{i-1/2}, \hat{\mathbf{U}}^+_{i-1/2}\right)\right\} - \frac{1}{\Delta x}\left\{\mathbf{F}\left(\hat{\mathbf{U}}^-_{i+1/2}\right) - \mathbf{F}\left(\hat{\mathbf{U}}^+_{i-1/2}\right)\right\}$$

$$-\frac{1}{\Delta x}\mathbf{C}(\mathbf{U}_i)\left(\hat{\mathbf{U}}^-_{i+1/2} - \hat{\mathbf{U}}^+_{i-1/2}\right)$$

$$-\frac{1}{\Delta x}\mathbf{C}(\mathbf{U}_i)\left\{\begin{array}{l}\left[\begin{array}{l}-\frac{1}{24}(\Delta x)^2\left[\partial_x^2\hat{\mathbf{U}}\right]_{i+1/2} + \frac{7}{5760}(\Delta x)^4\left[\partial_x^4\hat{\mathbf{U}}\right]_{i+1/2} \\ -\frac{31}{967680}(\Delta x)^6\left[\partial_x^6\hat{\mathbf{U}}\right]_{i+1/2} + \frac{127}{154828800}(\Delta x)^8\left[\partial_x^8\hat{\mathbf{U}}\right]_{i+1/2}\end{array}\right] \\ -\left[\begin{array}{l}-\frac{1}{24}(\Delta x)^2\left[\partial_x^2\hat{\mathbf{U}}\right]_{i-1/2} + \frac{7}{5760}(\Delta x)^4\left[\partial_x^4\hat{\mathbf{U}}\right]_{i-1/2} \\ -\frac{31}{967680}(\Delta x)^6\left[\partial_x^6\hat{\mathbf{U}}\right]_{i-1/2} + \frac{127}{154828800}(\Delta x)^8\left[\partial_x^8\hat{\mathbf{U}}\right]_{i-1/2}\end{array}\right]\end{array}\right\}$$

$$-\frac{1}{\Delta x}\left\{\begin{array}{l}\left[\begin{array}{l}-\frac{1}{24}(\Delta x)^2\left[\partial_x\left(\mathbf{B}\partial_x\hat{\mathbf{U}}\right)\right]_{i+1/2} + \frac{7}{5760}(\Delta x)^4\left[\partial_x^3\left(\mathbf{B}\partial_x\hat{\mathbf{U}}\right)\right]_{i+1/2} \\ -\frac{31}{967680}(\Delta x)^6\left[\partial_x^5\left(\mathbf{B}\partial_x\hat{\mathbf{U}}\right)\right]_{i+1/2} + \frac{127}{154828800}(\Delta x)^8\left[\partial_x^7\left(\mathbf{B}\partial_x\hat{\mathbf{U}}\right)\right]_{i+1/2}\end{array}\right] \\ -\left[\begin{array}{l}-\frac{1}{24}(\Delta x)^2\left[\partial_x\left(\mathbf{B}\partial_x\hat{\mathbf{U}}\right)\right]_{i-1/2} + \frac{7}{5760}(\Delta x)^4\left[\partial_x^3\left(\mathbf{B}\partial_x\hat{\mathbf{U}}\right)\right]_{i-1/2} \\ -\frac{31}{967680}(\Delta x)^6\left[\partial_x^5\left(\mathbf{B}\partial_x\hat{\mathbf{U}}\right)\right]_{i-1/2} + \frac{127}{154828800}(\Delta x)^8\left[\partial_x^7\left(\mathbf{B}\partial_x\hat{\mathbf{U}}\right)\right]_{i-1/2}\end{array}\right]\end{array}\right\} \quad (A.1)$$

It slightly reduces the cost of the interpolation that is needed in the last curly bracket. But it requires the evaluation of the matrix multiplication $\left(\mathbf{B}\partial_x\hat{\mathbf{U}}\right)_i$ at each zone center "$i$" in order to interpolate its derivatives to the zone boundaries.

**Appendix B) Phase-volume based Flattener function for the Baer-Nunziato model**

Flattening algorithms are used in simulations to detect location of the strong shocks within the computational domain. Their purpose is to enhance the accuracy of the simulation by reducing numerical artifacts near discontinuities. For Euler flow, a flattening function that is based on the divergence of the velocity field and characteristic speed has been presented in Colella and Woodward [17], and Balsara [6] and it can also be used here. However, we have found that a sudden change in the volume fraction can also be a source of spurious oscillations for the Baer-Nunziato system. Therefore, we propose a one-dimensional phase-volume based flattener function for the Baer-Nunziato model. The flattener algorithm trivially extends to the two- and three-dimensional case. The method begins by computing the sum and difference of extremal values (denoted by $\Delta_i^+\phi$ and $\Delta_i^-\phi$, respectively) of the solid volume fractions in cell $i$ over the neighbouring cells $i$-$1$, $i$ and $i$+$1$. Mathematically, we define the quantities $\Delta_i^+\phi$ and $\Delta_i^-\phi$ as follows.

$$\Delta_i^+\phi = \max[\phi_1^{i-1}, \phi_1^i, \phi_1^{i+1}] + \min[\phi_1^{i-1}, \phi_1^i, \phi_1^{i+1}],$$
$$\Delta_i^-\phi = \max[\phi_1^{i-1}, \phi_1^i, \phi_1^{i+1}] - \min[\phi_1^{i-1}, \phi_1^i, \phi_1^{i+1}],$$

where $\phi_1^i$ is the solid volume fraction within the cell $i$. In each zone $i$, the flattener function (denoted by $\eta_i$) is defined as

$$\eta_i = \min\left[1, \max\left[0, 0.5(a_i - \kappa)\right]\right],$$

where $\kappa$ is a positive constant number which works best for us when we choose $\kappa = 0.1$. The expression for the quantity $a_i$ is given by

$$a_i = \max\left[\frac{|\Delta_i^-\phi|}{0.5\,\Delta_i^+\phi + \epsilon}, \frac{|\Delta_i^-\phi|}{|1 - 0.5\,\Delta_i^+\phi| + \epsilon}\right], \quad \text{where } \epsilon = 10^{-12} \text{ is a small number.}$$

The flattener function does not modify the reconstruction/interpolation when the flow is smooth or consists of rarefactions, and in that case, $\eta_i = 0$. However, the flattener function gradually increases from $\eta_i = 0$ to $\eta_i = 1$ when strong shocks are present in the volume fraction.

**References**

[1] N. Andrianov, G. Warnecke, *The Riemann problem for the Baer–Nunziato two-phase flow model*, Journal Computational Physics, 212 (2004) 434–464.

[2] M.R. Baer, J.W. Nunziato, *A two-phase mixture theory for the deflagration-to-detonation transition (DDT) in reactive granular materials*, Int. J. Multiph. Flow 12 (1986) 861–889.

[3] D.S. Balsara, *Total Variation Diminishing Algorithm for Adiabatic and Isothermal Magnetohydrodynamics*, Ap.J. Supp., Vol. 116, Pgs. 133-153 (1998)

[4] D. S. Balsara and C.-W. Shu, *Monotonicity Preserving Weighted Non-oscillatory schemes with increasingly High Order of Accuracy*, Journal of Computational Physics, 160 (2000) 405-452

[5] D. S. Balsara, T. Rumpf, M. Dumbser & C.-D. Munz, *Efficient, high-accuracy ADER-WENO schemes for hydrodynamics and divergence-free magnetohydrodynamics*, Journal Computational Physics, 228 (2009) 2480

[6] D. S. Balsara, *Self-Adjusting, Positivity Preserving High Order Schemes for Hydrodynamics and Magnetohydrodynamics*, Journal Computational Physics, Vol. 231 (2012) Pgs. 7504-7517

[7] D.S. Balsara, S. Garain and C.-W. Shu, *An efficient class of WENO schemes with adaptive order*, Journal of Computational Physics, 326 (2016) 780-804

[8] D.S. Balsara, S. Garain, V. Florinski, W. Boscheri, *An Efficient Class of WENO Schemes with Adaptive Order for Unstructured Meshes*, Journal of Computational Physics, 404 (2020) 109062

[9] D.S. Balsara, S. Samantaray and S. Subramanian, *Efficient WENO-Based Prolongation Strategies for Divergence-Preserving Vector Fields*, Communications on Applied Mathematics and Computation, (2022)

[10] D.S. Balsara, D. Bhoriya, C.-W. Shu, H., Kumar, *Efficient Finite Difference WENO Scheme for Hyperbolic Systems with Non-Conservative Products*, Communications on Applied Mathematics and Computation, (2023), to appear

[11] D.S. Balsara, D. Bhoriya, C.-W. Shu, H., Kumar, *Efficient Alternative Finite Difference WENO Schemes for Hyperbolic Conservation Laws*, Communications on Applied Mathematics and Computation, (2023), submitted

[12] Borges, R., Carmona, M., Costa, B., Don, W.S., *An improved weighted essentially non-oscillatory scheme for hyperbolic conservation laws*, Journal of Computational Physics, 227 (6) (2008) 3101–3211

[13] M. Castro, B. Costa, W.S. Don, *High order weighted essentially non-oscillatory WENO-Z schemes for hyperbolic conservation laws*, Journal of Computational Physics, 230 (2011) 1766–1792

[14] M.J. Castro, A. Pardo, C. Parés, E.F. Toro, *On some fast well-balanced first order solvers for nonconservative systems*, Math. Comput. 79 (2010) 1427–1472

[15] M. Castro, J. E. Gallardo and C. Pares, *High order finite volume schemes based on reconstruction of states for solving hyperbolic systems with nonconservative products, applications to shallow-water systems*, Mathematics of Computation, vol. 75, number 255 (2006) 1103-1134

[16] S. Chiochetti and C. Müller, A Solver for Stiff Finite-Rate Relaxation in Baer–Nunziato Two-Phase Flow Models, in preparation, (2022)

[17] P. Colella and P. R. Woodward, *The piecewise parabolic method (PPM) for gas-dynamical simulations*. Journal of computational physics, 54(1) (1984), 174-201

[18] F. Coquel, C. Marmignon, P. Rai and F. Renac, An entropy stable high-order discontinuous Galerkin spectral element method for the Baer-Nunziato two-phase flow model, Journal of Computational Physics, 431 (2021) 110135

[19] M. Dumbser, M. Castro, C. Parés, E.F. Toro, *ADER schemes on unstructured meshes for non-conservative hyperbolic systems: applications to geophysical flows*, Comput. Fluids 38 (2009) 1731–1748

[20] M. Dumbser, A. Hidalgo, M. Castro, C. Parés, E.F. Toro, *FORCE schemes on unstructured meshes II: non-conservative hyperbolic systems*, Comput. Methods Appl. Mech. Eng. 199 (2010) 625–647

[21] M. Dumbser and W. Boscheri, *High-order unstructured Lagrangian one-step WENO finite volume schemes for non-conservative hyperbolic systems: Applications to compressible multi-phase flows*, Computers and Fluids, 86 (2013) 405-432

[22] M. Dumbser, A. Hidalgo and O. Zanotti, *High-order space-time adaptive ADER-WENO finite volume schemes for non-conservative hyperbolic systems*, Comput. Methods Appl. Mech. Engg., 268 (2014) 359-387

[23] M. Dumbser, U. Iben, M. Ioriatti, *An efficient semi-implicit finite volume method for axially symmetric compressible flows in compliant tubes*, Appl. Numer. Math. 89 (2015) 24–44.

[24] M. Dumbser and D.S. Balsara, *A New, Efficient Formulation of the HLLEM Riemann Solver for General Conservative and Non-Conservative Hyperbolic Systems*, Journal of Computational Physics 304 (2016) 275-319

[25] Z. Gao, L.-L. Fang, B.-S. Wang, Y. Wang, W. S. Don, *Seventh and ninth orders characteristic-wise alternative WENO finite difference schemes for hyperbolic conservation laws*, Computers and Fluids 202 (2020) 104519

[26] G.A. Gerolymos, D. Sénéchal and I. Vallet, *Very high order WENO schemes, Journal of Computational Physics*, 228 (2009) 8481-8524

[27] A. Harten, B. Engquist, S.Osher and S. Chakravarthy, *Uniformly high order essentially non-oscillatory schemes III*, Journal of Computational Physics, 71 (1987) 231-303

[28] A.K. Henrick, T.D. Aslam and J.M. Powers, *Mapped weighted essentially non-oscillatoriy schemes: Achieving optimal order near critical points*, Journal of Computational Physics 207 (2006) 542-567

[29] G.-S. Jiang and C.-W. Shu, *Efficient implementation of weighted ENO schemes*, Journal of Computational Physics, 126 (1996) 202-228

[30] Y. Jiang, C-W Shu, and M. Zhang, *An alternative formulation of finite difference ENO schemes with Lax-Wendroff time discretization for conservation laws*, SIAM J. Sci. Comput., 35(2) (2013) A1137–A1160

[31] Y. Jiang, C-W Shu, and M. Zhang, *Free-stream preserving finite-difference schemes on curvilinear meshes*, Methods and Applications of Analysis, 21(1), (2014), 001–030

[32] F. Kupka, N. Happenhofer, I. Higueras, O. Koch, *Total-variation-diminishing implicit–explicit Runge–Kutta methods for the simulation of double-diffusive convection in astrophysics*, Journal of Computational Physics, 231 (2012) 3561-3586

[33] X.-D. Liu, S. Osher and T. Chan, *Weighted essentially non-oscillatory schemes*, Journal of Computational Physics 115 (1994) 200-212

[34] B. Merriman, *Understanding the Shu–Osher Conservative Finite Difference Form*, Journal of Scientific Computing, 19(1–3), (2003), 309

[35] L. Pareschi and G. Russo, *Implicit-explicit Runge-Kutta schemes and applications to hyperbolic systems with relaxation*, Journal of Scientific Computing, 25 (2005) 129

[36] M. Pelanti, F. Bouchut, A. Mangeney, *A Roe-type scheme for two-phase shallow granular flows over variable topography*, Math. Model. Numer. Anal. 42 (2008) 851–885

[37] E.B. Pitman, L. Le, *A two-fluid model for avalanche and debris flows*, Philos. Trans. R. Soc. A 363 (2005) 1573–1601

[38] S. Rhebergen, O. Bokhove, J.J.W. van der Vegt, *Discontinuous Galerkin finite element methods for hyperbolic nonconservative partial differential equations*, Journal of Computational Physics, 227 (2008) 1887–1922

[39] R. Saurel, R. Abgrall, *A multiphase Godunov method for compressible multifluid and multiphase flows*, Journal of Computational Physics, 150 (1999) 425–467.

[40] C.-W. Shu and S. J. Osher, *Efficient implementation of essentially non-oscillatory shock capturing schemes*, Journal of Computational Physics, 77 (1988) 439-471

[41] C.-W. Shu and S. J. Osher, *Efficient implementation of essentially non-oscillatory shock capturing schemes II*, Journal of Computational Physics, 83 (1989) 32-78

[42] Shu, C.-W., *High order weighted essentially non-oscillatory schemes for convection dominated problems*, SIAM Review, 51 (2009) 82-126

[43] Shu, C.-W., *Essentially non-oscillatory and weighted essentially non-oscillatory schemes*, Acta Numerica, v29 (2020), pp.701-762

[44] Spiteri, R.J. and Ruuth, S.J., *A new class of optimal high-order strong-stability-preserving time-stepping schemes*, SIAM Journal of Numerical Analysis, 40 (2002), pp. 469–491

[45] Spiteri, R.J. and Ruuth, S.J., *Non-linear evolution using optimal fourth-order strong-stability-preserving Runge-Kutta methods*, Mathematics and Computers in Simulation 62 (2003) 125-135

[46] S.A. Tokareva, E.F. Toro, *HLLC-type Riemann solver for the Baer–Nunziato equations of compressible two-phase flow*, Journal of Computational Physics, 229 (2010) 3573–3604.

[47] P. Woodward and P. Colella, *The numerical simulation of two-dimensional fluid flow with strong shocks*, Journal of Computational Physics 54 (1984), 115-173

[48] F. Zheng, C.-W. Shu, J. Qiu, *A high order conservative finite difference scheme for compressible two-medium flows*, Journal of Computational Physics 445 (2021) 110597

[49] J. Zhu and J. Qiu, *A new fifth order finite difference WENO scheme for solving hyperbolic conservation laws*, accepted, Journal of Computational Physics (2016)

[50] J. Zhu and C.-W. Shu, A new type of multi-resolution WENO schemes with increasingly higher order of accuracy, Journal of Computational Physics, 375 (2018), 659-683

**Figure Captions**

*Fig. 1 shows part of the mesh around zone "i". The mesh functions are collocated at the zone centers, as shown by the thick dots. The zone boundaries are shown by the vertical lines. The figure also shows the stencils associated with the zone "i" for the third and fifth order **pointwise** WENO-AO interpolation strategies. We have three smaller third order stencils and a large fifth order stencil. For third order WENO-AO, only the three smaller stencils are used, whereas the larger stencil is also used for fifth order WENO-AO. The interpolated variables at the zone boundaries are shown with a caret. The variables with a superscript star are resolved states obtained by the pointwise application of a Riemann solver at the zone boundaries. The process described here can also be adapted for Multiresolution WENO interpolation.*

*Fig. 2 shows part of the mesh around zone boundary "i+1/2". The fluxes are evaluated pointwise at the zone centers, as shown by the thick dots. The zone boundaries are shown by the vertical lines. The figure also shows the stencils associated with the zone boundary "i+1/2" for the third and fifth order AFD-WENO schemes. We have two smaller third order stencils and a large sixth order stencil. For a third order AFD-WENO scheme, the two smaller stencils can be non-linearly hybridized. In that case, the second derivatives of the flux can be obtained at the zone boundary when the smoothness in the solution warrants it. For fifth order AFD-WENO, the two smaller stencils can be non-linearly hybridized along with the larger stencil. In that case, the second and fourth derivatives of the flux can be obtained at the zone boundary when the smoothness in the solution warrants it. The process described here can be done for Adaptive Order and Multiresolution WENO interpolation.*

*Fig. 3 shows part of the mesh around zone boundary "i+1/2". The state vectors are available pointwise at the zone centers, as shown by the thick dots. The zone boundaries are shown by the vertical lines. The figure also shows the stencils associated with the zone boundary "i+1/2" for the third and fifth order AFD-WENO schemes. We have two smaller third order stencils and a large sixth order stencil. For a third order AFD-WENO scheme, the two smaller stencils can be non-linearly hybridized. In that case, the second derivatives of the state can be obtained at the zone boundary when the smoothness in the solution warrants it. For fifth order AFD-WENO, the two smaller stencils can be non-linearly hybridized along with the larger stencil. In that case, the second and fourth derivatives of the state can be obtained at the zone boundary when the*

*smoothness in the solution warrants it. The process described here can be done for Adaptive Order and Multiresolution WENO interpolation.*

*Fig. 4) Baer-Nunziato: Abgrall Problem using the 5$^{th}$ order accurate LLF-based AFD-WENO scheme with 200 zones. Fig. a) shows the solid volume fraction and Fig. b) shows the velocities and pressures for both the phases. The 7$^{th}$ and 9$^{th}$ order AFD-WENO schemes also show identical results and therefore they are not shown here.*

*Fig. 5) Baer-Nunziato: Riemann Problem-1 using the 5$^{th}$ order accurate HLL-based AFD-WENO scheme with 200 zones. Fig. a) shows the density for the solid phase and Fig. b) shows the density for the gas phase. The 7$^{th}$ and 9$^{th}$ order AFD-WENO schemes also show identical results and therefore they are not shown here.*

*Fig. 6) Baer-Nunziato: Riemann Problem-2 using the 7$^{th}$ order accurate HLL-based AFD-WENO scheme with 200 zones. Fig. a) shows the density for the solid phase and Fig. b) shows the density for the gas phase. The 5$^{th}$ and 9$^{th}$ order AFD-WENO schemes also show identical results and therefore they are not shown here.*

*Fig. 7) Baer-Nunziato: Riemann Problem-3 using the 9$^{th}$ order accurate HLL-based AFD-WENO scheme with 200 zones. Fig. a) shows the density for the solid phase and Fig. b) shows the density for the gas phase. The 5$^{th}$ and 7$^{th}$ order AFD-WENO schemes also show identical results and therefore they are not shown here.*

*Fig. 8) Baer-Nunziato: Riemann Problem-4 using the 5$^{th}$ order accurate LLF-based AFD-WENO scheme with 200 zones. Fig. a) shows the density for the solid phase and Fig. b) shows the density for the gas phase. The 7$^{th}$ and 9$^{th}$ order AFD-WENO schemes also show identical results and therefore they are not shown here.*

*Fig. 9) Baer-Nunziato: Riemann Problem-5 using the 7$^{th}$ order accurate LLF-based AFD-WENO scheme with 200 zones. Fig. a) shows the density for the solid phase and Fig. b) shows the density for the gas phase. We have used flattener for this problem. The 5$^{th}$ and 9$^{th}$ order AFD-WENO schemes also show identical results and therefore they are not shown here.*

*Fig. 10) Baer-Nunziato: Riemann Problem-6 using the 9$^{th}$ order accurate LLF-based AFD-WENO scheme with 200 zones. Fig. a) shows the density for the solid phase and Fig. b) shows the density*

*for the gas phase. The 5$^{th}$ and 7$^{th}$ order AFD-WENO schemes also show identical results and therefore they are not shown here.*

*Fig. 11) Two-Layer Shallow water: Riemann problem-1 using the 5$^{th}$ order accurate HLLI-based AFD-WENO scheme with 200 zones. The 7$^{th}$ and 9$^{th}$ order AFD-WENO schemes also show identical results and therefore they are not shown here. The top row shows the results from the scheme developed here. The bottom row of three figures shows the results from the same scheme when central differences are used for the higher order derivatives. We see that using just the central differences results in unphysical oscillations. Similar oscillations would be found if central differences were used for a 7$^{th}$ or 9$^{th}$ order scheme.*

*Fig. 12) Two-Layer Shallow water: Fig. a) shows the results of Riemann problem-2 using the 7$^{th}$ order accurate HLLI-based AFD-WENO scheme and Fig. b) shows the results of Riemann problem-3 using the 9$^{th}$ order accurate HLLI-based AFD-WENO scheme with 200 zones. The 5$^{th}$ and 9$^{th}$ order scheme for the Riemann Problem-2, and the 5$^{th}$ and 7$^{th}$ order scheme for the Riemann Problem-3 also show identical results and therefore they are not shown here. Fig. c) shows the results of Riemann problem-3 using the 9$^{th}$ order scheme when only central differences are used for the higher order derivatives. We clearly see the formation of unphysical oscillations when only central differences are used. Similar oscillations would be found if central differences were used for a 5$^{th}$ or 7$^{th}$ order scheme.*

*Fig. 13a-e) Debris Flow: a to d show the jump in the linearly degenerate fields using the 5$^{th}$ order accurate HLLI-based AFD-WENO scheme with 200 zones. e shows the velocity. The 7$^{th}$ and 9$^{th}$ order AFD-WENO schemes also show identical results and therefore they are not shown here.*

*Fig. 13f-j) Debris Flow: f to i show the jump in the linearly degenerate fields using the 5$^{th}$ order accurate HLLI-based AFD-WENO scheme with 200 zones when only central differences are used for the higher order derivatives. j shows the velocity. The 7$^{th}$ and 9$^{th}$ order AFD-WENO schemes also show identical results and therefore they are not shown here.*

*Fig. 14) Debris Flow: Riemann problem-2 using the 7$^{th}$ order accurate HLLI-based AFD-WENO scheme with 200 zones. The 5$^{th}$ and 9$^{th}$ order AFD-WENO schemes also show identical results and therefore they are not shown here.*

*Fig. 15) Debris Flow: Riemann problem-3 using the 9$^{th}$ order accurate HLLI-based AFD-WENO scheme with 200 zones. The 5$^{th}$ and 7$^{th}$ order AFD-WENO schemes also show identical results and therefore they are not shown here.*

*Fig. 16a) Baer-Nunziato model: Shock-bubble interaction problem using the 3$^{rd}$ order accurate HLL-based AFD WENO-AO-3 scheme with 700×300 zones. The solid density profiles have been shown.*

*Fig. 16b) Baer-Nunziato model: Shock-bubble interaction problem using the 5$^{th}$ order accurate HLL-based AFD WENO-AO-(5,3) scheme with 700×300 zones. The solid density profiles have been shown.*

*Fig. 16c) Baer-Nunziato model: Shock-bubble interaction problem using the 7$^{th}$ order accurate HLL-based AFD WENO-AO-(7,5,3) scheme with 700×300 zones. The solid density profiles have been shown. The 9$^{th}$ order scheme also shows identical result and therefore it is not shown here.*

*Fig. 17) Baer-Nunziato: Shock-Vortex Interaction using the 5$^{th}$ order accurate HLL-based AFD WENO-AO-(5,3) scheme with 600×600 zones at time levels t=0.0, 0.23 and 0.84. Figs. 17a, 17b and 17c show solid volume fraction at times t=0.0, 0.23, 0.84. Figs. 17d, 17e, 17f show solid x-velocity at times t=0.0, 0.23, 0.84. For the solid volume fraction, 30 contours were fit between a range of 0.33 and 0.530. For the solid x-velocity, 30 contours were fit between a range of -0.5 and 1.95.*

*Fig. 18a) Two-Layer Shallow water: Shock-bubble interaction problem using the 3$^{rd}$ order accurate HLL-based AFD WENO-AO-3 scheme with 700×300 zones. The solid density profiles has been shown.*

*Fig. 18b) Two-Layer Shallow water: Shock-bubble interaction problem using the 5$^{th}$ order accurate HLL-based AFD WENO-AO-(5,3) scheme with 700×300 zones. The solid density profiles has been shown.*

*Fig. 18c) Two-Layer Shallow water: Shock-bubble interaction problem using the 7$^{th}$ order accurate HLL-based AFD WENO-AO-(7,5,3) scheme with 700×300 zones. The solid density profiles has been shown. The 9$^{th}$ order scheme also shows identical result and therefore it is not shown here.*

*Fig. 19) Two-Layer Shallow water: Shock-Vortex Interaction using the 7$^{th}$ order accurate HLL-based AFD WENO-AO-(7,3) scheme with 600×600 zones at time levels t=0.0, 0.06 and 0.24. Figs. 19a, 19b and 19c show height of the upper fluid at times t=0.0, 0.06, 0.24. Figs. 19d, 19e, 19f show x-velocity of the upper fluid at times t=0.0, 0.06, 0.24. For the height, 40 contours were fit between a range of 1.0 and 3.4. For the velocity, 40 contours were fit between a range of 1.8 and 5.6.*

*Fig. 20a) Debris Flow: Shock-bubble interaction problem using the 3$^{rd}$ order accurate HLL-based AFD WENO-AO-3 scheme with 700×500 zones. The solid density profiles have been shown.*

*Fig. 20b) Debris Flow: Shock-bubble interaction problem using the 5$^{th}$ order accurate HLL-based AFD WENO-AO-(5,3) scheme with 700×500 zones. The solid density profiles have been shown.*

*Fig. 20c) Debris Flow: Shock-bubble interaction problem using the 7$^{th}$ order accurate HLL-based AFD WENO-AO-(7,5,3) scheme with 700×500 zones. The solid density profiles have been shown. The 9$^{th}$ order accurate scheme also shows identical results and therefore it is not shown here.*

*Fig. 21) Debris Flow: Shock-Vortex Interaction using the 9$^{th}$ order accurate HLL-based AFD WENO-AO-(9,3) scheme with 600×600 zones at time levels t=0.0, 0.06 and 0.24. Figs. 21a, 21b and 21c show the solid height at times t=0.0, 0.06, 0.24. Figs. 21d, 21e, 21f show the solid x-velocity at times t=0.0, 0.06, 0.24. For the height, 30 contours were fit between a range of 1.0 and 3.0. For the velocity, 30 contours were fit between a range of 1.5 and 5.6.*

*Fig. 22) Baer-Nunziato with stiff source: Results for the one-dimensional Riemann Problem using the 5$^{th}$ order accurate HLL-based AFD WENO-AO-(5,3) scheme with 400 zones. Figs. 22a, 22b and 22c show the solid density, solid x-velocity and solid pressure.*

*Fig. 23) Baer-Nunziato with stiff source: Results for the two-dimensional Riemann Problem using the 7$^{th}$ order accurate HLL-based AFD WENO-AO-(7,3) scheme with 400×400 zones. Fig. 23a shows the solid density, Fig. 23b shows the gas density and Fig. 23c shows the solid volume fraction. 30 equidistant contour lines are shown over the color plots.*

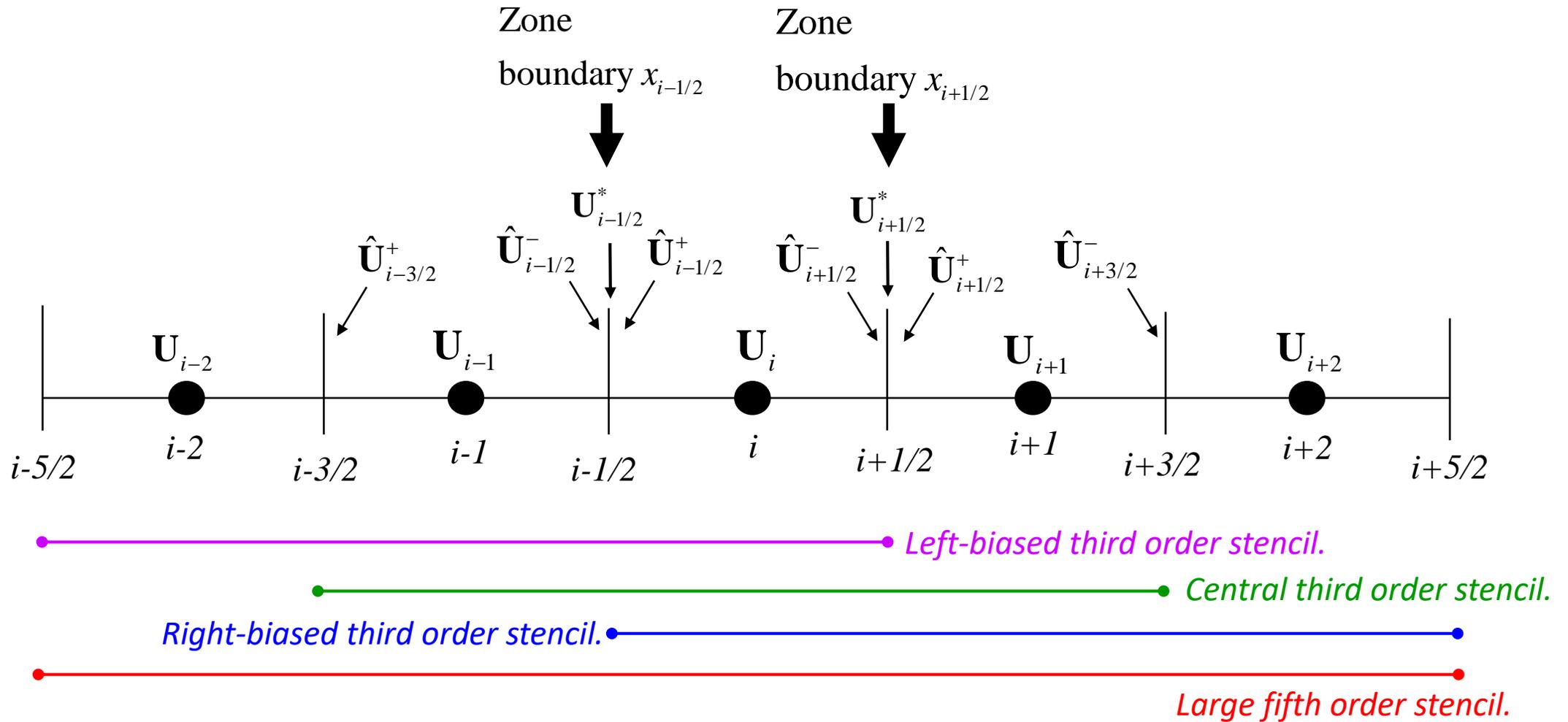

Fig. 1 shows part of the mesh around zone "i". The mesh functions are collocated at the zone centers, as shown by the thick dots. The zone boundaries are shown by the vertical lines. The figure also shows the stencils associated with the zone "i" for the third and fifth order **pointwise** WENO-AO interpolation strategies. We have three smaller third order stencils and a large fifth order stencil. For third order WENO-AO, only the three smaller stencils are used, whereas the larger stencil is also used for fifth order WENO-AO. The interpolated variables at the zone boundaries are shown with a caret. The variables with a superscript star are resolved states obtained by the pointwise application of a Riemann solver at the zone boundaries. The process described here can also be adapted for Multiresolution WENO interpolation.

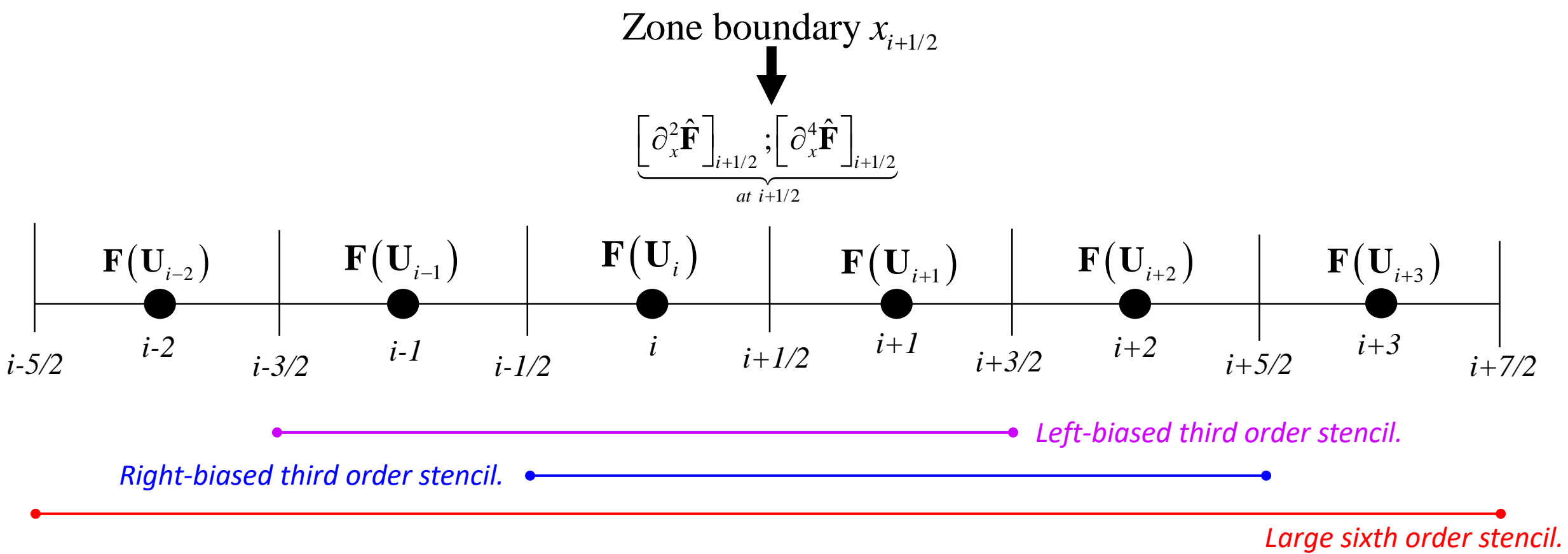

Fig. 2 shows part of the mesh around zone boundary "i+1/2". The fluxes are evaluated pointwise at the zone centers, as shown by the thick dots. The zone boundaries are shown by the vertical lines. The figure also shows the stencils associated with the zone boundary "i+1/2" for the third and fifth order AFD-WENO schemes. We have two smaller third order stencils and a large sixth order stencil. For a third order AFD-WENO scheme, the two smaller stencils can be non-linearly hybridized. In that case, the second derivatives of the flux can be obtained at the zone boundary when the smoothness in the solution warrants it. For fifth order AFD-WENO, the two smaller stencils can be non-linearly hybridized along with the larger stencil. In that case, the second and fourth derivatives of the flux can be obtained at the zone boundary when the smoothness in the solution warrants it. The process described here can be done for Adaptive Order and Multiresolution WENO interpolation.

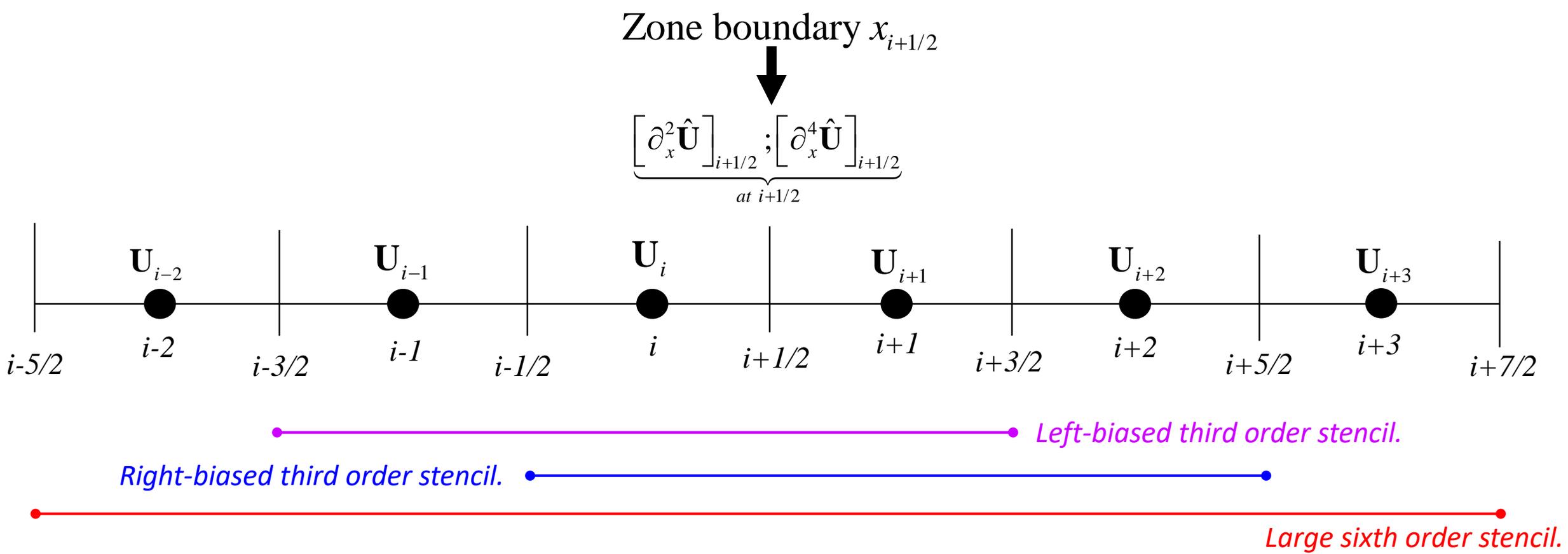

Fig. 3 shows part of the mesh around zone boundary "i+1/2". The state vectors are available pointwise at the zone centers, as shown by the thick dots. The zone boundaries are shown by the vertical lines. The figure also shows the stencils associated with the zone boundary "i+1/2" for the third and fifth order AFD-WENO schemes. We have two smaller third order stencils and a large sixth order stencil. For a third order AFD-WENO scheme, the two smaller stencils can be non-linearly hybridized. In that case, the second derivatives of the state can be obtained at the zone boundary when the smoothness in the solution warrants it. For fifth order AFD-WENO, the two smaller stencils can be non-linearly hybridized along with the larger stencil. In that case, the second and fourth derivatives of the state can be obtained at the zone boundary when the smoothness in the solution warrants it. The process described here can be done for Adaptive Order and Multiresolution WENO interpolation.

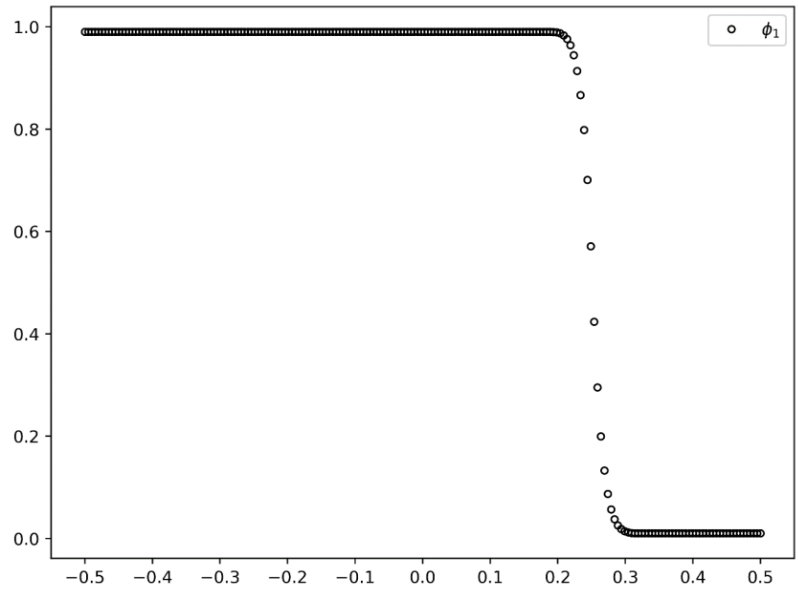 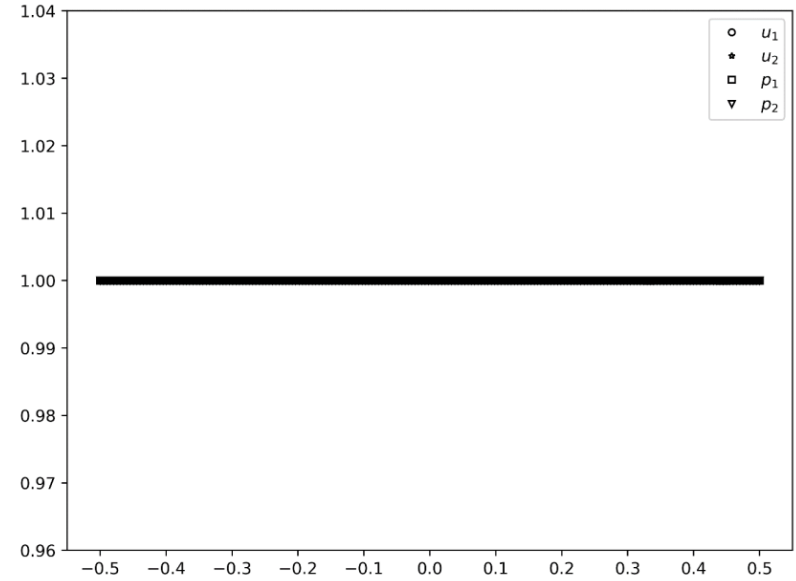

a) b)

*Fig. 4) Baer-Nunziato: Abgrall Problem using the 5$^{th}$ order accurate LLF-based AFD-WENO scheme with 200 zones. Fig. a) shows the solid volume fraction and Fig. b) shows the velocities and pressures for both the phases. The 7$^{th}$ and 9$^{th}$ order AFD-WENO schemes also show identical results and therefore they are not shown here.*

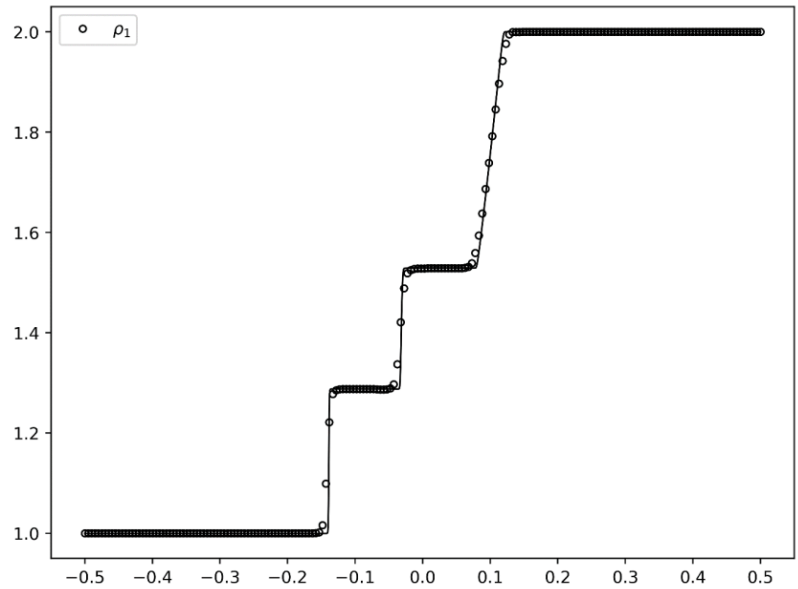 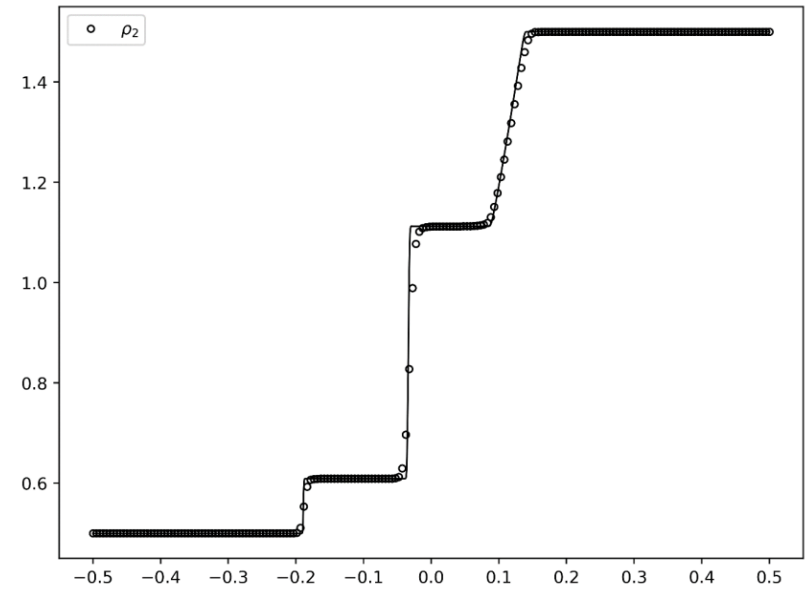

a) b)

*Fig. 5) Baer-Nunziato: Riemann Problem-1 using the 5$^{th}$ order accurate HLL-based AFD-WENO scheme with 200 zones. Fig. a) shows the density for the solid phase and Fig. b) shows the density for the gas phase. The 7$^{th}$ and 9$^{th}$ order AFD-WENO schemes also show identical results and therefore they are not shown here.*

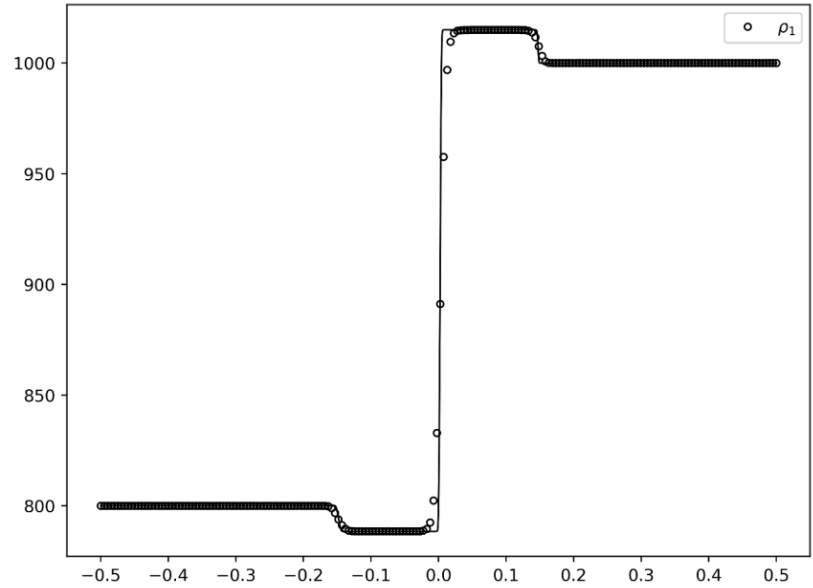 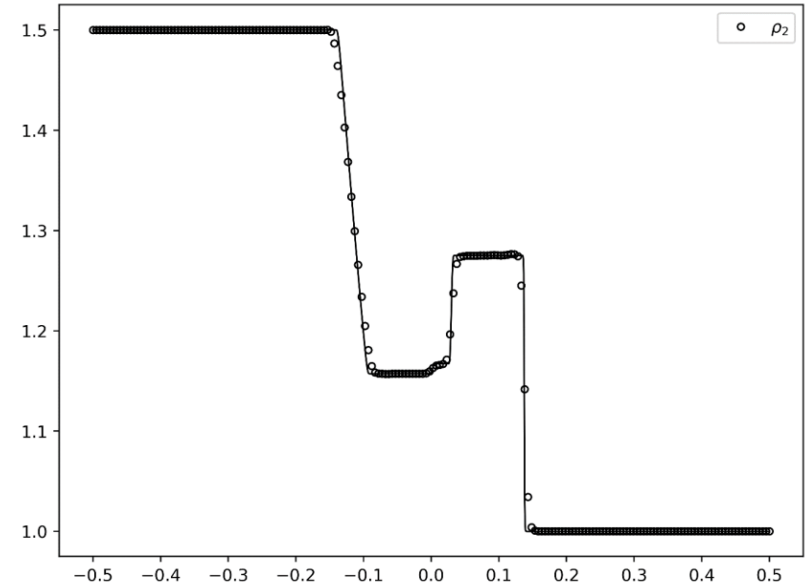

a)  b)

*Fig. 6) Baer-Nunziato: Riemann Problem-2 using the 7$^{th}$ order accurate HLL-based AFD-WENO scheme with 200 zones. Fig. a) shows the density for the solid phase and Fig. b) shows the density for the gas phase. The 5$^{th}$ and 9$^{th}$ order AFD-WENO schemes also show identical results and therefore they are not shown here.*

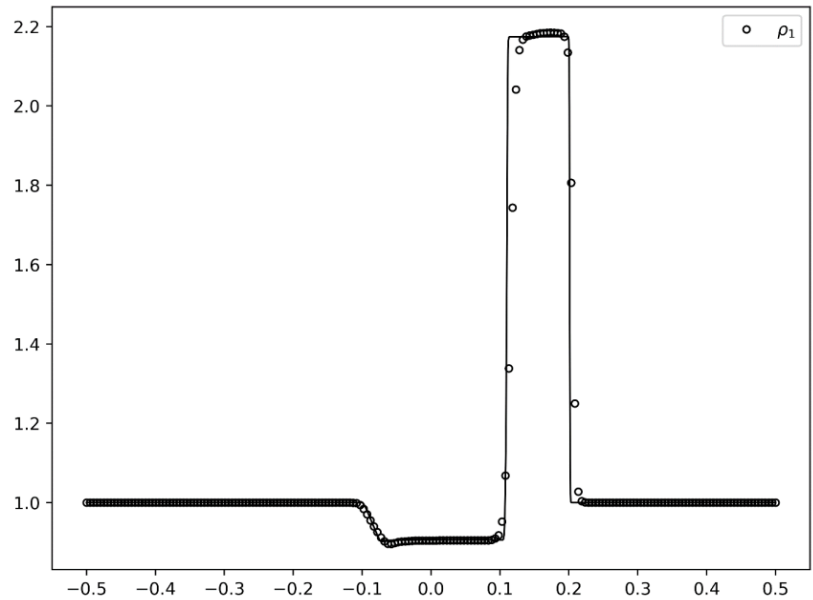 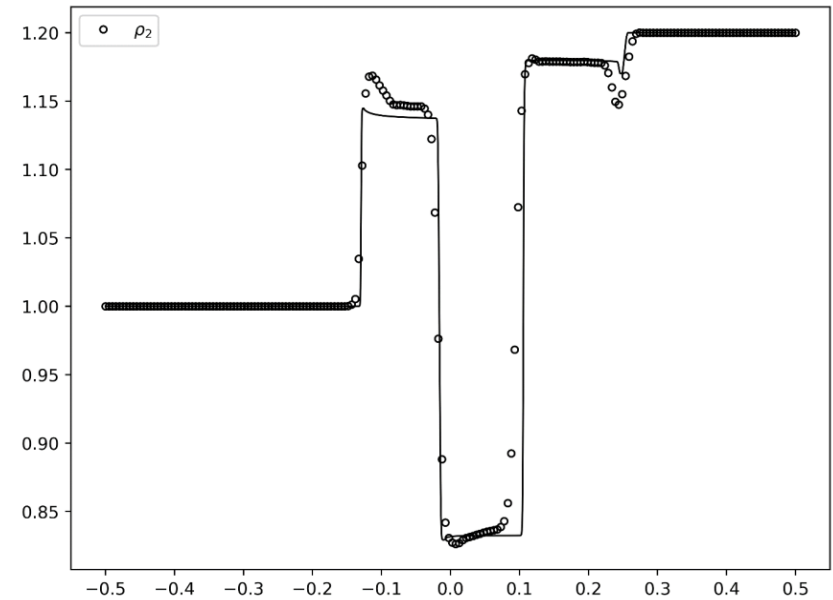

a)  b)

*Fig. 7) Baer-Nunziato: Riemann Problem-3 using the 9th order accurate HLL-based AFD-WENO scheme with 200 zones. Fig. a) shows the density for the solid phase and Fig. b) shows the density for the gas phase. The 5th and 7th order AFD-WENO schemes also show identical results and therefore they are not shown here.*

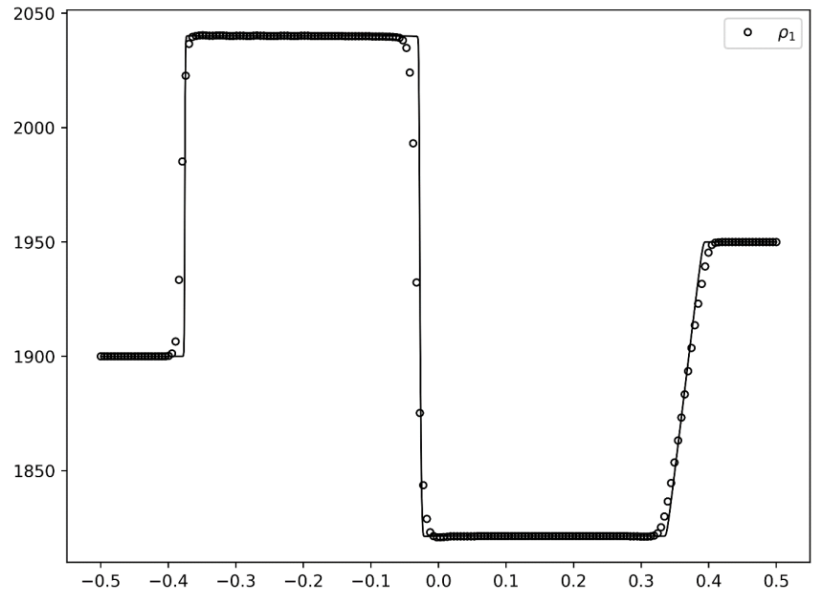
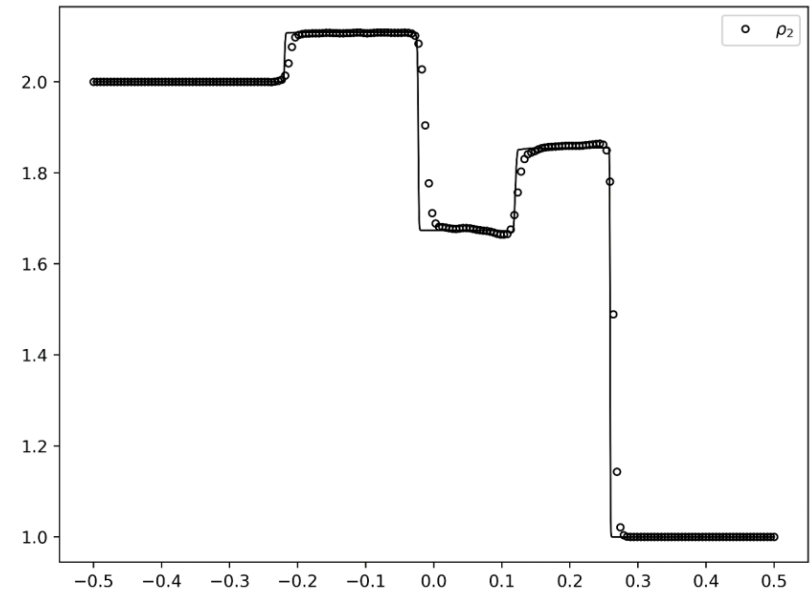

a)

b)

*Fig. 8) Baer-Nunziato: Riemann Problem-4 using the 5th order accurate LLF-based AFD-WENO scheme with 200 zones. Fig. a) shows the density for the solid phase and Fig. b) shows the density for the gas phase. The 7th and 9th order AFD-WENO schemes also show identical results and therefore they are not shown here.*

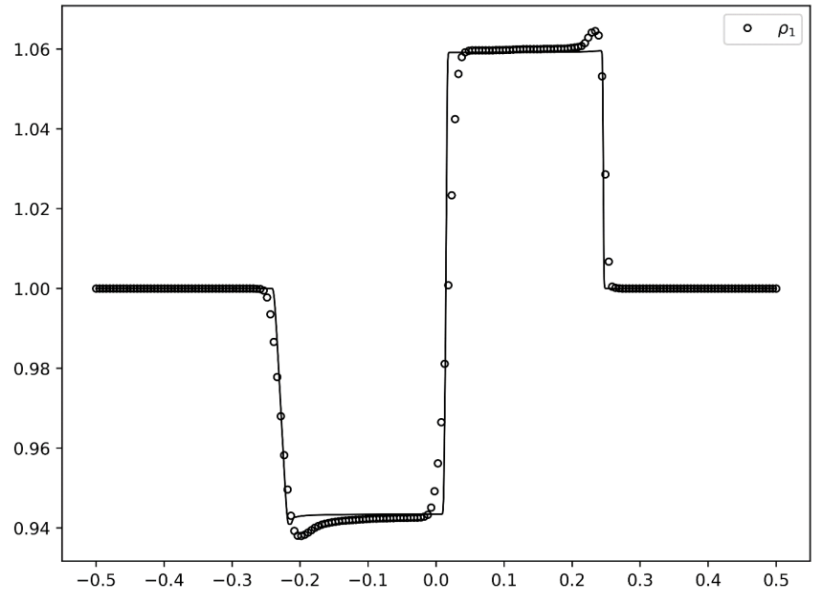 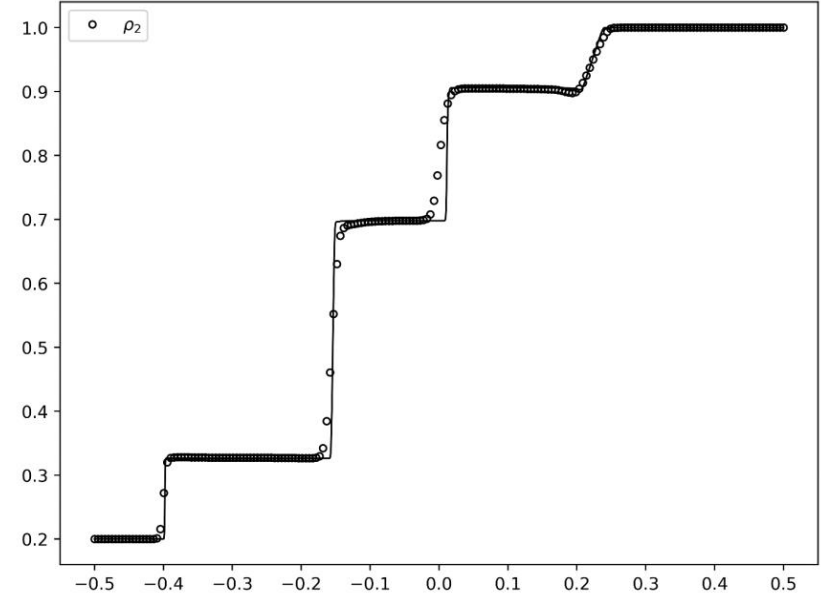

a)    b)

*Fig. 9) Baer-Nunziato: Riemann Problem-5 using the 7th order accurate LLF-based AFD-WENO scheme with 200 zones. Fig. a) shows the density for the solid phase and Fig. b) shows the density for the gas phase. We have used flattener for this problem. The 5th and 9th order AFD-WENO schemes also show identical results and therefore they are not shown here.*

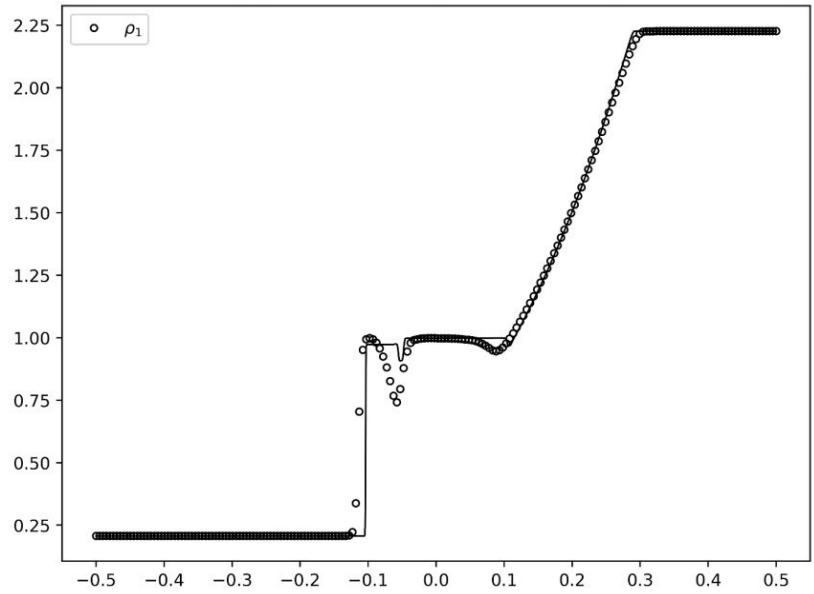 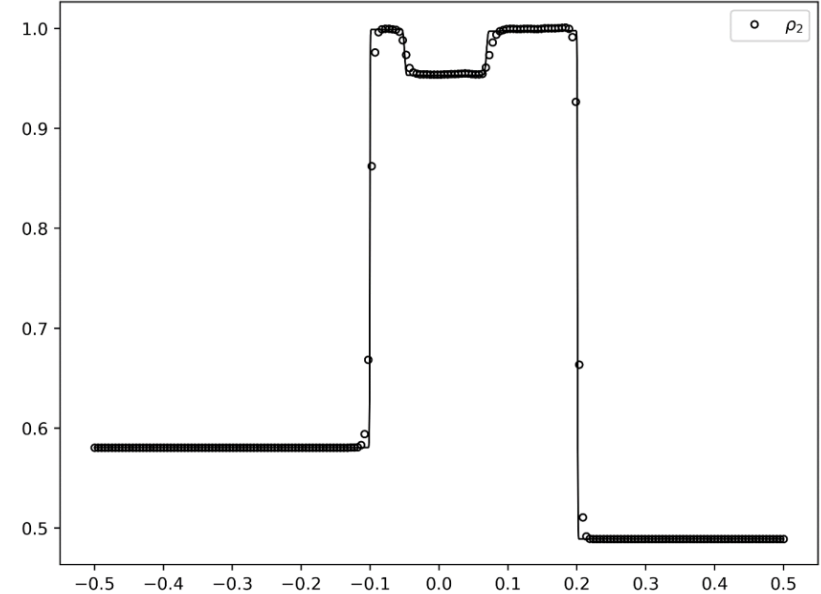

a) b)

*Fig. 10) Baer-Nunziato: Riemann Problem-6 using the 9th order accurate LLF-based AFD-WENO scheme with 200 zones. Fig. a) shows the density for the solid phase and Fig. b) shows the density for the gas phase. The 5th and 7th order AFD-WENO schemes also show identical results and therefore they are not shown here.*

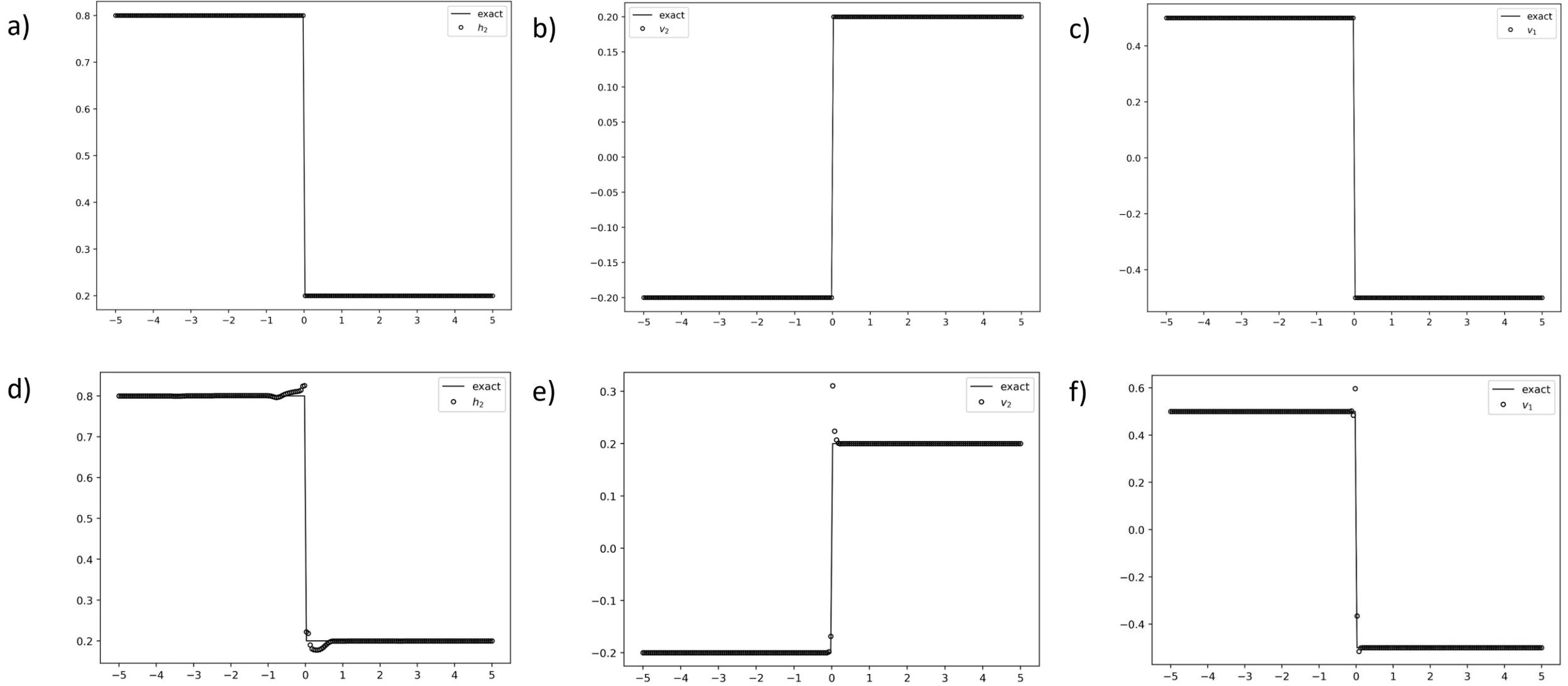

*Fig. 11) Two-Layer Shallow water: Riemann problem-1 using the 5$^{th}$ order accurate HLLI-based AFD-WENO scheme with 200 zones. The 7$^{th}$ and 9$^{th}$ order AFD-WENO schemes also show identical results and therefore they are not shown here. The top row shows the results from the scheme developed here. The bottom row of three figures shows the results from the same scheme when central differences are used for the higher order derivatives. We see that using just the central differences results in unphysical oscillations. Similar oscillations would be found if central differences were used for a 7$^{th}$ or 9$^{th}$ order scheme.*

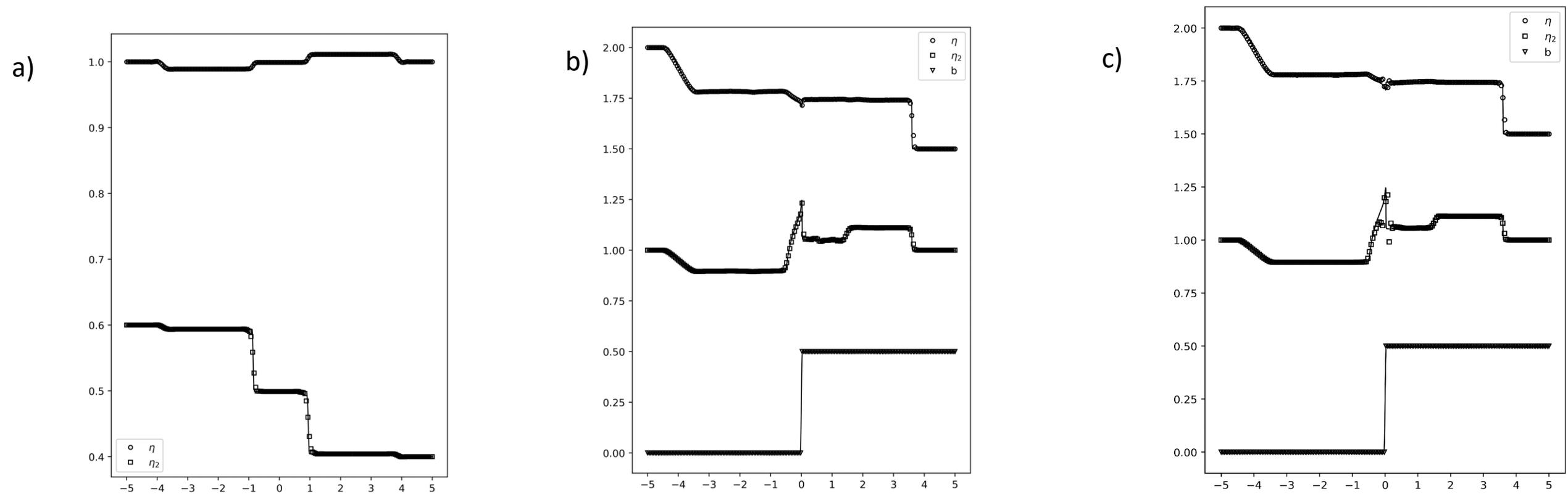

*Fig. 12) Two-Layer Shallow water: Fig. a) shows the results of Riemann problem-2 using the 7th order accurate HLLI-based AFD-WENO scheme and Fig. b) shows the results of Riemann problem-3 using the 9th order accurate HLLI-based AFD-WENO scheme with 200 zones. The 5th and 9th order scheme for the Riemann Problem-2, and the 5th and 7th order scheme for the Riemann Problem-3 also show identical results and therefore they are not shown here. Fig. c) shows the results of Riemann problem-3 using the 9th order scheme when only central differences are used for the higher order derivatives. We clearly see the formation of unphysical oscillations when only central differences are used. Similar oscillations would be found if central differences were used for a 5th or 7th order scheme.*

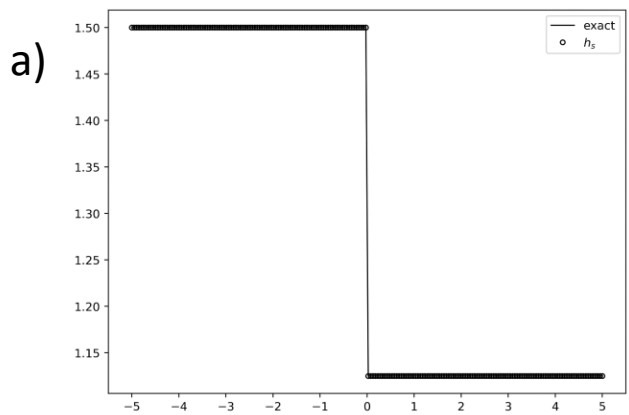
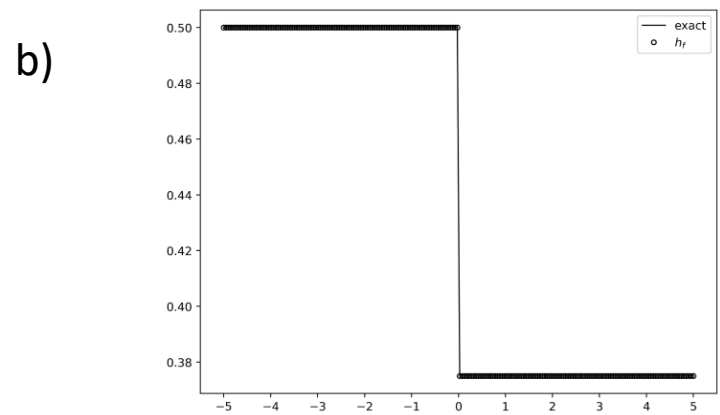
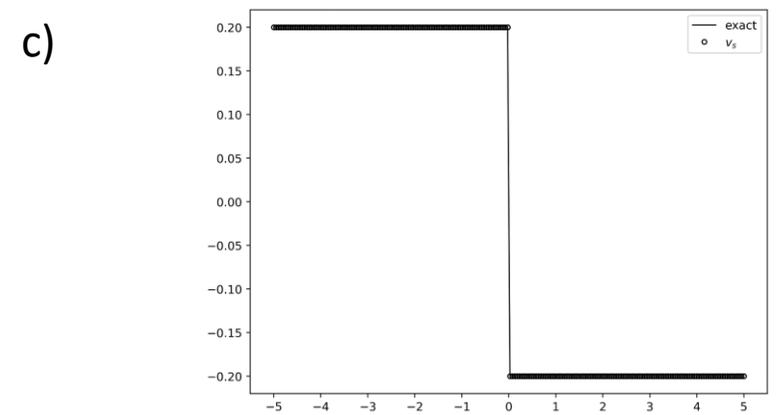
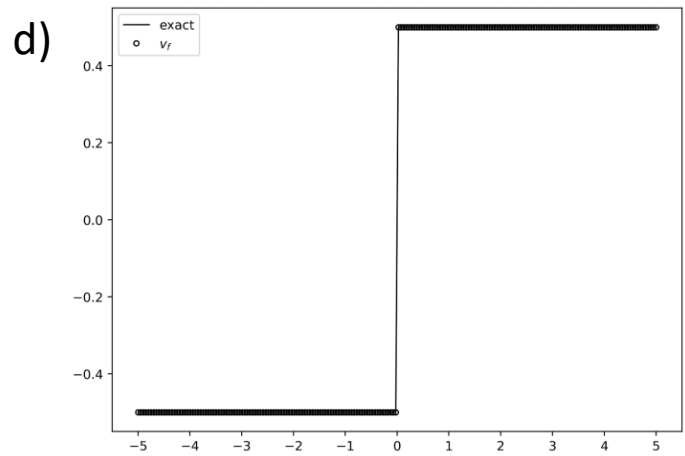
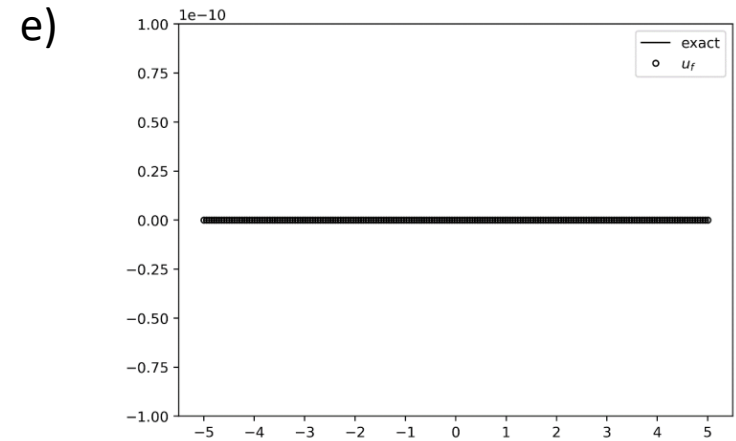

*Fig. 13a-e) Debris Flow: a to d show the jump in the linearly degenerate fields using the 5$^{th}$ order accurate HLLI-based AFD-WENO scheme with 200 zones. e shows the velocity. The 7$^{th}$ and 9$^{th}$ order AFD-WENO schemes also show identical results and therefore they are not shown here.*

f) 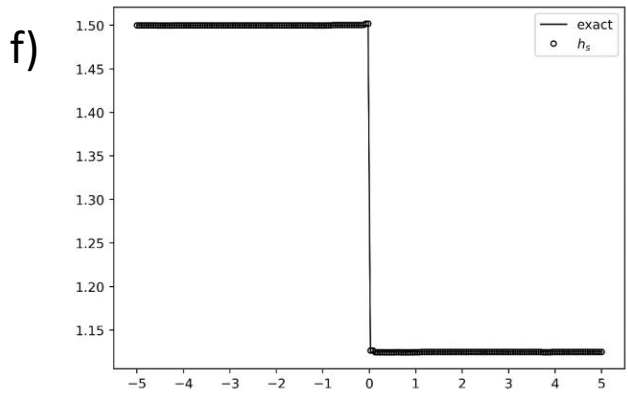
g) 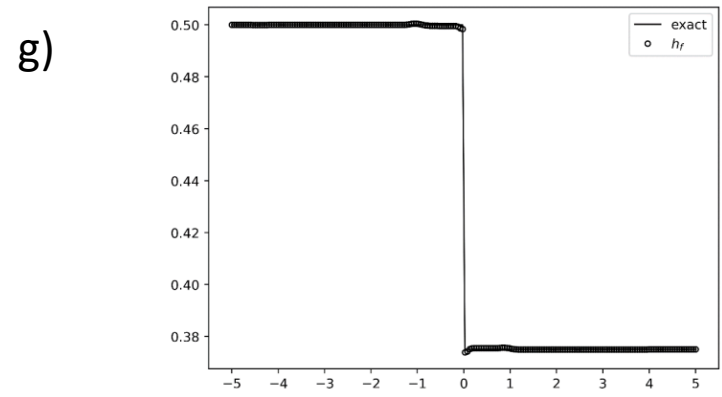
h) 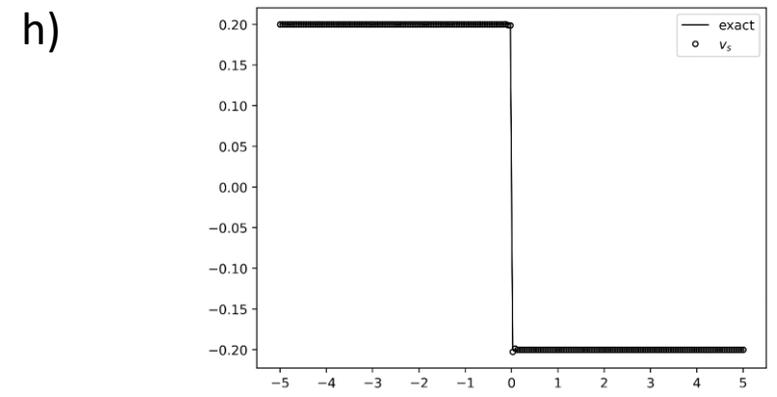

i) 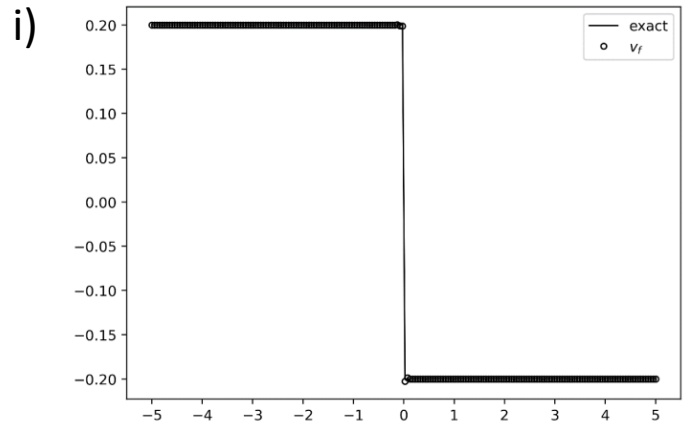
j) 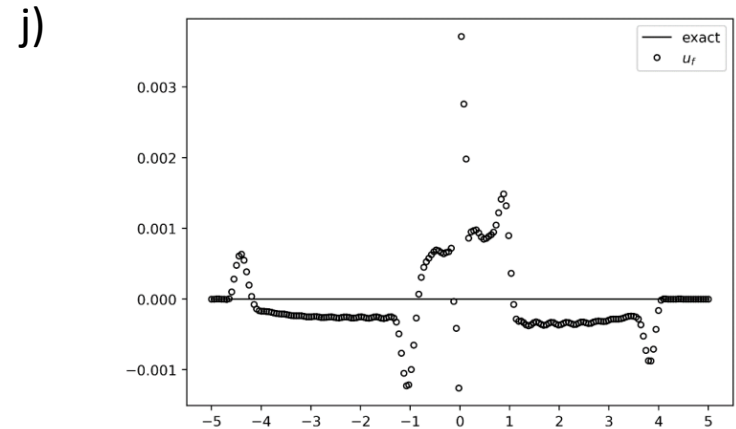

*Fig. 13f-j) Debris Flow: f to i show the jump in the linearly degenerate fields using the 5$^{th}$ order accurate HLLI-based AFD-WENO scheme with 200 zones when only central differences are used for the higher order derivatives. j shows the velocity. The 7$^{th}$ and 9$^{th}$ order AFD-WENO schemes also show identical results and therefore they are not shown here.*

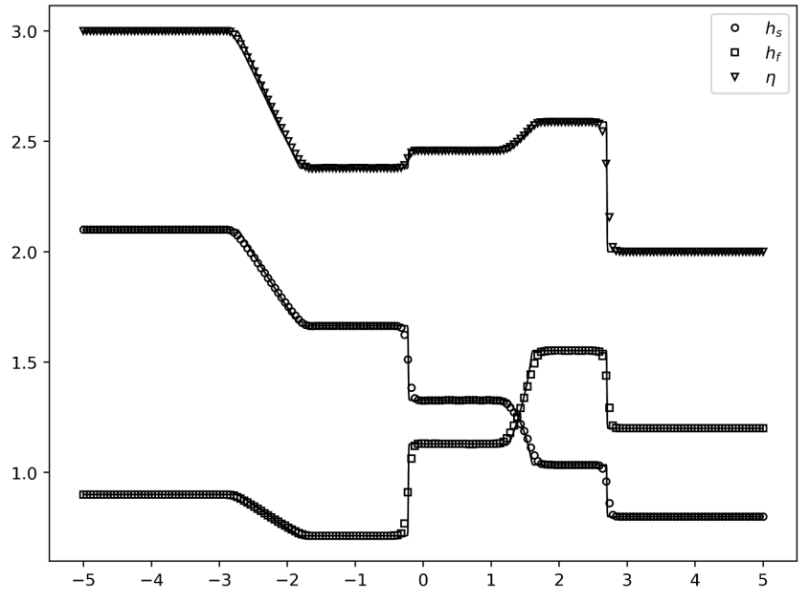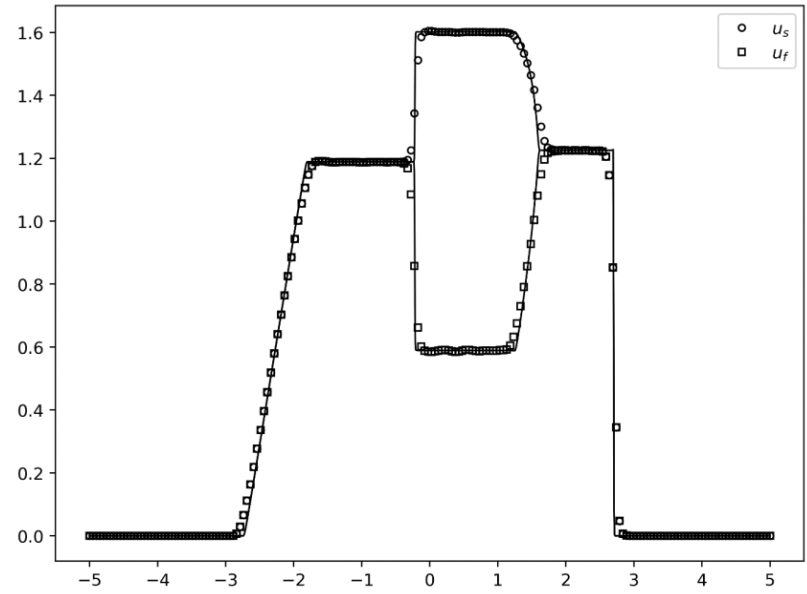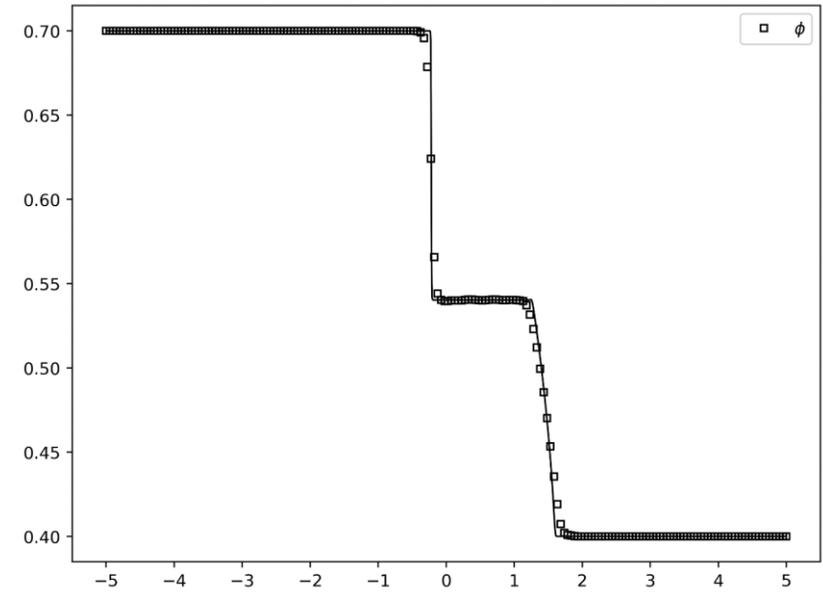

*Fig. 14) Debris Flow: Riemann problem-2 using the 7th order accurate HLLI-based AFD-WENO scheme with 200 zones. The 5th and 9th order AFD-WENO schemes also show identical results and therefore they are not shown here.*

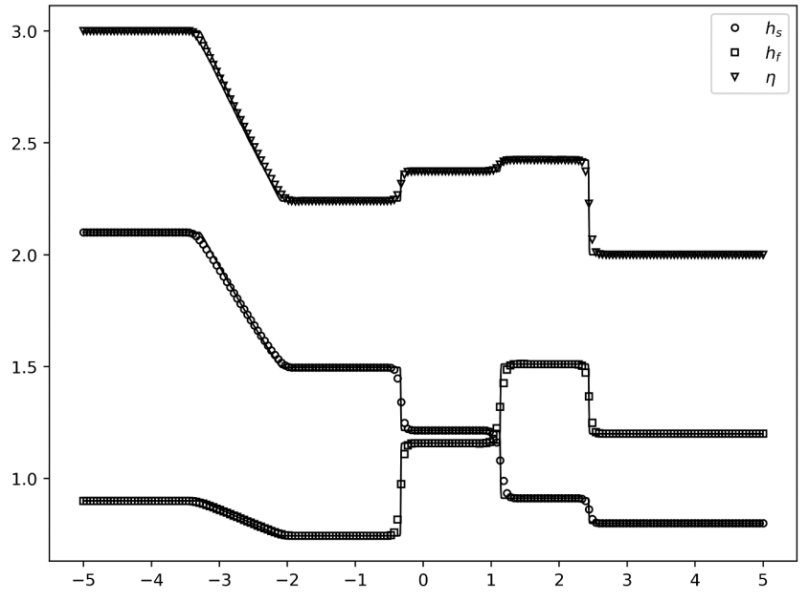 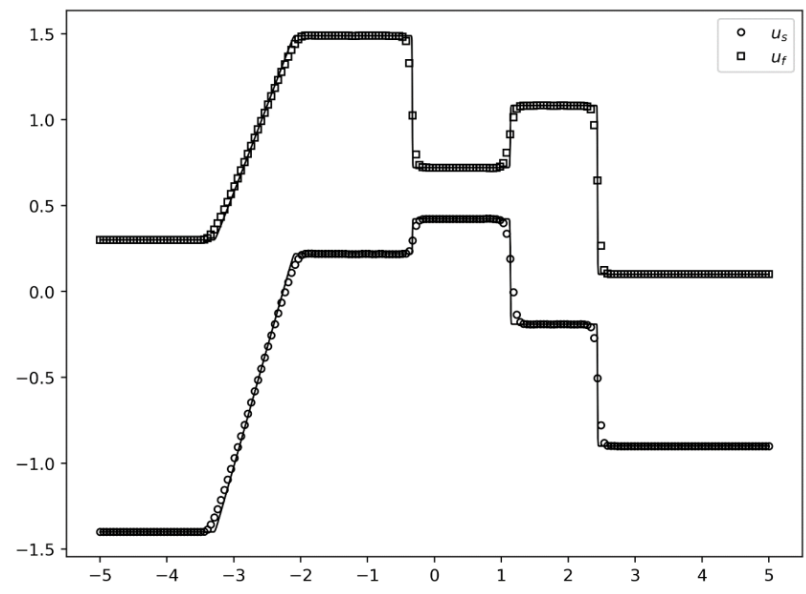 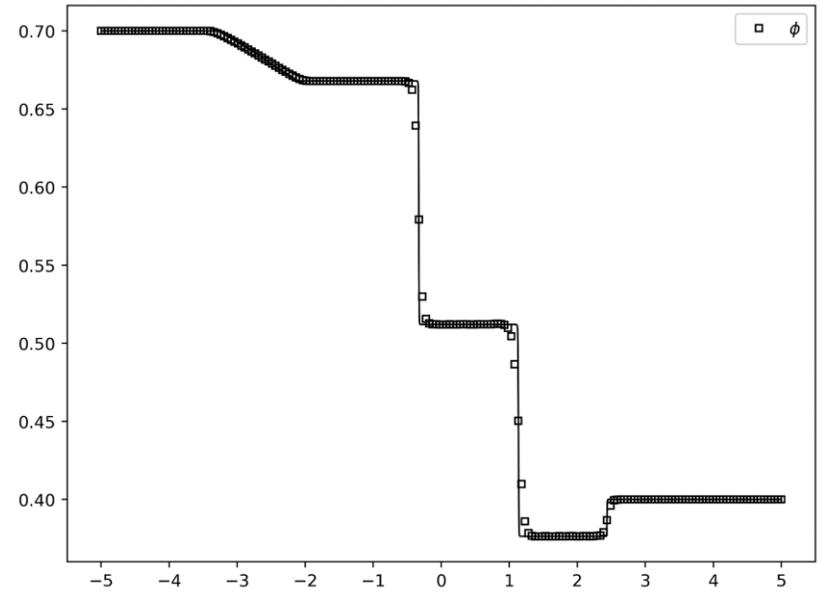

a)          b)          c)

*Fig. 15) Debris Flow: Riemann problem-3 using the 9th order accurate HLLI-based AFD-WENO scheme with 200 zones. The 5th and 7th order AFD-WENO schemes also show identical results and therefore they are not shown here.*

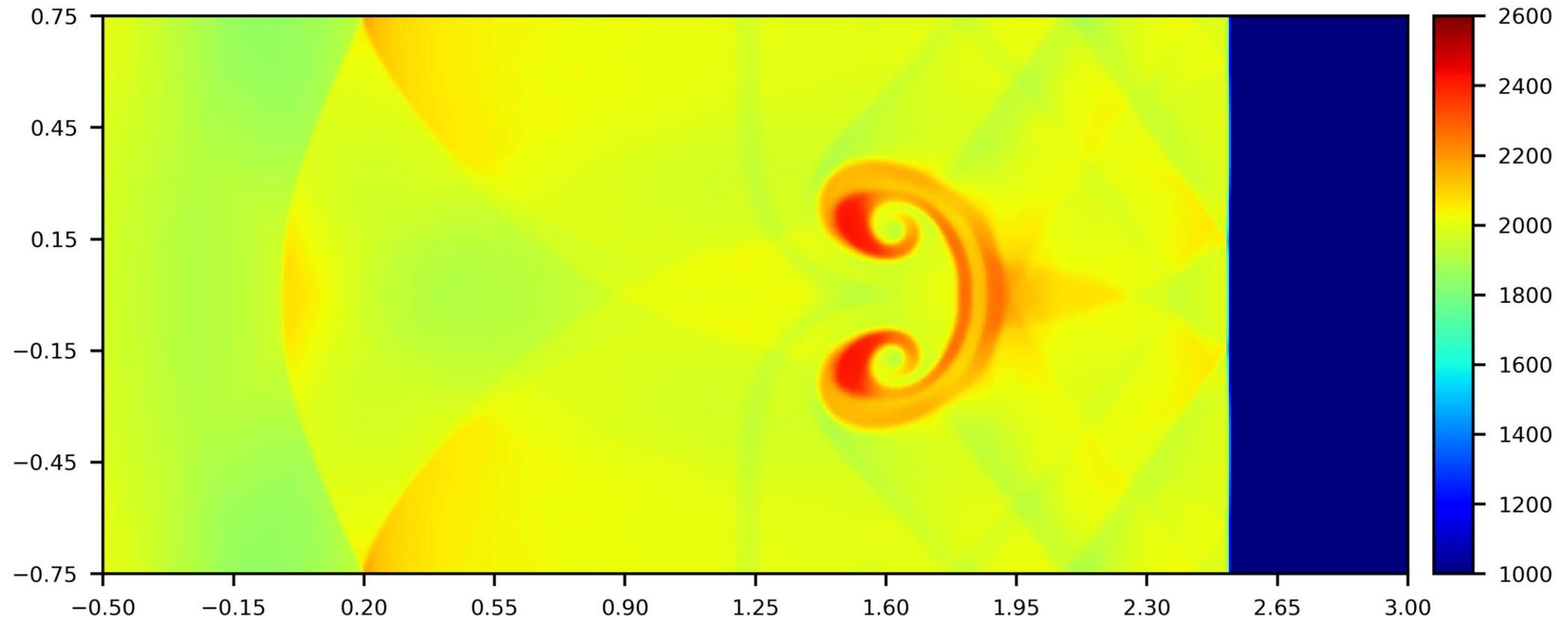

*Fig. 16a) Baer-Nunziato model: Shock-bubble interaction problem using the 3rd order accurate HLL-based AFD WENO-AO-3 scheme with 700×300 zones. The solid density profiles have been shown.*

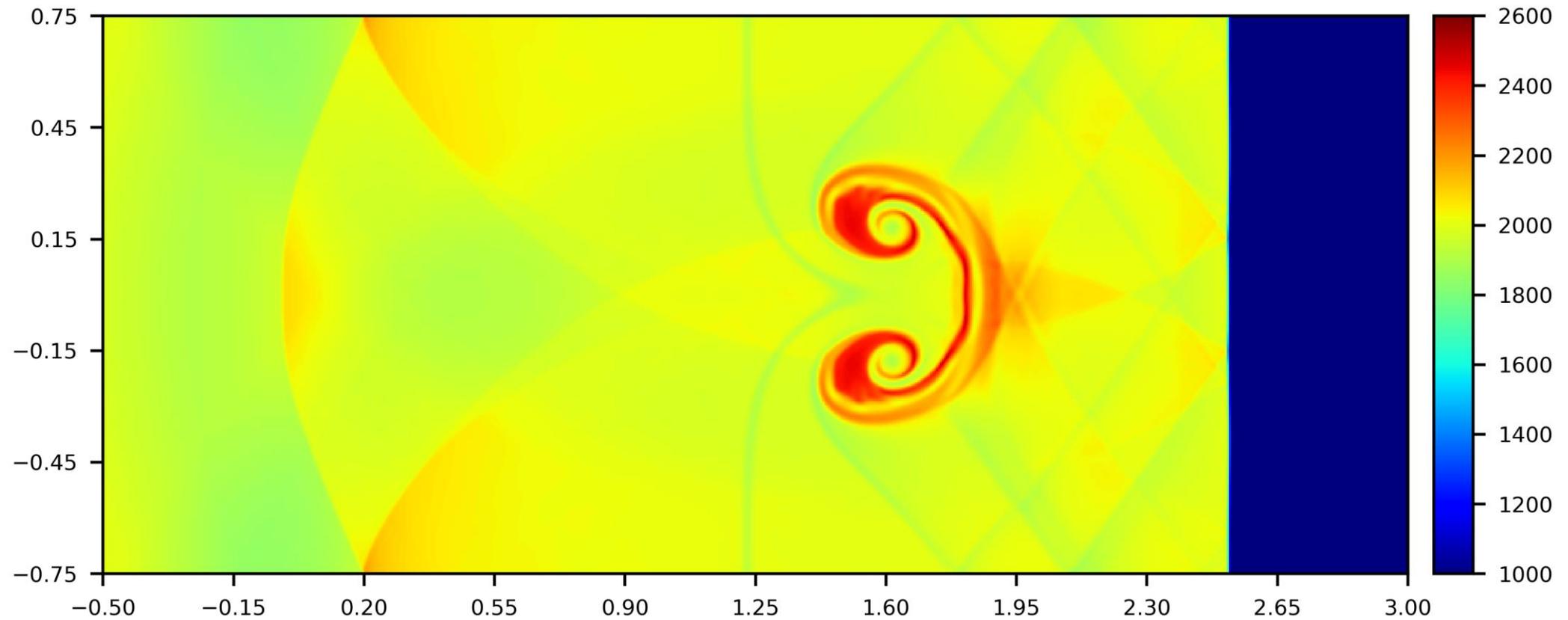

*Fig. 16b) Baer-Nunziato model: Shock-bubble interaction problem using the 5$^{th}$ order accurate HLL-based AFD WENO-AO-(5,3) scheme with 700 ×300 zones. The solid density profiles have been shown.*

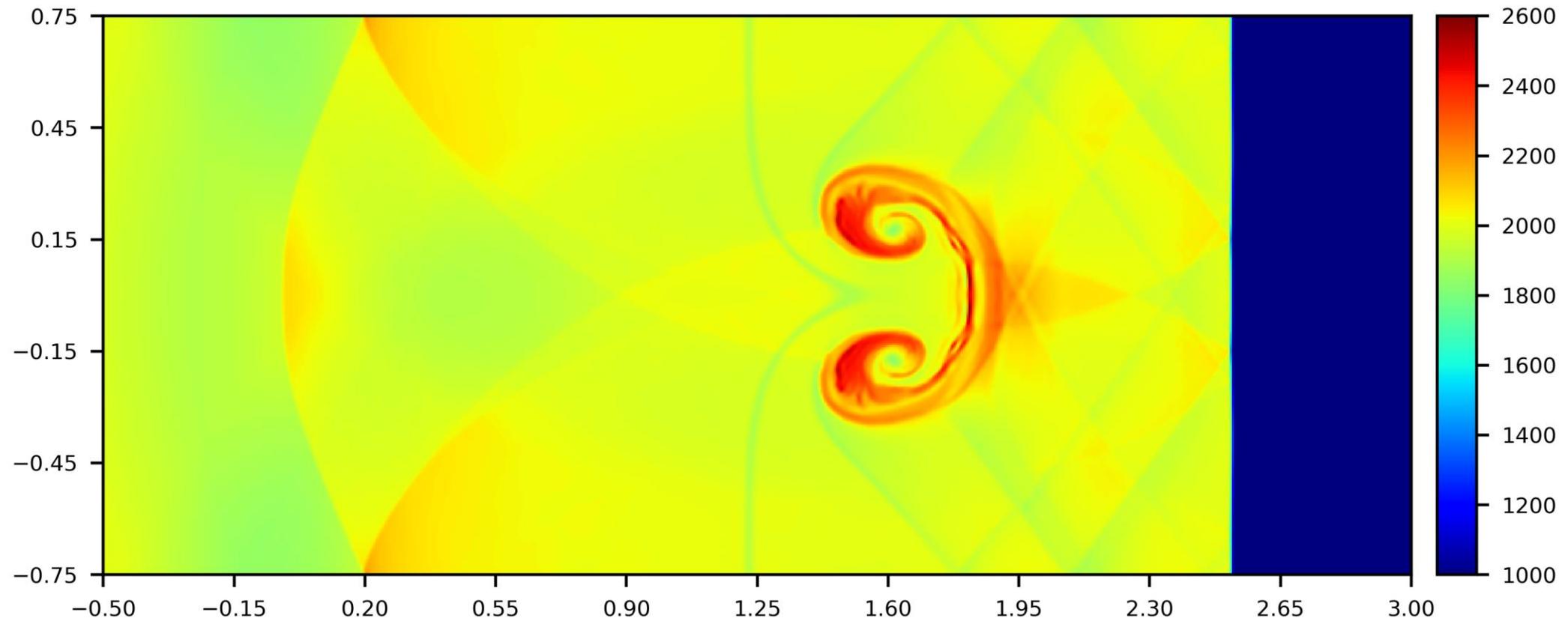

*Fig. 16c) Baer-Nunziato model: Shock-bubble interaction problem using the 7$^{th}$ order accurate HLL-based AFD WENO-AO-(7,5,3) scheme with 700×300 zones. The solid density profiles have been shown. The 9$^{th}$ order scheme also shows identical result and therefore it is not shown here.*

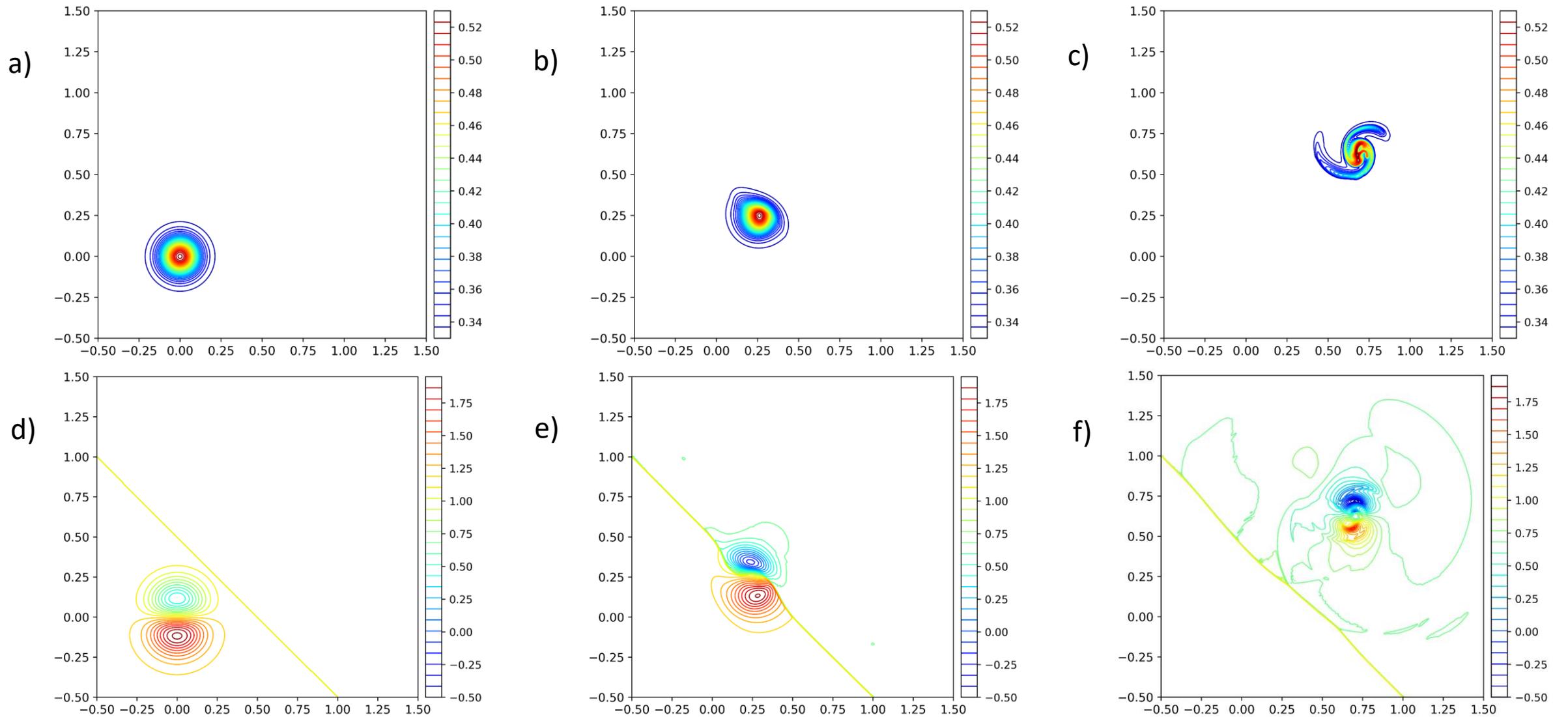

*Fig. 17) Baer-Nunziato: Shock-Vortex Interaction using the 5$^{th}$ order accurate HLL-based AFD WENO-AO-(5,3) scheme with 600×600 zones at time levels t=0.0, 0.23 and 0.84. Figs. 17a, 17b and 17c show solid volume fraction at times t=0.0, 0.23, 0.84. Figs. 17d, 17e, 17f show solid x-velocity at times t=0.0, 0.23, 0.84. For the solid volume fraction, 30 contours were fit between a range of 0.33 and 0.530. For the solid x-velocity, 30 contours were fit between a range of -0.5 and 1.95.*

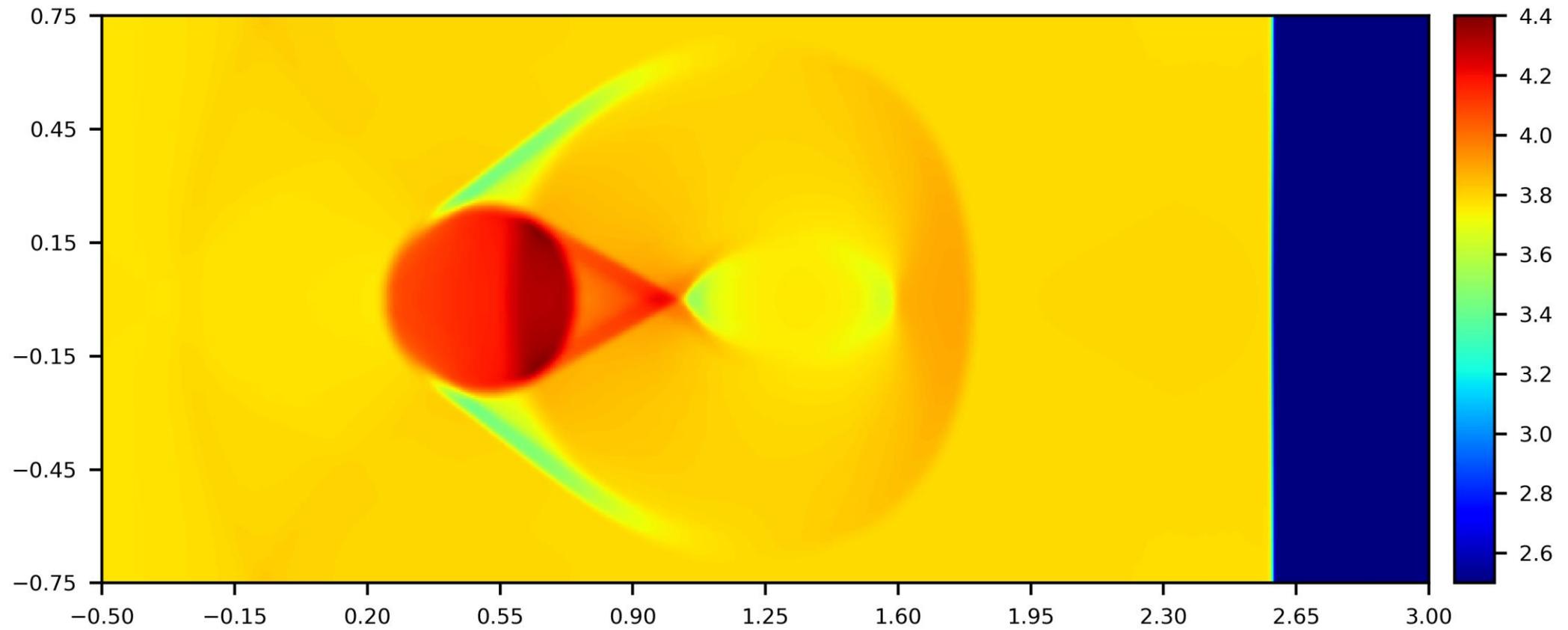

*Fig. 18a) Two-Layer Shallow water: Shock-bubble interaction problem using the 3rd order accurate HLL-based AFD WENO-AO-3 scheme with 700 ×300 zones. The solid density profiles have been shown.*

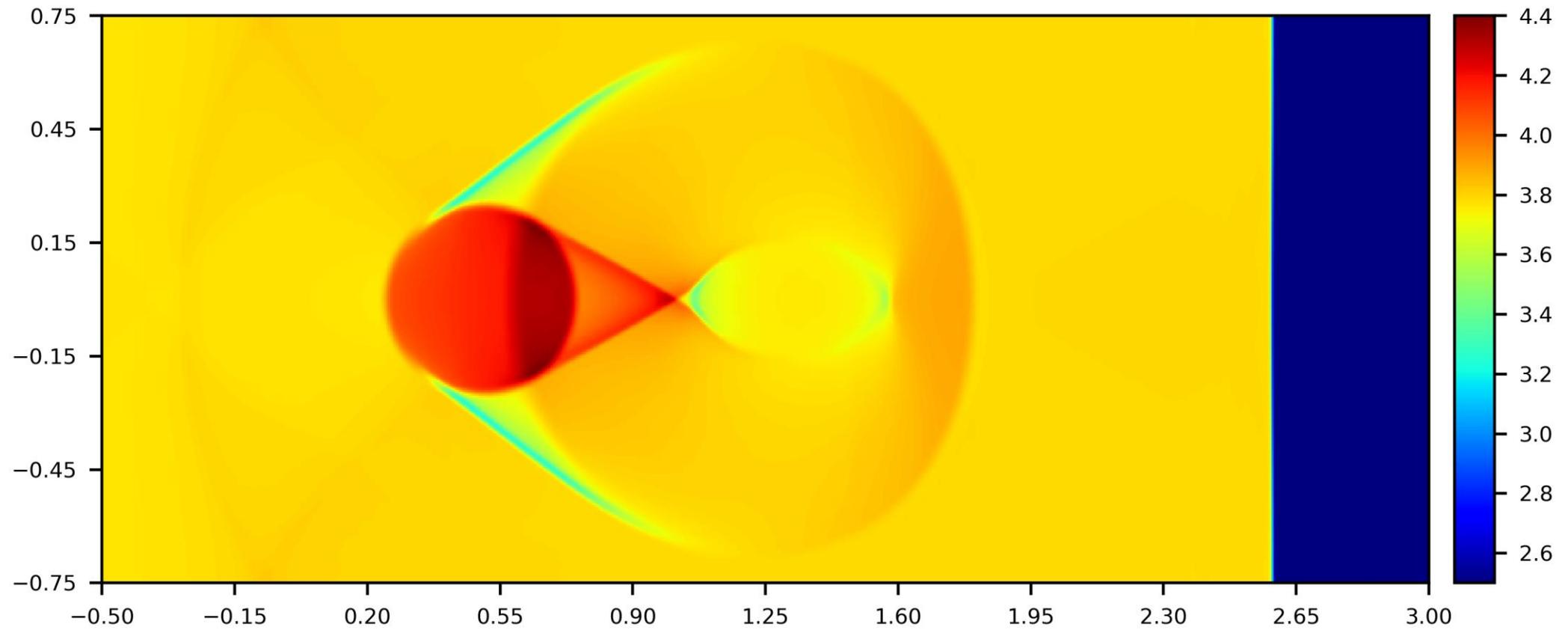

*Fig. 18b) Two-Layer Shallow water: Shock-bubble interaction problem using the 5th order accurate HLL-based AFD WENO-AO-(5,3) scheme with 700 ×300 zones. The solid density profiles have been shown.*

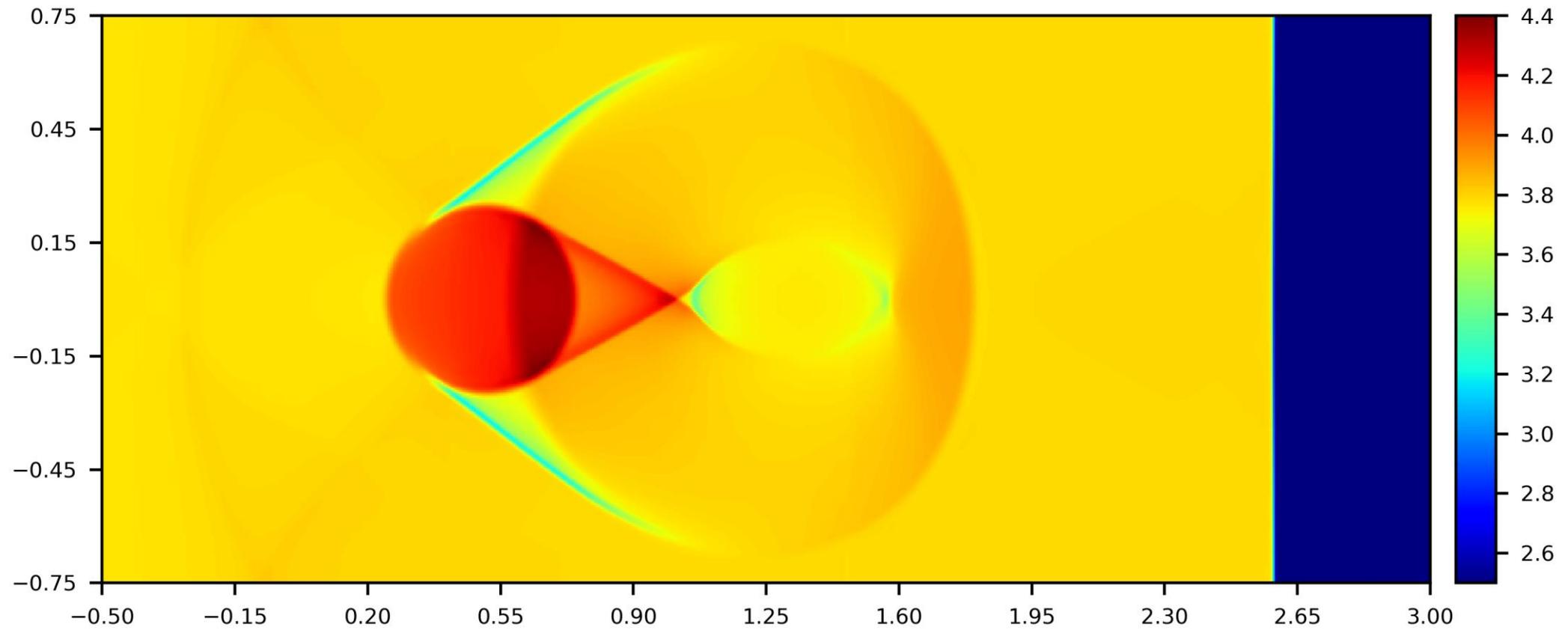

*Fig. 18c) Two-Layer Shallow water: Shock-bubble interaction problem using the 7th order accurate HLL-based AFD WENO-AO-(7,5,3) scheme with 700×300 zones. The solid density profiles have been shown. The 9th order scheme also shows identical result and therefore it is not shown here.*

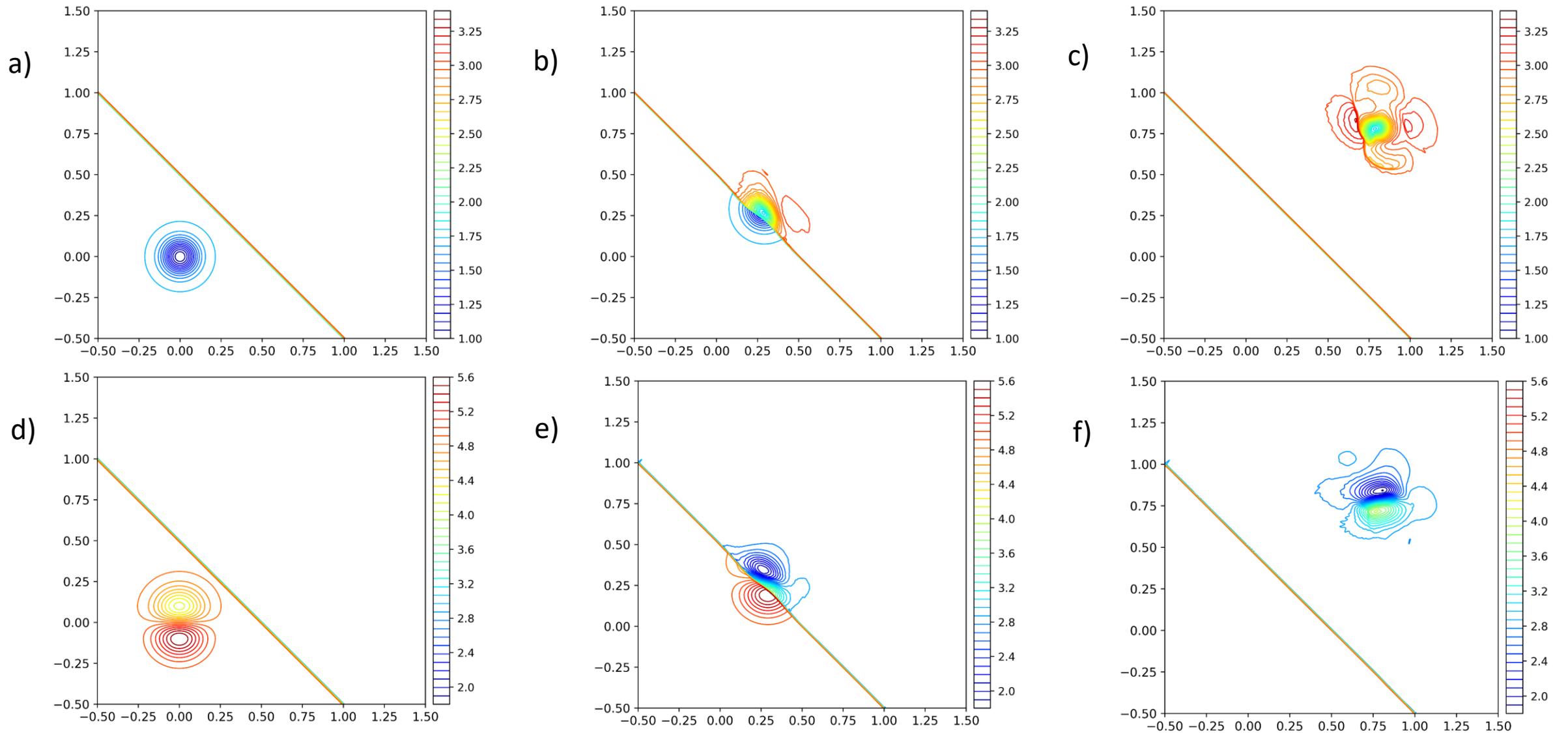

Fig. 19) Two-Layer Shallow water: Shock-Vortex Interaction using the 7$^{th}$ order accurate HLL-based AFD WENO-AO-(7,3) scheme with 600×600 zones at time levels t=0.0, 0.06 and 0.24. Figs. 19a, 19b and 19c show height of the upper fluid at times t=0.0, 0.06, 0.24. Figs. 19d, 19e, 19f show x-velocity of the upper fluid at times t=0.0, 0.06, 0.24. For the height, 40 contours were fit between a range of 1.0 and 3.4. For the velocity, 40 contours were fit between a range of 1.8 and 5.6.

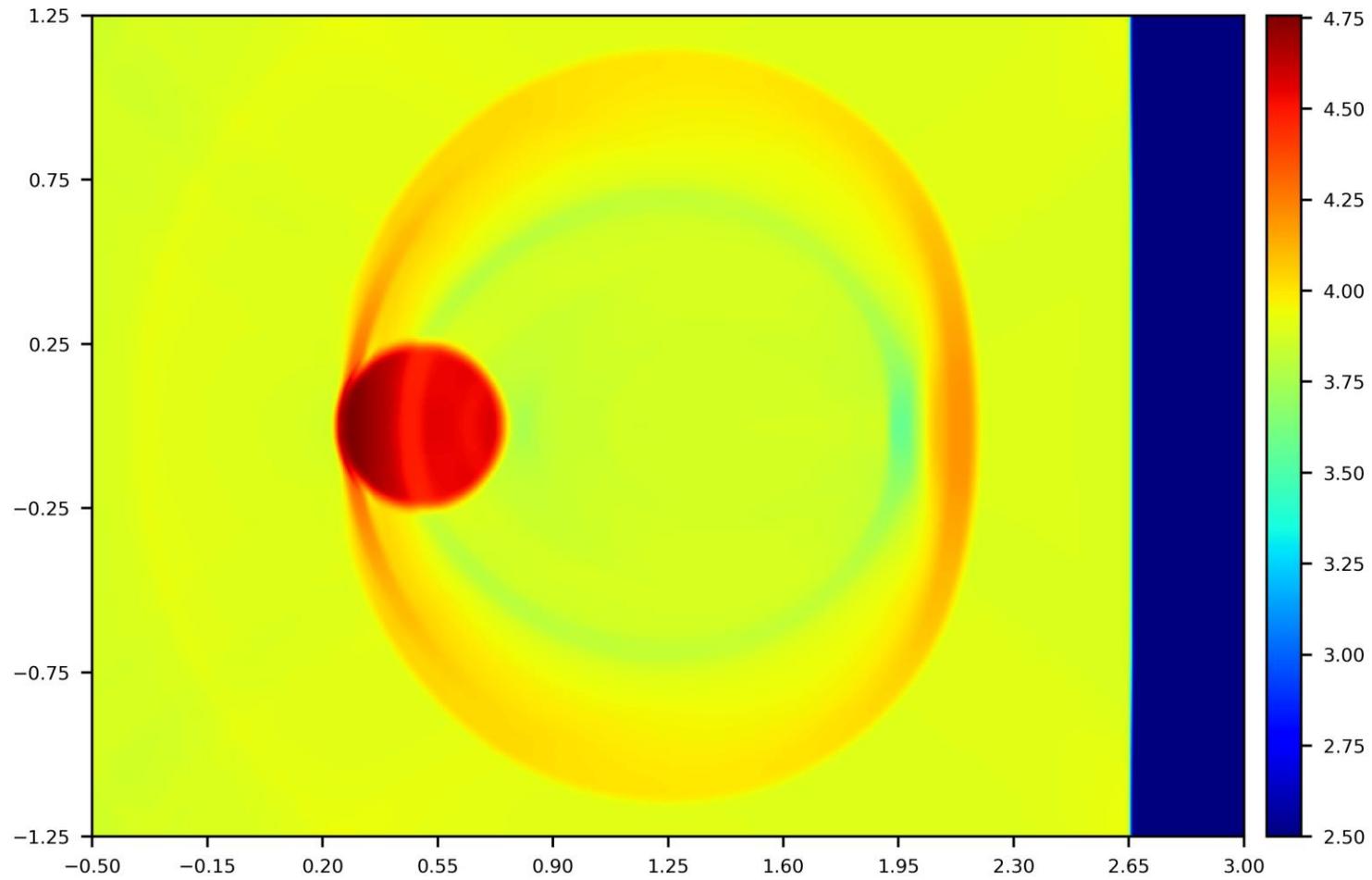

*Fig. 20a) Debris Flow: Shock-bubble interaction problem using the 3$^{rd}$ order accurate HLL-based AFD WENO-AO-3 scheme with 700 ×500 zones. The solid density profiles have been shown.*

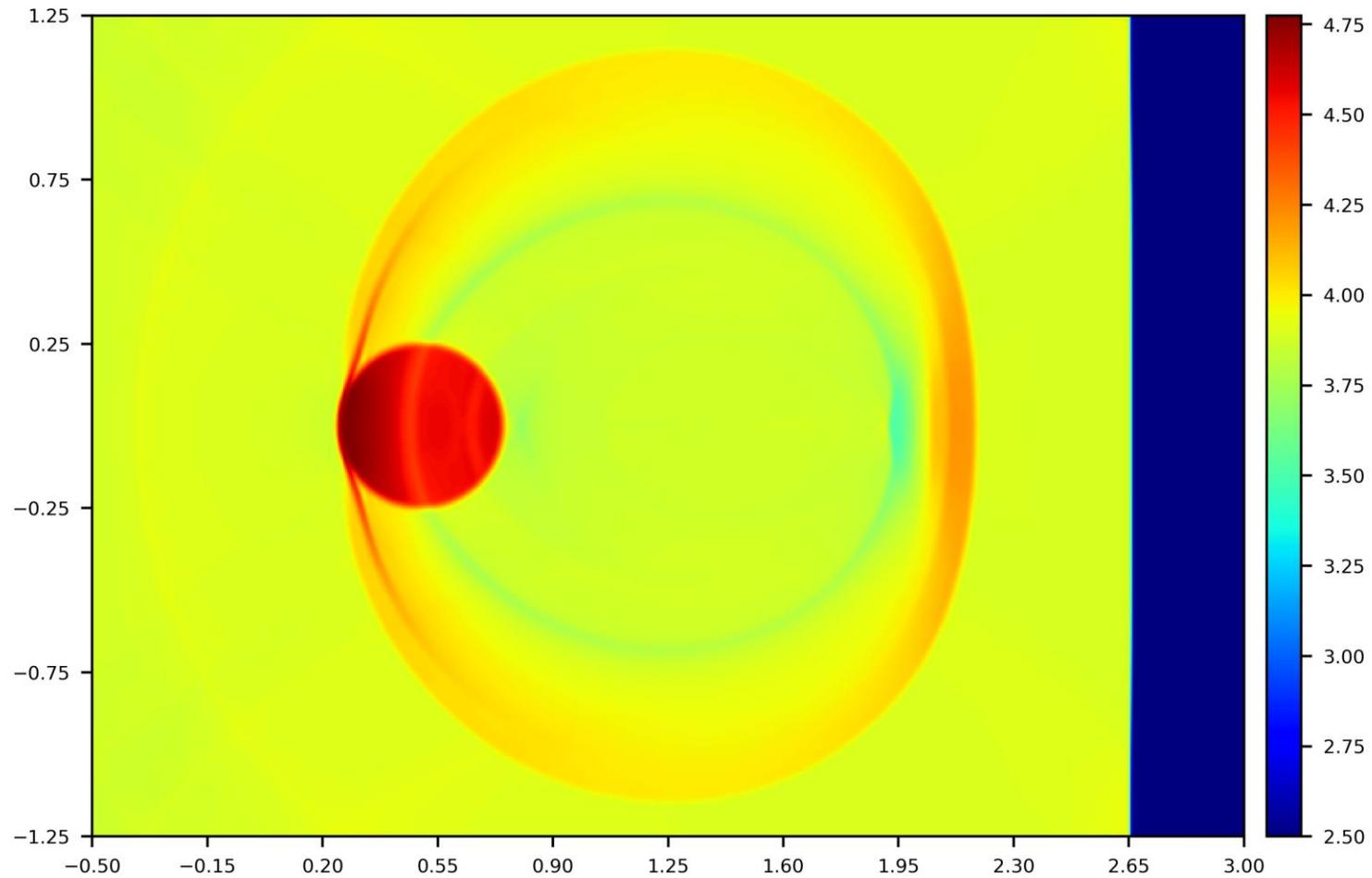

*Fig. 20b) Debris Flow: Shock-bubble interaction problem using the 5$^{th}$ order accurate HLL-based AFD WENO-AO-(5,3) scheme with 700 ×500 zones. The solid density profiles have been shown.*

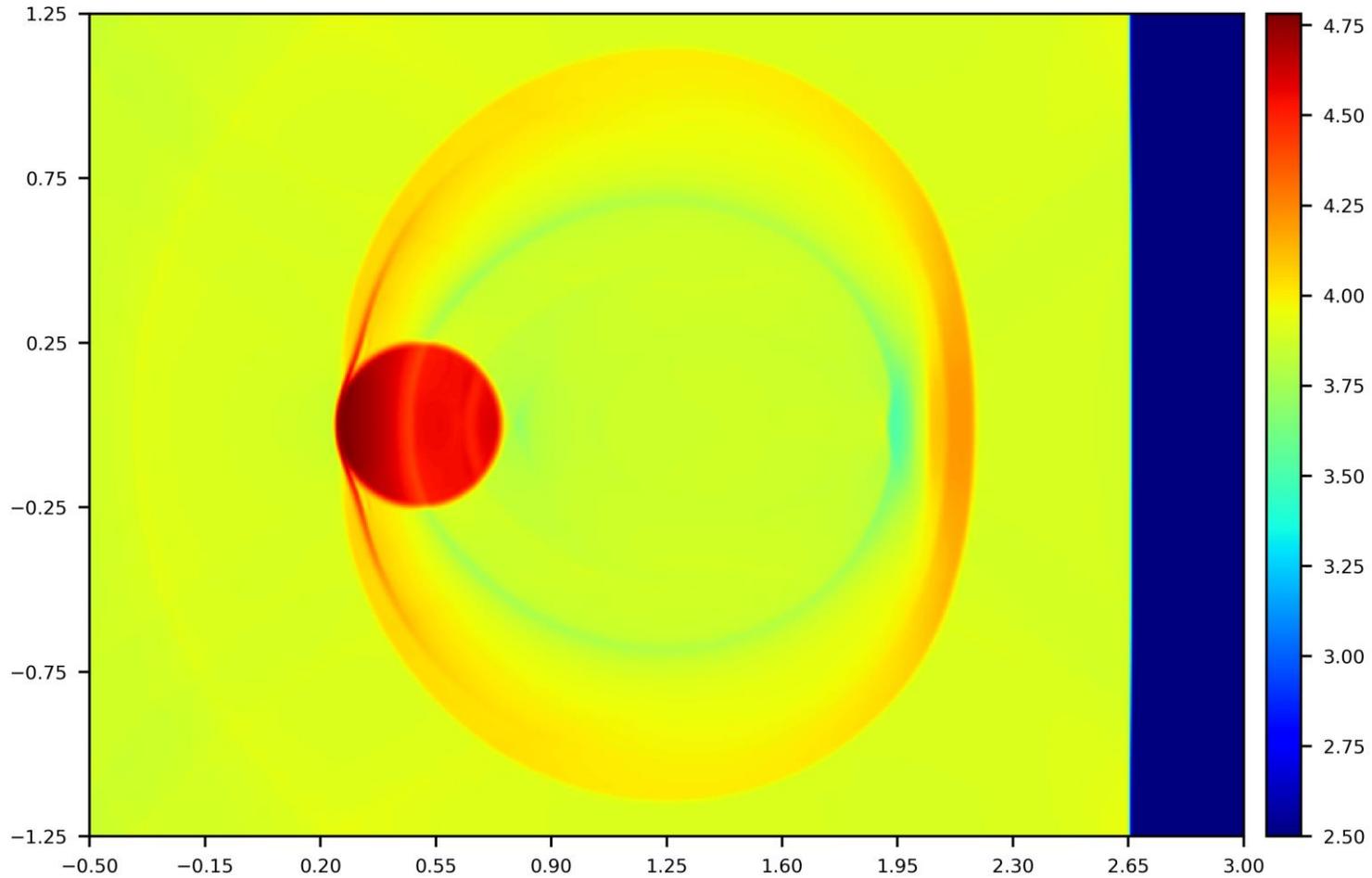

*Fig. 20c) Debris Flow: Shock-bubble interaction problem using the 7th order accurate HLL-based AFD WENO-AO-(7,5,3) scheme with 700×500 zones. The solid density profiles have been shown. The 9th order accurate scheme also shows identical results and therefore it is not shown here.*

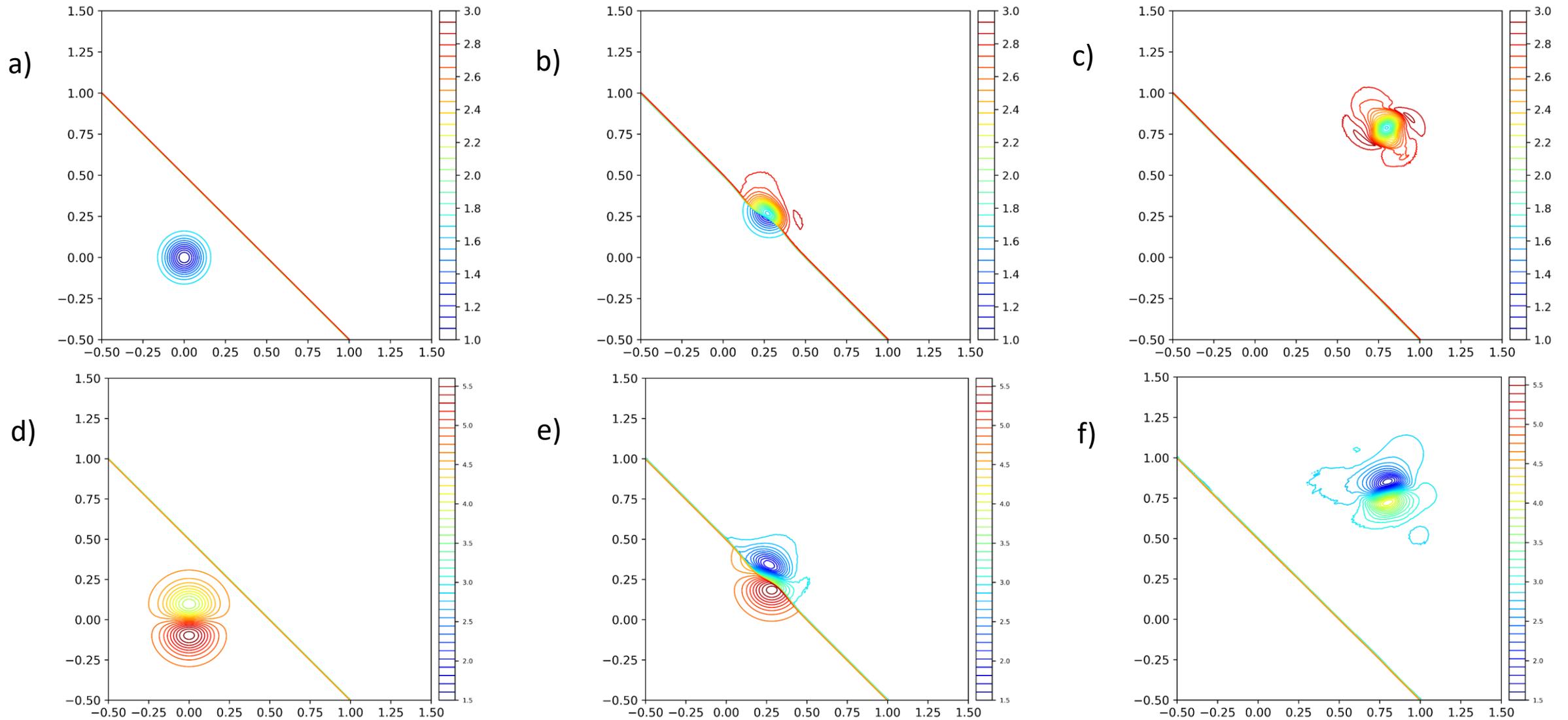

Fig. 21) Debris Flow: Shock-Vortex Interaction using the 9$^{th}$ order accurate HLL-based AFD WENO-AO-(9,3) scheme with 600×600 zones at time levels t=0.0, 0.06 and 0.24. Figs. 21a, 21b and 21c show the solid height at times t=0.0, 0.06, 0.24. Figs. 21d, 21e, 21f show the solid x-velocity at times t=0.0, 0.06, 0.24. For the height, 30 contours were fit between a range of 1.0 and 3.0. For the velocity, 30 contours were fit between a range of 1.5 and 5.6.

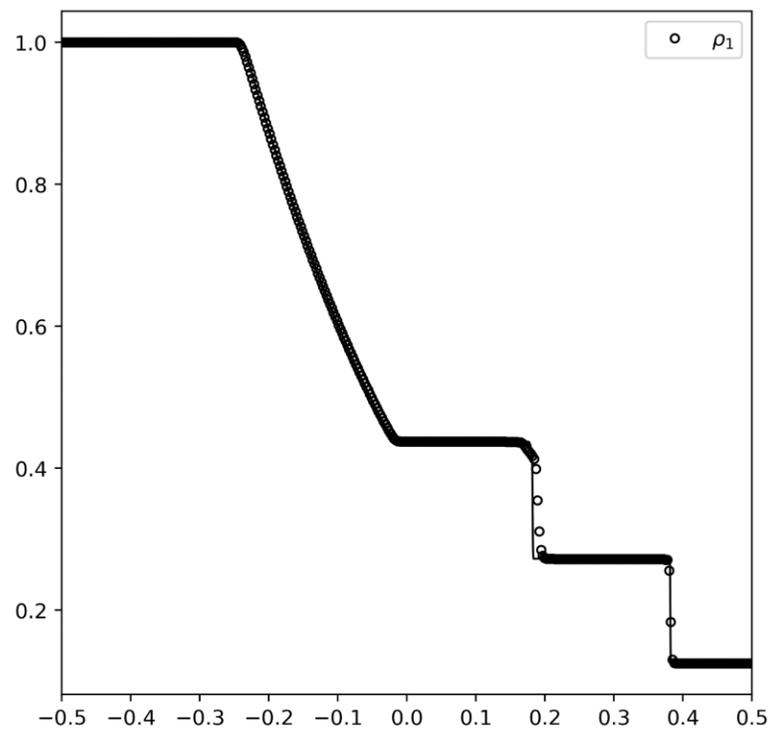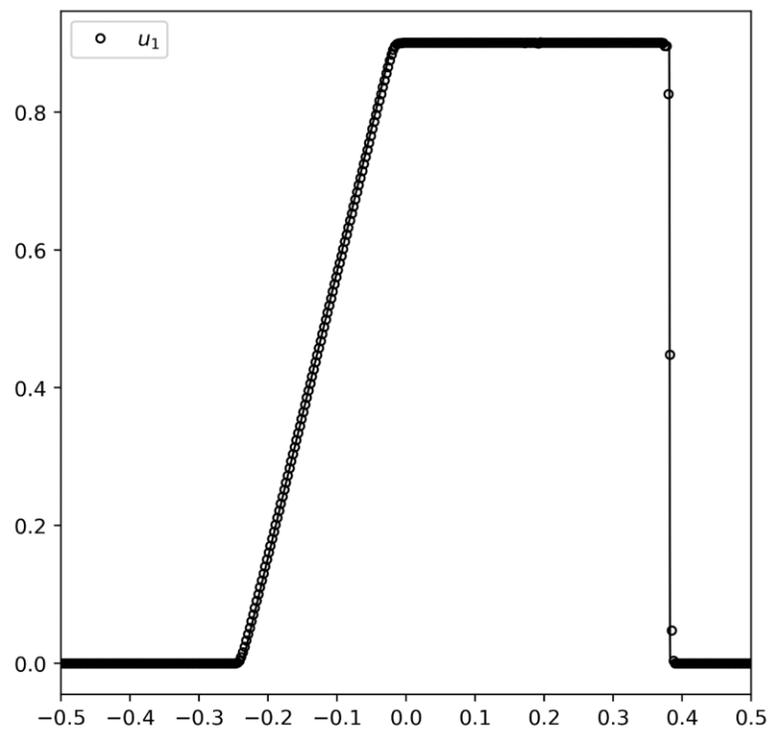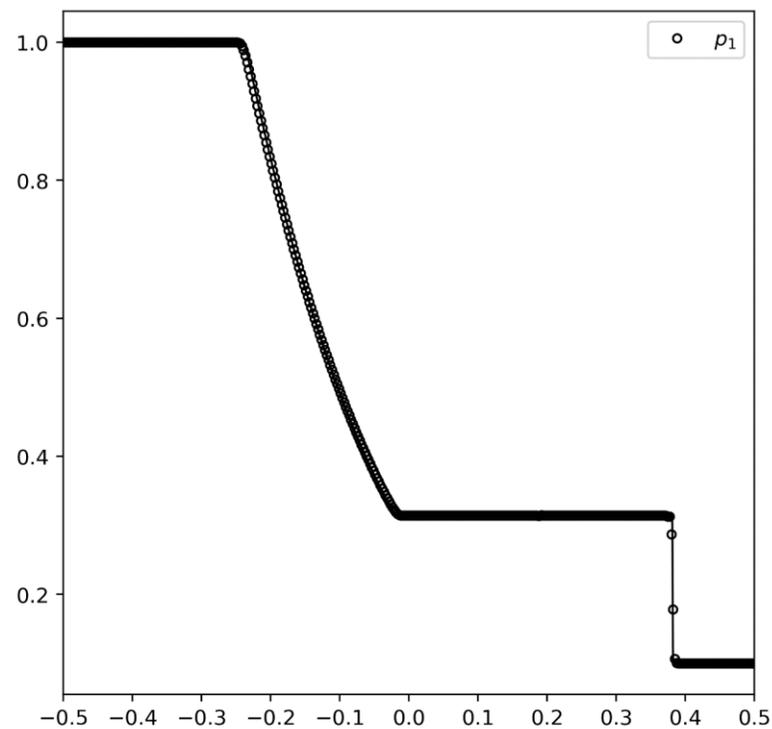

a)          b)          c)

*Fig. 22) Baer-Nunziato with stiff source: Results for the one-dimensional Riemann Problem using the 5$^{th}$ order accurate HLL-based AFD WENO-AO-(5,3) scheme with 400 zones. Figs. 22a, 22b and 22c show the solid density, solid x-velocity and solid pressure.*

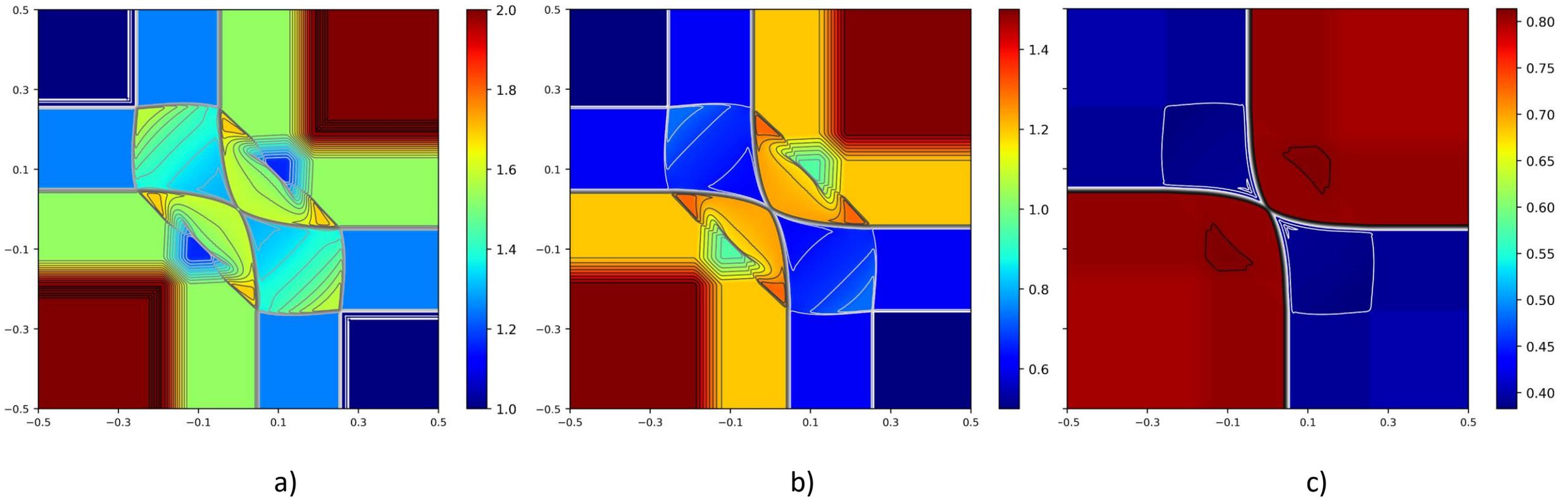

*Fig. 23) Baer-Nunziato with stiff source: Results for the two-dimensional Riemann Problem using the 7$^{th}$ order accurate HLL-based AFD WENO-AO-(7,3) scheme with 400 ×400 zones. Fig. 23a shows the solid density, Fig. 23b shows the gas density and Fig. 23c shows the solid volume fraction. 30 equidistant contour lines are shown over the color plots.*